\documentclass{article}

\usepackage{amsmath, amssymb} 
\usepackage{enumerate}
\usepackage{pgf}
\usepackage{tikz}
\usepackage{amsthm}
\newcommand{\bc}{\begin{C}}
\newcommand{\ec}{\end{C}}
\newcommand{\be}{\begin{equation}}
\newcommand{\ee}{\end{equation}}

\newcommand{\claim}{\begin{Cl}}
\newcommand{\eclaim}{\end{Cl}}
\newcommand{\nb}{\begin{Nb}}
\newcommand{\nbe}{\end{Nb}}
\newcommand{\bl}{\begin{LE}}
\newcommand{\el}{\end{LE}}

\newtheorem{Cl}{Claim}

\newcommand{\bd}{\begin{Def}}
\newcommand{\ed}{\end{Def}}
\newcommand{\bt}{\begin{Th}}
\newcommand{\et}{\end{Th}}
\newcommand{\bp}{\begin{Prop}}
\newcommand{\ep}{\end{Prop}}

\newtheorem{Prop}{Proposition}
\newtheorem{Th}{Theorem}
\newtheorem{LE}{Lemma}
\newtheorem{C}{Corollary}
\newtheorem{Nb}{Note}
\theoremstyle{definition}
\newtheorem{Def}{Definition}
\theoremstyle{definition}
\newtheorem*{remark}{Remark}
\theoremstyle{definition}
\newtheorem*{ex}{Example}

\def\id{\rm id}
\def\ev{\rm ev}

\DeclareMathSymbol{\mlq}{\mathord}{operators}{``}
\DeclareMathSymbol{\mrq}{\mathord}{operators}{`'}

\begin{document}
 \title{Classical linear logic, cobordisms and  categorial grammars
}
\author{Sergey Slavnov
\\  National Research University Higher School of Economics
\\ sslavnov@yandex.ru\\} \maketitle

\begin{abstract}
We propose a categorial grammar based on classical multiplicative linear logic.

This can be seen as an extension of abstract categorial grammars (ACG) and is at least as expressive. However, constituents of {\it linear logic grammars (LLG)} are not abstract ${\lambda}$-terms, but simply tuples of words with labeled endpoints
and supplied with specific {\it plugging instructions}: the sets of endpoints are subdivided into the {\it incoming} and the {\it outgoing} parts. We call such objects {\it word cobordisms}.

A key observation is that word cobordisms  can be organized in a category, very similar to the familiar category of topological cobordisms. This category is symmetric monoidal closed and compact closed and thus is a model of linear $\lambda$-calculus and classical, as well as intuitionistic linear logic. This allows us using linear logic as a typing system for word cobordisms.

 At least, this gives a concrete and intuitive representation of ACG.

We think, however, that the category of word cobordisms, which has a rich  structure and is  independent of any grammar, might be  interesting on its own right.
\end{abstract}

\section{Introduction}
A prototypical example of categorial grammar is Lambek grammars \cite{Lambek}. These are based on  logical {\it Lambek calculus}, which is, speaking in modern terms,  a noncommutative variant of (intuitionistic) linear logic \cite{Girard}. It is well known  that Lambek grammars  generate exactly the same class of languages as context-free grammars \cite{Pentus}.

However, it is agreed that context-free grammar are, in general, not sufficient for modeling natural language. Linguists consider various more expressive formalisms. Lambek calculus is extended to different complex  {\it multimodal}, {\it mixed commutative} and {\it mixed nonassociative} systems, see \cite{Moortgat}.
Many grammars operate with more complex constituents than just words. For example {\it displacement grammars} \cite{Morrill_Displacement},   extending Lambek grammars, operate on discontinuous tuples of words.

 {\it Abstract categorial grammars (ACG)} \cite{deGroote},  on the other hand,  are based on a more intuitive and familiar {\it commutative} logic, namely,  the implicational fragment of  linear logic. There are also other formalisms very close to ACG,  such as  {\it $\lambda$-grammars}  \cite{Muskens} or {\it linear grammars}  \cite{PollardMichalicek} developed in this setting. (We should mention {\it hybrid type logical grammars} \cite{KubotaLevine} as well; these extend ACG, mixing them with Lambek grammars.)

 ACG are very strong in their expressive power  \cite{YoshinakaKanazawa} and show remarkable flexibility in many aspects. (Yet, they have their own deficiencies  from the point of view of natural language modeling as well, see \cite{Moot_inadequacy}). However, in a striking difference from other formalisms,  basic constituents of ACG used for syntax generation are just linear $\lambda$-terms. Unfortunately it is not so easy to identify $\lambda$-terms  with any elements of language.

We should also mention  another,  very interesting approach based on commutative logic. It  is the unifying approach of  \cite{Moot_comparing}. It turns out that many grammatical formalisms can be faithfully  represented as fragments of first order multiplicative intuitionistic linear logic {\bf MILL1}. This provides some common ground on which different systems can be compared.

 In this work we propose one more categorial grammar based on a commutative system, namely on {\it classical} linear logic.

 {\it Linear logic grammars (LLG)} of this paper can be seen as an extension of ACG to the full multiplicative fragment. Although, the list of different formalisms is already sufficiently long, we think that our work deserves some interest at least for two reasons.

 First, unlike the case of ACG, constituents of LLG are very simple. They are tuples of words with labeled endpoints, we call them {\it multiwords}. Multiwords are directly identified as basic elements of language, and apparently they are somewhat easier to deal with than abstract $\lambda$-terms. Each multiword of LLG is
  supplied with specific {\it plugging instructions}: the set of its endpoints is subdivided into the {\it incoming} and the {\it outgoing} parts (not to be confused with the subdivision into left and right endpoints). We call a multiword supplied with such a subdivision a {\it word cobordism}, and we usually abbreviate this title as {\it cowordism}, for a joke.

   Word cobordisms can be composed by gluing outgoing parts (outputs) to the incoming ones (inputs).
    Word endpoints can also be moved from input to output and vice versa producing  new word cobordisms; an operation directly corresponding to $\lambda$-abstraction.

     ACG embed into LLG, so at least we give a concrete and intuitive {\it representation} of ACG. We don't know if LLG have stronger expressive power as ACG, or just the same.
     It can be observed though that the cowordism representation has been  widely used in ACG anyway, but {\it implicitly}, under the title of {\it proof-nets}.

 Second,
  we identify on the class of word cobordisms  a fundamental algebraic structure, namely the structure of a {\it category} (in the mathematical, rather than linguistic sense of the word).

\subsection{Abstract algebra point of view}
   The algebraic structure underlying linguistic interpretations of Lambek calculus is that of a monoid:
 the set of words over a given alphabet is a free monoid under concatenation, and Lambek calculus can be interpreted as a logic of the poset of this monoid subsets (i.e. of formal languages).

When constituents of a grammar are more complicated, such as word tuples, there is no unique concatenation, since tuples can be glued together in many ways. Thus the algebra is more complex.

   When the constituents are word cobordisms, the underlying algebraic  structure is a {\it category} (in the mathematical, rather than linguistic sense of the word).

Word cobordisms, just like ordinary topological cobordisms \cite{Baez}  form a  category, which
is {\it symmetric monoidal closed} and {\it compact closed}.

Thus, we shift from a non-commutative (``nonsymmetric'') monoid  of words  to a symmetric (``commutative'') monoidal category of  {\it cowordisms}.

It is this categorical structure that allows us representing linear $\lambda$-calculus and ACG, as well as classical linear logic.


An  LLG consists of the category of cowordisms, seen as a (degenerate) model of  linear logic, supplied with
 a {\it lexicon}, which is a finite set of non-logical axioms, i.e.  cowordisms together with their typing specifications.

Syntactic derivations from the lexicon directly translate to  {\it cowordisms  generated by the grammar}.


Comparing with Lambek calculus, we shift    from a {\it poset} of formal languages to a {\it category} of cowordism types.
If sequents in Lambek calculus simply represent  inclusions between languages, typing judgements in LLG represent concrete cowordisms  computed from syntactic derivations,  in other words, nontrivial instructions  how to glue (multi-)words together.

It should be observed that the category of cowordisms itself is independent of any grammar.
Apparently, at least some other formalisms can be represented in this setting as well. Possibly,  this can give some common reference for different systems.
In this connection, it might be interesting to study relationship between LLG, as well as the category of cowordisms in general, and the  first order linear logic framework for categorial grammars of \cite{Moot_comparing}. (Apparently, cowordism endpoints correspond to first order variables, and gluing to variable unification.)

\subsection{Syntax vs. semantics}
Somewhat ironically, the word ``semantics'' in the context of categorial grammars can have different semantics .

On one hand, we have {\it denotational} semantics of logical systems. For example, Lambek grammars generate denotational models for Lambek calculus, and LLG of this work are models of linear logic and linear $\lambda$-calculus.

On the other hand, we can have models for  semantics of languages generated by a grammar.

One of the main features making categorial grammars interesting is that they allow a bridge between language   syntax and  language semantics (see \cite{MootRetore}).

In particular, ACG, apart from generating {\it syntactic} representation (which corresponds to  string or tree object signatures), can also generate {\it meaning representation} by changing the object signature, provided that meaning is modeled by means of the same typing system, i.e., intuitionistic linear logic, as in \cite{Dalrymple}.
(We do not discuss here the case of {\it non-linear} ACG.)

As for LLG,  obviously, they do not generate any meaning representation.   Word cobordisms  are just specifically labeled tuples of words. (Speaking in more linguistic terms,  LLG are related directly only to {\it surface} structures.)
However, translating from ACG to LLG does not lead to any information loss: any semantic analysis of a language provided by an ACG remains available.

 Indeed, ACG represent language generation in terms of   $\lambda$-calculus typing derivations. These are mapped  to a chosen meaning representation, thus providing a bridge with semantics.  But ACG for syntactic representation translate to LLG {\it isomorphically}. All typing derivations are preserved and can be mapped to meaning representation equally well. LLG  becomes  then a faithful and  convenient representation of the syntactic (``surface'') part of ACG

Since the category of word cobordisms is compact closed, we would like to mention an  interesting approach to language semantics, where the  meaning representation model is    compact closed as well, moreover, compactness plays a crucial role.

 In {\it categorical compositional distributional models of meaning (DisCoCat)} \cite{CoeckeSadrzadehClark}, \cite{Coecke_Lmabek_vs_Lambek}  it is proposed to model and analyze language semantics by a functorial mapping (``quantization'') of syntactic derivations in a categorial grammar to the (symmetric) compact closed category {\bf FDVec} of finite-dimensional vector spaces. The approach has been developed so far mainly  on the base of Lambek grammars or pregroup grammars (see \cite{Lambek_pregroups}), which are, from the category-theoretical point of view, {\it non-symmetric monoidal closed} (non-symmetric compact closed, in the case of pregroups). On the other hand,  the cowordism category is symmetric and compact closed, and in this sense it is a better mirror  of {\bf FDVec}. Thus it seems a more natural candidate for quantization. Possibly, cowordism representation may help to apply  ideas of DisCoCat to LLG or ACG, thus going beyond context-free languages.

\subsection{Background}
We assume that the reader is familiar with $\lambda$-calculus (see \cite{Barendregt}) and has at least some basic idea of sequent calculus and cut-elimination (see, for example, \cite{TroelstraSchwichtenberg} or \cite{GirardProofsAndTypes} for introduction.)

We also assume some basic acquaintance with categories, in particular, with monoidal categories, see \cite{MacLane} for background. However, we tried to separate concrete examples and definitions from  general categorical discussion, so that the reader uncomfortable with ``abstract nonsense'' might hopefully still grasp what is going on. Categorical formalization is essential in order to guarantee soundness of our constructions
without going through (rather routine) proofs by induction on syntactic derivations.


\section{Word cobordisms}
\subsection{Multiwords}
Let $T$ be a finite alphabet.
We denote the set of all finite words in $T$ as $T^*$, and the empty word as $\epsilon$.

For a set $X$ of natural numbers and $n\in{\bf N}$ we denote
$$n+X=\{n+m|~m\in X\},\quad n-X=\{n-m|~m\in X\}.$$
For a positive integer $n$, we denote ${\bf I}(n)=\{1,\ldots,n\}$.

We are going to define {\it multiwords} over $T$ as oriented bipartite graphs with edges labeled with words in $T^*$.

First we define {\it boundaries}, which represent ordered vertex sets.

\bd
A {\it boundary} $X$ consists of a natural number $|X|$, the {\it cardinality} of $X$, and a subset $X_l\subseteq{\bf I}(|X|)$, called the {\it left boundary} of $X$.
\ed
We will sometimes denote ${\bf I}(|X|)={\bf I}(X)$.

The complement $X_r$ of $X_l$ in ${\bf I}(X)$ is called the {\it right boundary} of $X$.
Elements of $X_l$ are called {\it left endpoints} of $X$ and are said to have {\it left polarity}. Accordingly, elements of $X_r$ are {\it right endpoints} of $X$ and  have {\it right polarity}.

The boundary $X$ can be thought of as the set ${\bf I}(X)$ equipped with the partition into left and right endpoints. It is implied that $X$ parameterizes vertices of a directed bipartite graph whose edges are directed from left to right.

Before turning to graphs, we define some operations on boundaries.

Let $X,Y$ be boundaries and $i\in{\bf N}$ be such that $i+|Y|\leq|X|$.
We say that
$i+Y$ is a {\it subboundary} of $X$ if $i+Y_l\subseteq X_l$ and $i+Y_r\subseteq X_r$.

Given two boundaries $X$ and $Y$, we define the {\it tensor product} $X\otimes Y$ as the boundary with
$$|{X\otimes Y}|=|X|+|Y|,\quad (X\otimes Y)_l=X_l\cup(|X|+Y_l).$$

Given a boundary $X$, we define the {\it dual boundary} $X^\bot$ as the boundary with
 $$|{X^\bot}|=|X|,\quad (X^\bot)_l=|X|+1-X_r.$$

Thus,  tensor product concatenates two boundaries, while  duality reverses the vertex ordering and polarity

Note that we have
\be\label{tensor dual}
(X\otimes Y)^\bot=Y^\bot\otimes X^\bot.
\ee
\begin{remark}
A reader familiar with compact categories (which will be discussed shortly) or cyclic linear logic might anticipate from formula (\ref{tensor dual}) some sort of noncommutativity in the category of cowordisms. This is misleading. We will have a natural isomorphism between $(A\otimes B)^\bot$ and $A^\bot\otimes B^\bot$, just not equality.
\end{remark}

\bd
A {\it regular multiword} $M$ with boundary $X$ over an alphabet $T$ is a  directed   graph on the set ${\bf I}(X)$ of vertices, whose  edges are labelled with words in $T^*$,   such that each vertex is adjacent
to exactly one edge (so that it is a perfect matching), and for every edge its left endpoint is in $X_l$ and its right endpoint is in $X_r$.
\ed
We use notation $(x,w,y)$ for an edge from $x$ to $y$ labeled with $w$. We will identify a regular multiword with the set of its labeled edges.

A basic operation on multiwords will be gluing. Since, in principle, such gluing can be done cyclically, it can result in objects that are not  graphs at all.
For consistency of definitions we will also have to consider {\it cyclic words} and {\it singular multiwords}.

We say that  two words in $T^*$ are {\it cyclically equivalent} if they differ by a cyclic permutation of letters.
A {\it cyclic word} over  $T$ is an equivalence class of cyclically equivalent words in $T^*$.
For $w\in T^*$  we denote the corresponding cyclic word  as $[w]$.

\bd
A {\it multiword} $M$ with boundary $X$ over the alphabet $T$ is a pair $M=(M_0,M_c)$, where $M_0$, the {\it regular} part, is a regular multiword over $T$ with the boundary $X$, and $M_c$, the {\it singular} or {\it cyclic} part, is a finite multiset of cyclic words over $T$.
\ed
In the following we denote the sum (i.e. disjoint union) of multisets as $+$.

A multiword is {\it acyclic} or {\it regular} if its singular part is empty. Otherwise it is {\it singular}.

A general multiword $M$ can be pictured geometrically as the edge-labelled graph $M_0$ and a disjoint union of isolated loops (i.e. closed curves) labeled with elements of $M_c$.
The underlying geometric object is no longer  a graph, but it is a topological space,  even a topological manifold with boundary.


Tensor product readily extends from boundaries to multiwords.

 Given two multiwords $M=(M_0,M_c)$ and $N=(N_0,N_c)$  with boundaries $X$ and $Y$ respectively, the {\it tensor product} $M\otimes N$ is the multiword with the boundary $X\otimes Y$ defined by
  $$(M\otimes N)_c=M_c+ N_c$$
 $$(M\otimes N)_0=\{(i,w,j)|\mbox{ } (i,w,j)\in M_0\}\cup\{(|X|+i,w,|X|+j)|\mbox{ } (i,w,j)\in N_0\}.$$

\subsubsection{Contractions}
Let $M$ be a multiword with  boundary $X$, and let $n<|X|$   be such that $n$ and $n+1$ have opposite polarity.

We define the {\it elementary contraction} $\langle M\rangle_{n,n+1}$ of $M$ as a multiword obtained by gluing edges of $M$ along  $n$ and $n+1$.
An accurate definition follows.

Define an embedding $$\phi:{\bf I}( |X|-2)\to {\bf I}(X)$$ by
$$
\phi(i)=
\begin{cases}
i,\mbox{ if }i<n,\\
i+2\mbox{ otherwise.}
\end{cases}
$$
Define the boundary $X'$
by
$$\rm{card}_{X'}=|X|-2,\quad (X')_l=\phi^{-1}(X_l).$$

Now let $x$ be the element of the pair $(n,n+1)$ that belongs to $X_r$, and $y\in X_l$   be the other element.
We define a new  multiword $ M'$ with the boundary $X'$  by cases.
\begin{itemize}
\item
If $x$ and $y$ are not connected by an edge in $M_0$, then
$$M_c'=M_c,$$
$$M_0'=\{(i,w,j)|\mbox{ }(\phi(i),w,\phi(j)\in M_0)\}\cup\{(i,uv,j)|\mbox{ }(\phi(i),u,x),(y,v,\phi(j))\in M_0\}.$$
\item If there is an edge $(y,w,x)\in M_0$, then
$$M_c'=M_c+\{[w]\},$$
$$M_0'=\{(i,w,j)|\mbox{ }(\phi(i),w,\phi(j)\in M_0)\}.$$
\end{itemize}
It is easy to see that in both cases $M_0'$ is a perfect matching and its edges start at left endpoints of $X_l'$.
We put $\langle M\rangle_{n,n+1}=M'$.
\smallskip

Elementary contractions can be iterated.

Let $X$, $Y$ be  boundaries, $i\in {\bf N}$ and assume that $i+Y^\bot\otimes Y$ is a subboundary of $X$.
Let $n=|Y|=\rm{card}_{Y^\bot}$

Then for any multiword $M$ with the boundary $X$  we define {\it iterated contraction} $\langle M\rangle_{i+Y^\bot\otimes Y}$ of $M$ by
$$\langle M\rangle_{i+Y^\bot\otimes Y}=\langle\ldots\langle~ \langle~\langle M\rangle_{i+n,i+n+1}\rangle_{i+n-1,i+n}\rangle\ldots\rangle_{i+1,i+2}.$$
(It is easy to check that the above is well defined.)

\subsubsection{Word cobordisms}
We want to consider multiwords as morphisms between boundaries and compose them by gluing. So we want to organize multiwords in a category.

We  remarked above that multiwords can be represented geometrically as very simple manifolds with boundary. Manifolds with boundary give rise to the category of {\it cobordisms}, see \cite{Baez}. We will construct a similar category  of {\it word cobordisms}. We find it amusing to abbreviate the latter term as {\it cowordism}, and we will do so.

In order to  treat a multiword as a morphism consistently, we need to subdivide its boundary into the {\it incoming} and the {\it outgoing} parts and to fix this subdivision. (It should not be confused with the subdivision into left and right endpoints.) A multiword equipped with such a subdivision
will be called a cowordism.
\bd
Given two boundaries $X,Y$, a {\it cowordism} $$\sigma:X\to Y$$
 over an alphabet $T$ from $X$ to $Y$ is  a multiword over $T$ with  boundary $Y\otimes X^\bot$.
\ed
In the above setting we say that $Y$ is the {\it outgoing boundary} of $\sigma$, and $X$ is the {\it incoming   boundary}. The outgoing boundary can be seen as a subboundary of $\sigma$ in an obvious way. The incoming boundary $X$ can be identified with the subboundary $|Y|+X^\bot$ by means of an order and polarity reversing bijection.

We say that a cowordism is {\it regular} if its underlying multiword is regular. Otherwise the cowordism is {\it singular}.

Matching cowordisms are composed  by gluing incoming and outgoing boundaries.

Given boundaries $X,Y,Z$ and cowordisms
$$\sigma:X\to Y,\quad \tau:Y\to Z,$$
the {\it composition} $$\tau\circ\sigma:X\to Z$$
of $\tau$ and $\sigma$
is the cowordism defined by the iterated contraction
$$\tau\circ\sigma=\langle\tau\otimes\sigma\rangle_{|Z|+Y^\bot\otimes Y}.$$
\smallskip

In order to discuss cowordisms in a greater detail it will be convenient to develop some systematic conventions for representing them geometrically. Indeed, cowordisms are, by definition, geometric objects, and a graphical language makes their properties very transparent.

\section{Graphical language}
When  depicting a cowordism
\be\label{sample cowordism}
\sigma:X\to Y,
\ee
it is convenient to put the boundary vertices
on two parallel straight lines, reserving one line for the outgoing subboundary $Y$ and the other one for the incoming subboundary $\rm{card_Y}+X^\bot$, and to place all edges and loops (if there is a singular part) that comprise $\sigma$ between the two lines.

We will use either the {\it horizontal} or the {\it vertical} representation, the horizontal one being default.
In both cases we will depict left endpoints of the boundary of $\sigma$ as small filled-in circles and right endpoints, as arrow heads.

In the horizontal representation
of  cowordism (\ref{sample cowordism}),
we put elements $$1,\ldots,|Y|$$ of the outgoing subboundary $Y$ on one vertical line, with the  increasing order  corresponding to the direction {\it up}, and
we put the elements $$\rm{card_Y}+1,\ldots,|Y|+|X|$$ of the incoming subboundary $|Y|+X^\bot$ on a parallel line to the left, in the increasing order corresponding to the direction {\it down}.

For example, if the boundaries $X,Y$ are given by
\be\label{graph.lang.example}
|X|=4,\quad X_l=\{2\},\quad |Y|=4,\quad Y_l=\{3\},
\ee
then  cowordism (\ref{sample cowordism}) will be depicted as the following picture (we put vertex numbers for clarity).
$$
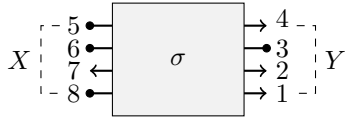
\begin{tikzpicture}
\begin{scope}[shift={(-.7,0)}]
\draw[dashed](0.25,1.2)--(0,1.2)--(0,.3)--(0.25,.3);
\node[left]at(0,.75){$X$};
\end{scope}
\draw[thick](-.05,0.3)--(0.25,0.3);
\draw[thick,<-](-.05,0.6)--(0.25,0.6);
\draw[thick,-](-.05,.9)--(0.25,.9);
\draw[thick](-.05,1.2)--(0.25,1.2);

\node[left] at(-.05,0.3){$8$};
\node[left] at(-.05,0.6){$7$};
\node[left] at(-.05,0.9){$6$};
\node[left] at(-.05,1.2){$5$};

\node[right] at(2.3,0.3){$1$};
\node[right] at(2.3,0.6){$2$};
\node[right] at(2.3,0.9){$3$};
\node[right] at(2.3,1.3){$4$};

\draw [fill] (2.3,.9) circle [radius=0.05];
\draw [fill] (-.05,0.3) circle [radius=0.05];
\draw [fill] (-.050,.9) circle [radius=0.05];
\draw [fill] (-.05,1.2) circle [radius=0.05];

\draw[draw=black,fill=gray!10](0.25,0)rectangle(2,1.5);\node at(1.125,.75){$\sigma$};
\draw[thick,->](2,0.3)--(2.3,0.3);
\draw[thick,->](2,.6)--(2.3,0.6);
\draw[thick,-](2,.9)--(2.3,.9);
\draw[thick,->](2,1.2)--(2.3,1.2);
\begin{scope}[shift={(2.7,0)}]
\draw[dashed](0.,1.2)--(0.25,1.2)--(0.25,.3)--(0.,.3);
\node[right]at(0.25,.75){$Y$};
\end{scope}
\end{tikzpicture}$$
The vertices $1,2,3,4$ in the boundary of $\sigma$ correspond to the outgoing boundary $Y$, and $3$ is the only left endpoint. The vertices $5,6,7,8$ can be identified with the incoming boundary $X$ by means of an order and polarity reversing bijection, and the only left endpoint $2$ of $X$ corresponds to the right endpoint $7$ of $\sigma$.

In the vertical representation we put the vertices of the outgoing subboundary $|Y|$ on one horizontal line, with the increasing order  corresponding to the direction {\it from right to left}, and we put the elements of the incoming boundary $|Y|+{\bf I}(X^\bot)$ on a parallel line below, with the increasing order corresponding to the direction {\it from left to right}.

Thus,  if $X$, $Y$ are as in (\ref{graph.lang.example}),
then  cowordism (\ref{sample cowordism}) in the vertical representation will be depicted as
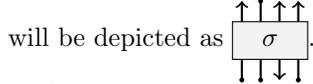
\begin{tikzpicture}[baseline=5pt,scale=.5]
  \draw[thick](-.25,-.55)--(-.25,0);
\draw[thick,<-](-.75,-.55)--(-0.75,0);
\draw[thick,-](-1.25,-.55)--(-1.25,-0);
\draw[thick](-1.75,-.55)--(-1.75,0);

%

\draw [fill] (-1.25,1.55) circle [radius=0.05];
\draw [fill] (-.25,-.55) circle [radius=0.05];
\draw [fill] (-1.25,-.55) circle [radius=0.05];
\draw [fill] (-1.75,-.55) circle [radius=0.05];

\draw[draw=black,fill=gray!10](-2,0)rectangle(0,1);\node at(-1,.5){$\sigma$};
 \draw[thick,<-](-.25,1.55)--(-.25,1.0);
\draw[thick,<-](-.75,1.55)--(-0.75,1.0);
\draw[thick,-](-1.25,1.55)--(-1.25,1);
\draw[thick,<-](-1.75,1.55)--(-1.75,1);
\end{tikzpicture}.

In general, when the detailed structure of the boundary is not important, we ``squeeze'' parallel edges into one and represent (\ref{sample cowordism}) schematically as a box with an incoming wire labeled with $X$ and an outgoing wire labeled with $Y$.
\begin{tikzpicture}[baseline=2.5pt,scale=.7]
\node[left] at(-1,.25){$\sigma:$};
\draw[thick](0,.25)--(-0.5,.25);

\node[left] at(-0.5,.25){$X$};

\draw[thick](1,.25)--(1.5,.25);
\node[right] at(1.5,0.25){$Y$};

\draw[draw=black,fill=gray!10](0.,0)rectangle(1,.5);\node at(.5,.25){$\sigma$};
\end{tikzpicture}

\subsection{Composition}
It is easy to see that, with our conventions, the composition of cowordisms $$\sigma:X\to Y,\quad\tau:Y\to Z$$ corresponds to the following schematic picture.
$$\begin{tikzpicture}[scale=.7]
\node[left] at(-1,.25){$\sigma:$};
\draw[thick](0,.25)--(-0.5,.25);

\node[left] at(-0.5,.25){$X$};

\draw[thick](1,.25)--(1.5,.25);
\node[right] at(1.5,0.25){$Y$};

\draw[draw=black,fill=gray!10](0.,0)rectangle(1,.5);\node at(.5,.25){$\sigma$};

\begin{scope}[shift={(4.5,0)}]
\node[left] at(-1,.25){$\tau:$};
\draw[thick](0,.25)--(-0.5,.25);

\node[left] at(-0.5,.25){$Y$};

\draw[thick](1,.25)--(1.5,.25);
\node[right] at(1.5,0.25){$Z$};

\draw[draw=black,fill=gray!10](0.,0)rectangle(1,.5);\node at(.5,.25){$\tau$};
\end{scope}

\begin{scope}[shift={(1.5,0)}]
\begin{scope}[shift={(8.5,0)}]
\node[left] at(-1,.25){$\tau\circ\sigma:$};
\draw[thick](0,.25)--(-0.5,.25);

\node[left] at(-0.5,.25){$X$};

\draw[thick](1,.25)--(1.5,.25);

\draw[draw=black,fill=gray!10](0.,0)rectangle(1,.5);\node at(.5,.25){$\sigma$};

\begin{scope}[shift={(1.5,0)}]

\draw[thick](1,.25)--(1.5,.25);
\node[right] at(1.5,0.25){$Z$};

\draw[draw=black,fill=gray!10](0.,0)rectangle(1,.5);\node at(.5,.25){$\tau$};
\end{scope}
\end{scope}

\end{scope}
\end{tikzpicture}
$$
We get a detailed, ``full'' picture by expanding each edge into as many  parallel edges as there are points in the corresponding boundary.

For example, if $X,Y$ are as in (\ref{graph.lang.example}), then the above picture translates to the following.
$$
\begin{tikzpicture}[scale=.5]
\node[left]at(-.7,.75){$\sigma:$};
\draw[dashed](0.25,1.2)--(0,1.2)--(0,.3)--(0.25,.3);
\node[left]at(0,.75){$X$};
\begin{scope}[shift={(.8,0)}]
\draw[thick](-.3,0.3)--(0.25,0.3);
\draw[thick,<-](-.3,0.6)--(0.25,0.6);
\draw[thick,-](-.3,.9)--(0.25,.9);
\draw[thick](-.3,1.2)--(0.25,1.2);

%

\draw [fill] (2.55,.9) circle [radius=0.05];
\draw [fill] (-.3,0.3) circle [radius=0.05];
\draw [fill] (-.30,.9) circle [radius=0.05];
\draw [fill] (-.3,1.2) circle [radius=0.05];

\draw[draw=black,fill=gray!10](0.25,0)rectangle(2,1.5);\node at(1.125,.75){$\sigma$};
\draw[thick,->](2,0.3)--(2.55,0.3);
\draw[thick,->](2,.6)--(2.55,0.6);
\draw[thick,-](2,.9)--(2.55,.9);
\draw[thick,->](2,1.2)--(2.55,1.2);
\end{scope}

\begin{scope}[shift={(3.5,0)}]
\draw[dashed](0.,1.2)--(0.25,1.2)--(0.25,.3)--(0.,.3);
\node[right]at(0.25,.75){$Y$};
\end{scope}

\begin{scope}[shift={(7,0)}]
\node[left]at(-.7,.75){$\tau:$};
\draw[dashed](0.25,1.2)--(0,1.2)--(0,.3)--(0.25,.3);
\node[left]at(0,.75){$Y$};
\begin{scope}[shift={(.8,0)}]
\draw[thick](-.3,0.3)--(0.25,0.3);
\draw[thick,-](-.3,0.6)--(0.25,0.6);
\draw[thick,<-](-.3,.9)--(0.25,.9);
\draw[thick](-.3,1.2)--(0.25,1.2);


\node[right] at(2.3,.75){$Z$};

\draw [fill] (-.3,0.3) circle [radius=0.05];
\draw [fill] (-.3,.6) circle [radius=0.05];
\draw [fill] (-.3,1.2) circle [radius=0.05];

\draw[draw=black,fill=gray!10](0.25,0)rectangle(2,1.5);\node at(1.125,.75){$\tau$};
\draw[thick,-](2,0.75)--(2.55,0.75);
\end{scope}

\end{scope}
\begin{scope}[shift={(15,0)}]
\node[left]at(-.8,.75){$\tau\circ\sigma:$};
\draw[dashed](0.25,1.2)--(0,1.2)--(0,.3)--(0.25,.3);
\node[left]at(0,.75){$X$};
\begin{scope}[shift={(.8,0)}]
\draw[thick](-.3,0.3)--(0.25,0.3);
\draw[thick,<-](-.3,0.6)--(0.25,0.6);
\draw[thick,-](-.3,.9)--(0.25,.9);
\draw[thick](-.3,1.2)--(0.25,1.2);

%

\draw [fill] (-.3,0.3) circle [radius=0.05];
\draw [fill] (-.3,.9) circle [radius=0.05];
\draw [fill] (-.3,1.2) circle [radius=0.05];

\draw[draw=black,fill=gray!10](0.25,0)rectangle(2,1.5);\node at(1.125,.75){$\sigma$};
\draw[thick,-](2,0.3)--(3.5,0.3);
\draw[thick,-](2,.6)--(3.5,0.6);
\draw[thick,-](2,.9)--(3.5,.9);
\draw[thick,-](2,1.2)--(3.5,1.2);
\end{scope}

\begin{scope}[shift={(3,0)}]

\begin{scope}[shift={(.8,0)}]
\draw[thick](-.05,0.3)--(0.25,0.3);
\draw[thick,-](-.05,0.6)--(0.25,0.6);
\draw[thick,-](-.05,.9)--(0.25,.9);
\draw[thick](-.05,1.2)--(0.25,1.2);

\node[right] at(2.3,.75){$Z$};

\draw[draw=black,fill=gray!10](0.25,0)rectangle(2,1.5);\node at(1.125,.75){$\tau$};
\draw[thick,-](2,0.75)--(2.5,0.75);
\end{scope}

\end{scope}

\end{scope}
\end{tikzpicture}
$$
%
%
%
%
%
%
%
%
%
%
%
%
\bp
Composition of cowordisms is associative.
\ep
\begin{proof}
  Evident from geometric representation.
\end{proof}

\subsubsection{Identities}
In order to have a category of cowordisms, we also need identities

Given a boundary $X$, the {\it identity cowordism} $$\id_X:X\to X$$
is  the regular multiword with the boundary $X\otimes X^\bot$ defined by the set of edges
$$\{(i,\epsilon,2|X|-i+1)|\mbox{ }i\in X_l\}\cup\{(2|X|-i+1,\epsilon,i)|\mbox{ }i\in X_r\}.$$

In a schematic, ``squeezed'' picture, the identity cowordism corresponds to a single wire
\begin{tikzpicture}[baseline=8pt,scale=.5]
\node[left]at(0.,.75){$\id_X:$};
\node[left]at(0.75,.75){$X$};
\draw[thick](.75,0.75)--(1.75,0.725);
\node[right]at(1.75,.75){$X$};
\end{tikzpicture}.

In the full picture there are as may parallel wires as there are points in $X$.

If $X$ is as in (\ref{graph.lang.example}), then the  full picture is
\begin{tikzpicture}[baseline=8pt,scale=.5]
\begin{scope}[shift={(-7.55,0)}]
\node[left]at(-.7,.75){$\id_X:$};
\draw[dashed](0.25,1.2)--(0,1.2)--(0,.3)--(0.25,.3);
\node[left]at(0,.75){$X$};
\begin{scope}[shift={(.4,0)}]
\draw[thick,->](-.05,0.3)--(1.45,0.3);
\draw[thick,<-](-.05,0.6)--(1.45,0.6);
\draw[thick,->](-.05,.9)--(1.45,.9);
\draw[thick,->](-.05,1.2)--(1.45,1.2);

\draw [fill] (1.45,.6) circle [radius=0.05];
\draw [fill] (-.05,0.3) circle [radius=0.05];
\draw [fill] (-.050,.9) circle [radius=0.05];
\draw [fill] (-.05,1.2) circle [radius=0.05];
\draw[dashed](1.45,1.2)--(1.7,1.2)--(1.7,.3)--(1.45,.3);
\node[right]at(1.7,.75){$X$};
\end{scope}
\end{scope}
\end{tikzpicture}.
\bp
For all boundaries $Y$ and cowordisms $\sigma:X\to Y$ and $\tau:Y\to X$ we have $$\sigma\circ\id_X=\sigma,\quad\id_X\circ\tau=\tau.$$
\ep
\begin{proof}
  Evident from geometric representation.
\end{proof}

We say that a cowordism $\sigma:X\to Y$ is {\it invertible} if there exists a cowordism $\tau:Y\to X$ such that
$$\tau\circ\sigma=\id_X,\quad\sigma\circ\tau=\id_Y.$$

\subsection{Tensor product}
We have already defined  tensor product on boundaries (and multiwords). Now we are going to establish
conventions for depicting tensor product in the graphical language and then to extend it consistently to cowordisms.

Given boundaries $$X_1,\ldots,X_n,Y_1,\ldots Y_m,$$ we represent a cowordism
\be\label{sample tensor}
\sigma:X_1\otimes\ldots\otimes X_n\to Y_1\otimes\ldots\otimes Y_m
\ee
schematically as a box whose $n$ incoming wires are labeled with $X_i$'s and $m$ outgoing wires, labeled with $Y_i$'s.
$$
\begin{tikzpicture}[scale=.8]
\node[left] at(-1.25,-0.){$\sigma:$};
\draw[thick](0,.25)--(-0.5,.25);

\node[left] at(-0.5,.25){$X_n$};

\draw[thick](1,.25)--(1.5,.25);
\node[right] at(1.5,0.25){$Y_m$};

\node[right] at(1,0){$\cdots$};
\node[left] at(0,0){$\cdots$};

\draw[draw=black,fill=gray!10](0.,-.5)rectangle(1,.5);\node at(.5,0){$\tau$};

\begin{scope}[shift={(0,-.5)}]

\draw[thick](0,.25)--(-0.5,.25);

\node[left] at(-0.5,.25){$X_1$};

\draw[thick](1,.25)--(1.5,.25);
\node[right] at(1.5,0.25){$Y_1$};

\end{scope}
\end{tikzpicture}
$$
Note that the above ``squeezed'' picture is consistent with the full picture. If we ``expand'' each edge into  parallel edges adjacent to  points in the corresponding subboundary, we obtain a detailed picture of (\ref{sample tensor}) according to our conventions.

The case when $n$ or $m$ in the above picture is $0$ also makes sense.
\bd
The {\it unit boundary} ${\bf 1}$ is defined by $|{\bf 1}|=0$, ${\bf 1}_l=\emptyset$.
\ed
\nb
For any boundary $X$ we have $X\otimes {\bf 1}={\bf 1}\otimes X=X$. $\Box$
\nbe

When depicting a cowordism $\sigma:{\bf 1}\to X$, respectively,  $\tau:X\to{\bf 1}$ we do not have wires on the left, respectively, right.

Now
let boundaries $X,Y,Z,T$ and cowordisms
$$\sigma:X\to Y,\quad \tau:Z\to T$$ be given.

Let us write $\sigma_0$, respectively, $\tau_0$ for the regular part of (the underlying multiword of) $\sigma$, respectively $\tau$, and let us write $\sigma_c$, $\tau_c$ for the respective singular parts.

Let $S=(Y\otimes T)\otimes(X\otimes Z)^\bot=Y\otimes T\otimes Z^\bot\otimes X^\bot$.

Define the embeddings  $$\phi:{\bf I}(|Y|+|X|)\to{\bf I}(|Y|+|T|+|Z|+|X|),\quad
\psi:{\bf I}(|T|+|Z|)\to{\bf I}(|Y|+|T|+|Z|+|X|)$$
 by
$$\phi(i)=
\begin{cases}
i,\mbox{ if }i\leq |Y|,\\
i+|T|+|Z|\mbox{ otherwise},
\end{cases}
\quad
 \psi(i)=i+|Y|.
$$

The {\it tensor product} $\sigma\otimes \tau$ is the multiword with the boundary $S$ defined by
  $$(\sigma\otimes \tau)_c=\sigma_c+ \tau_c,$$
 $$(\sigma\otimes \tau)_0=\{(\phi(i),w,\phi(j)|\mbox{ }(i,w,j)\in \sigma_0\}\cup
 \{(\psi(i),w,\psi(j)|\mbox{ }(i,w,j)\in \tau_0\}.$$

In the graphical language, tensor product of cowordisms corresponds simply to putting two boxes side by side.
$$\begin{tikzpicture}[scale=.7]
\node[left] at(-1,-0.125){$\sigma\otimes\tau:$};
\draw[thick](0,.25)--(-0.5,.25);

\node[left] at(-0.5,.25){$Z$};

\draw[thick](1,.25)--(1.5,.25);
\node[right] at(1.5,0.25){$T$};

\draw[draw=black,fill=gray!10](0.,0)rectangle(1,.5);\node at(.5,.25){$\tau$};

\begin{scope}[shift={(0,-.75)}]

\draw[thick](0,.25)--(-0.5,.25);

\node[left] at(-0.5,.25){$X$};

\draw[thick](1,.25)--(1.5,.25);
\node[right] at(1.5,0.25){$Y$};

\draw[draw=black,fill=gray!10](0.,0)rectangle(1,.5);\node at(.5,.25){$\sigma$};
\end{scope}
\end{tikzpicture}
$$
\bp
Tensor product of cowordisms is functorial:
$$(\tau_1\otimes\tau_2)\circ(\sigma_1\otimes\sigma_2)=(\tau_1\circ\sigma_1)\otimes(\tau_2\circ\sigma_2)$$
for all boundaries $X_i,Y_i,Z_i$ and cowordisms
$$\sigma_i:X_i\to Y_i,\quad \tau_i:Y_i\to Z_i,\quad i=1,2.$$
\ep
\begin{proof}
Evident from geometric representation.
\end{proof}

\subsubsection{Symmetries}
Let the boundaries $X,Y$ be given.

The {\it symmetry cowordism}
$$s_{X,Y}:X\otimes Y\to Y\otimes X$$
is the regular multiword with the boundary $Y\otimes X\otimes Y^\bot\otimes X^\bot$ defined by the set of edges
$$\{(|Y|+i,\epsilon,2|X|+2|Y|-i+1)|\mbox{ }i\in X_l\}\cup\{(2|X|+2|Y|-i+1,\epsilon,|Y|+i)|\mbox{ }i\in X_r\}\cup$$
$$\{(i,\epsilon,|X|+2|Y|-i+1)|\mbox{ }i\in Y_l\}\cup\{(|X|+2|Y|-i+1,\epsilon,i)|\mbox{ }i\in Y_r\}.$$

The schematic  picture of $s_{X,Y}$ is transparent.
$$
\begin{tikzpicture}[scale=.7]
\node[left] at(-1.25,0.5){$s_{XY}:$};
\draw[thick](-0.5,.25)--(0.,.25)--(0.5,.75)--(1,.75);
\draw[thick](-0.5,.75)--(0.,.75)--(0.5,.25)--(1,.25);

\node[left] at(-0.5,.25){$X$};
\node[left] at(-0.5,.75){$Y$};

\node[right] at(1,0.75){$X$};
\node[right] at(1,0.25){$Y$};
\end{tikzpicture}
$$
\bp\label{symmetry squared}
For all boudaries $X,Y$, the symmetry cowordism is invertible, with $s_{Y,X}\circ s_{X,Y}=\id _{X\otimes Y}$.
\ep
\begin{proof}
Obvious.
\end{proof}

\subsection{Duality}
We have already defined the operation of duality on boundaries. Now we extend it to cowordisms.

Let $X$, $Y$ be boundaries, and let $\sigma:X\to Y$ be a cowordism.

Let $\phi:{\bf I}(|Y|+|X|)\to{\bf I}(|X|+|Y|)$ be the cyclic permutation
$$\phi(i)=
\begin{cases}
i+|X|,\mbox{ if }i\leq |Y|,\\
i-|Y|\mbox{ otherwise}.
\end{cases}
$$


The {\it dual} cowordism $$\sigma^\bot:Y^\bot\to X^\bot$$
is the multiword with  the boundary
 $X^\bot\otimes Y^{\bot\bot}=X^\bot\otimes Y$
 defined by
 $$(\sigma^\bot)_c=\sigma_c,$$
 $$(\sigma^\bot)_0=\{(\phi(i),w,\phi(j))|\mbox{ }(i,w,j)\in \sigma_0\}.$$

In a schematic  picture, duality looks as follows.
$$
\begin{tikzpicture}[scale=.5]
\begin{scope}[shift={(-7,0)}]
\node[left]at(-.7,.25){$\sigma:$};
\node[left]at(0,.25){$X$};
\begin{scope}[shift={(.2,0)}]
\draw[thick](-.25,0.25)--(0.25,0.25);

\node[right] at(2.5,.25){$Y$};

\draw[draw=black,fill=gray!10](0.25,0)rectangle(2,.5);\node at(1.125,.25){$\sigma$};
\draw[thick,-](2,0.25)--(2.5,0.25);
\end{scope}

\end{scope}


\node[left]at(-.45,.25){$\sigma^\bot:$};

\draw[thick](0.25,0.25)
to[out=180,in=-90](-.25,.5)
to[out=90,in=180](.25,1)--(2.5,1)
;
\node[right] at(2.5,1){$X^\bot$};

\draw[draw=black,fill=gray!10](0.25,0)rectangle(2,.5);\node at(1.125,.25){$\sigma$};
\draw[thick,-](2,0.25)
to[out=0,in=90](2.5,0)
to[out=-90,in=0](2,-0.5)--(-.25,-0.5)
;

\node[left] at(-.25,-0.5){$Y^\bot$};

\end{tikzpicture}
$$
The full picture, again, can be recovered by expanding every wire  into a parallel cluster.

For example, if $X,Y$ are as in (\ref{graph.lang.example}), the above picture translates to the following.
$$
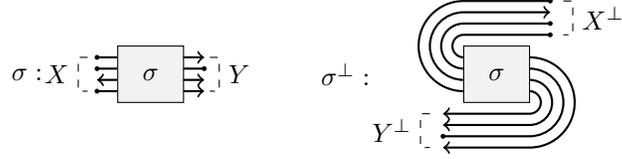
\begin{tikzpicture}[scale=.5]
\begin{scope}[shift={(-10,0)}]
\node[left]at(-.7,.75){$\sigma:$};
\draw[dashed](0.25,1.2)--(0,1.2)--(0,.3)--(0.25,.3);
\node[left]at(0,.75){$X$};
\begin{scope}[shift={(.8,0)}]
\draw[thick](-.3,0.3)--(0.25,0.3);
\draw[thick,<-](-.3,0.6)--(0.25,0.6);
\draw[thick,-](-.3,.9)--(0.25,.9);
\draw[thick](-.3,1.2)--(0.25,1.2);

\draw [fill] (2.55,.9) circle [radius=0.05];
\draw [fill] (-.3,0.3) circle [radius=0.05];
\draw [fill] (-.3,.9) circle [radius=0.05];
\draw [fill] (-.3,1.2) circle [radius=0.05];

\draw[draw=black,fill=gray!10](0.25,0)rectangle(2,1.5);\node at(1.125,.75){$\sigma$};
\draw[thick,->](2,0.3)--(2.55,0.3);
\draw[thick,->](2,.6)--(2.55,0.6);
\draw[thick,-](2,.9)--(2.55,.9);
\draw[thick,->](2,1.2)--(2.55,1.2);
\end{scope}

\begin{scope}[shift={(3.5,0)}]
\draw[dashed](0.,1.2)--(0.25,1.2)--(0.25,.3)--(0.,.3);
\node[right]at(0.25,.75){$Y$};
\end{scope}
\end{scope}

\begin{scope}[shift={(-0.9,-1.5)}]
\draw[dashed](0.25,1.2)--(0,1.2)--(0,.3)--(0.25,.3);
\node[left]at(0,.75){$Y^\bot$};
\end{scope}

\begin{scope}[shift={(2.9,1.5)}]
\draw[dashed](0,1.2)--(0.25,1.2)--(0.25,.3)--(0,.3);
\node[right]at(0.25,.75){$X^\bot$};
\end{scope}

\node[left]at(-2,.75){$\sigma^\bot:$};

\draw[thick](.25,0.3)
to[out=180,in=-90](-.95,1.5)
to[out=90,in=180](.25,2.7)--(2.55,2.7)
;
\draw[thick,->](0.25,0.6)
to[out=180,in=-90](-.65,1.5)
to[out=90,in=180](.25,2.4)--(2.55,2.4)
;
\draw[thick,-](0.25,.9)
to[out=180,in=-90](-.35,1.5)
to[out=90,in=180](.25,2.1)--(2.55,2.1)
;
\draw[thick](0.25,1.2)
to[out=180,in=-90](-.05,1.5)
to[out=90,in=180](.25,1.8)--(2.55,1.8)
;

\draw [fill] (2.55,2.7) circle [radius=0.05];
\draw [fill] (2.55,2.1) circle [radius=0.05];
\draw [fill] (2.55,1.8) circle [radius=0.05];

\draw[draw=black,fill=gray!10](0.25,0)rectangle(2,1.5);\node at(1.125,.75){$\sigma$};
\draw[thick,->](2,0.3)
to[out=0,in=90](2.3,0)
to[out=-90,in=0](2,-0.3)--(-0.3,-0.3)
;
\draw[thick,->](2,.6)
to[out=0,in=90](2.6,0)
to[out=-90,in=0](2,-0.6)--(-0.3,-0.6)
;
\draw[thick,-](2,.9)
to[out=0,in=90](2.9,0)
to[out=-90,in=0](2,-0.9)--(-0.3,-0.9)
;
\draw[thick,->](2,1.2)
to[out=0,in=90](3.2,0)
to[out=-90,in=0](2,-1.2)--(-0.3,-1.2)
;

\draw [fill] (-0.3,-.9) circle [radius=0.05];
\end{tikzpicture}
$$
\begin{remark}
We defined duality on boundaries as not only flipping  left endpoints with right endpoints, but also reversing the list of boundary elements precisely in order to have  this consistency with ``parallel wires substitution'' in the graphical language. The price to pay for that is the twist of tensor factors in formula (\ref{tensor dual}).
\end{remark}
\bp
Formula (\ref{tensor dual}) holds for cowordisms as well. For all cowordisms $\sigma$, $\tau$
we have $(\sigma\otimes\tau)^\bot=\tau^\bot\otimes\sigma$.
\ep
\begin{proof}
  Evident from geometric representation.
\end{proof}
\bp
Duality is a contravariant functor: for all boundaries $X,Y,Z$ and cowordisms $$\sigma:X\to Y,\quad \tau:Y\to Z$$ we have $(\tau\circ\sigma)^\bot=\sigma^\bot\circ\tau^\bot$.
\ep
\begin{proof}
  Evident from geometric representation.
\end{proof}

\section{Structure of cowordisms category}
Since composition of cowordisms is associative and there are identities, as shown in the preceding section, it follows that cowordisms and boundaries form a category.
\bd The {\it category} ${\bf Cow}_T$ of cowordisms over an alphabet $T$ has boundaries as objects and cowordisms over $T$ as morphisms.
\ed

\subsection{Over the empty alphabet}
Even when the alphabet $T$ is empty, the category of cowordisms is nontrivial.

In fact, it becomes equivalent to the category of  oriented 1-dimensional topological cobordisms.

In the sequel we will use  the term {\it cobordism} for a cowordism over the empty alphabet, and
 denote $${\bf Cow_\emptyset}={\bf Cob}.$$

Given two boundaries $X,Y$ and a cowordism $\sigma:X\to Y$ over some alphabet $T$, we define the {\it pattern} of $\sigma$ as the cobordism from $X$ to $Y$ obtained by erasing from $\sigma$ all letters.
Note that pattern  is a functor (from cowordisms to cobordisms).

The category of cobordisms equipped with tensor product and duality has a rich structure, and this structure is inherited by categories of cowordisms over non-empty alphabets.

In particular, categories of cowordisms are {\it symmetric monoidal closed, $*$-autonomous}, and {\it compact closed}. This makes them  models of linear $\lambda$-calculus and of classical multiplicative linear logic, which is most relevant for our discussion.

Now we discuss this structure in a greater detail.

\subsection{Zoo of monoidal closed categories}
Recall that a {\it monoidal category} ${\bf C}$ is a category equipped with a {\it tensor product} bifunctor
$$\otimes:{\bf C}\times{\bf C}\to{\bf C}$$
together with a {\it tensor unit} ${\bf 1}$ and natural {\it associativity} transformations (see \cite{MacLane} for details).

A monoidal category is {\it symmetric} if there exists a natural {\it symmetry} transformation
$$s_{X,Y}:X\otimes Y\to Y\otimes X$$
satisfying $$s^2=\id.$$

A symmetric monoidal  category  is {\it closed} if it is  equipped with a bifunctor $\multimap$, contravariant in the first entry and covariant in the second entry,
such that there exists a natural bijection
\be\label{monoidal closure}
Hom(X\otimes Y, Z)\cong Hom(X,Y\multimap Z).
\ee

The functor $\multimap$ in the above definition is called {\it internal homs functor}.

A  {\it $*$-autonomous} category \cite{Barr} is a symmetric monoidal category   equipped with a
contravariant functor $(.)^\bot$, such that there is a natural isomorphism
$$A^{\bot\bot}\cong A$$
and a natural bijection
\be\label{*-autonomy}
Hom(X\otimes Y, Z)\cong Hom(X,(Y\otimes Z^\bot)^\bot).
\ee

Duality $(.)^\bot$ equips a $*$-autonomous category with a second monoidal structure. The {\it cotensor product} $\wp$ is defined by
\be\label{cotensor}
X\wp Y=(X^\bot\otimes Y^\bot)^\bot.
\ee
The neutral object for the cotensor product is
$$\bot={\bf 1}^\bot.$$

Any $*$-autonomous category is monoidal closed. The internal homs functor is defined by
\be\label{internal homs}
X\multimap Y=X^\bot\wp Y.
\ee


A {\it compact closed} or, simply, {\it compact category}  is a $*$-autonomous category for which duality commutes with tensor, i.e. such that there exist natural isomorphisms
$$X\wp Y\cong X\otimes Y,\quad {\bf 1}\cong \bot.$$
\begin{remark}
Compact categories were defined in \cite{KellyLaplaza}. For the  definition used above and its equivalence to the original definition of \cite{KellyLaplaza} see \cite{AbramskyCoecke}.
\end{remark}

 A prototypical example of a compact category is the category of finite-dimensional vector spaces with the usual tensor product and algebraic duality. Note, however, that in this case, as, in general, in the algebraic setting, duality is denoted as a star $(.)^*$.

 The category of topological cobordisms is another widely used example of a compact category.

Compact categories are relevant for our discussion, because categories of cowordisms are compact,  and compact structure provides a lot of important maps and constructions. A short and readable introduction into this subject  can be found, for example, in \cite{AbramskyCoecke}.
We  emphasize  that our graphical language for cowordisms is adapted from the general {\it pictorial language} customarily used for general compact categories, see \cite{Selinger}.

\subsection{Cowordisms as a compact category}\label{abstract nonsense}
Tensor product  equips the category of cowordisms with a  monoidal structure, which is symmetric.

Indeed, tensor product is functorial and {\it strictly} associative, i.e.
$$(X\otimes Y)\otimes Z=X\otimes (Y\otimes Z),$$
 moreover,
 the empty boundary ${\bf 1}$ is a {\it strict} tensor unit.

The symmetry transformation $s_{X,Y}$ discussed in the preceding section is obviously natural and it squares to the identity (Proposition \ref{symmetry squared}), as required.

\bp
Duality $(.)^\bot$ equips the category of cowordisms with a compact closed (hence monoidal closed and $*$-autonomous) structure.
\ep
\begin{proof}
The set of cowordisms from $X\otimes Y$ to $Z$ is by definition the set of multiwords with the boundary
$Z\otimes(X\otimes Y)^\bot=Z\otimes Y^\bot\otimes X^\bot$.

The set of cowordisms from $X$ to $(Y\otimes Z^\bot)^\bot$ is easily seen to be precisely the same set.

Thus formula (\ref{monoidal closure}) holds and the category of cowordisms is $*$-autonomous.

Also, the symmetry transformation gives us a natural isomorphism
$$s_{Y^\bot,X^\bot}:(X\otimes Y)^\bot\to X^\bot\otimes Y^\bot,$$
so  the category of cowordisms is compact closed.
\end{proof}

The geometric meaning of bijection (\ref{*-autonomy}) is transparent.
$$
\begin{tikzpicture}[scale=.7]
\draw[thick](0,.25)--(-0.5,.25);

\node[left] at(-0.5,.25){$X$};

\draw[thick](0,.75)--(-0.5,.75);

\node[left] at(-0.5,.75){$Y$};

\draw[thick](1,.5)--(1.5,.5);
\node[right] at(1.5,0.25){$Z$};

\draw[draw=black,fill=gray!10](0.,0)rectangle(1,1);
\node at (2.5,.5){$\cong$};
\begin{scope}[shift={(4,-.6)}]
\draw[thick](0,.25)--(-0.5,.25);

\node[left] at(-0.5,.25){$X$};

\draw[thick](0,.75)
to[out=180,in=-90](-.5,1.25)
to[out=90,in=180](0,1.75)
--(1.5,1.75)
;

\node[right] at(1.5,1.75){$Y^\bot$};
\draw[thick](1,.5)--(1.5,.5);
\node[right] at(1.5,0.5){$Z$};

\draw[draw=black,fill=gray!10](0.,0)rectangle(1,1);
\end{scope}


\end{tikzpicture}
\bigskip
$$
Since the category of cowordisms is also $*$-autonomous, it has the cotensor product and internal homs defined by
(\ref{cotensor}) and (\ref{internal homs}).

It is easy to see that, in the category of cowordisms, we have
 $$X\wp Y=Y\otimes X,\quad X\multimap Y=Y\otimes X^\bot.$$ This formula will be much used in the sequel.

 As for the cotensor unit $\bot$, we have $\bot={\bf 1}^\bot={\bf 1}$.

 Anticipating (linear) logic interpretation, we note that correspondence (\ref{*-autonomy}) together with the above picture simply express  the possibility of moving formulas of a sequent between the right and the left sides of the turnstile.

We now discuss important constructions determined by the compact structure.

\subsubsection{Names and conames}
      In any $*$-autonomous category   there  are operations of {\it naming} and {\it conaming} (the first one exists also in any monoidal closed category).

Let
$\sigma:X\to Y$ be a morphism in a monoidal closed category ${\bf C}$.

Correspondence (\ref{monoidal closure}) together with the isomorphism
$$X\cong{\bf 1}\otimes X$$ yields the morphism
$$\ulcorner\sigma\urcorner:{\bf 1}\to X\multimap Y,$$
 called the {\it name}
of $\sigma$.

If ${\bf C}$ is furthermore $*$-autonomous, then  the name
$\ulcorner\sigma\urcorner$ can be written as a morphism
$$\ulcorner\sigma\urcorner:{\bf 1}\to X^\bot\wp Y,$$
 and the {\it coname}
$\llcorner\sigma\lrcorner$
of $\sigma$ is the morphism
$$\llcorner\sigma\lrcorner: X\otimes Y^\bot\to\bot,$$
defined by $\llcorner\sigma\lrcorner=(\ulcorner\sigma^\bot\urcorner)^\bot$.

In the case of cowordisms, the name  and the coname of a cowordism $\sigma:X\to Y$ are represented by the following picture.
$$
\begin{tikzpicture}[scale=.5]

\node[left] at(-1.2,.75){$\sigma:$};

\draw[thick](0,.75)--(-0.5,.75);

\node[left] at(-0.5,.75){$X$};

\draw[thick](1,.75)--(1.5,.75);
\node[right] at(1.5,0.75){$Y$};

\draw[draw=black,fill=gray!10](0.,.5)rectangle(1,1);\node at(.5,.75){$\sigma$};

\begin{scope}[shift={(6,-.6)}]
\node[left] at(-.5,1.25){$\ulcorner\sigma\urcorner:$};

\draw[thick](0,.75)
to[out=180,in=-90](-.5,1.25)
to[out=90,in=180](0,1.75)
--(1.5,1.75)
;

\node[right] at(1.5,1.75){$X^\bot$};
\draw[thick](1,.75)--(1.5,.75);
\node[right] at(1.5,0.75){$Y$};

\draw[draw=black,fill=gray!10](0.,.5)rectangle(1,1);\node at(.5,.75){$\sigma$};
\end{scope}

\begin{scope}[shift={(14,.4)}] 

\node[left] at(-1.5,.25){$\llcorner\sigma\lrcorner:$};

\draw[thick](0,.-.25)--(-0.5,-.25);

\node[left] at(-0.4,.75){$Y^\bot$};

\draw[thick](1,-.25)
to[out=0,in=-90](1.5,.25)
to[out=90,in=0](0,.75)
--(-.5,.75)
;
\node[left] at(-.6,-0.25){$X$};

\draw[draw=black,fill=gray!10](0.,-.5)rectangle(1,0);\node at(.5,-.25){$\sigma$};

\end{scope}

\end{tikzpicture}
$$
Again, in anticipation of linear logic interpretation, we note that naming and conaming of a morphism correspond simply to putting all formulas of a sequent to one side of the  turnstile.

Especially important are the name and the coname of the identity, respectively, the {\it pairing}
\be\label{pairing}
\llcorner\id_X\lrcorner:X^\bot\otimes X\to{ \bot},
\ee
and the {\it copairing} maps
\be\label{copairing}
\ulcorner\id_X\urcorner:{\bf 1}\to X^\bot\wp X
\ee (existing in any $*$-autonomous category).

In the category of cowordisms these are represented by  the following picture (note that wires are depicted in the correct order!).
$$
\begin{tikzpicture}[scale=.35]
\begin{scope}[shift={(0,-.6)}]
\node[left] at(.25,1.25){$\ulcorner\id_X\urcorner:$};

\draw[thick](1.5,.75)
to[out=180,in=-90](.5,1.25)
to[out=90,in=180](1.5,1.75)
;

\node[right] at(1.5,1.75){$X^\bot$};

\node[right] at(1.5,0.75){$X$};

\end{scope}

\begin{scope}[shift={(10,.4)}]

\node[left] at(-.55,.25){$\llcorner\id_X\lrcorner:$};

\node[left] at(1,.75){$X$};

\draw[thick](1,.75)
to[out=0,in=90](2,.25)
to[out=-90,in=0](1,-.25)
;
\node[left] at(1,-0.25){$X^\bot$};


\end{scope}

\end{tikzpicture}
$$

\subsubsection{Evaluation and linear distributivity}
In any monoidal closed category there exists a natural {\it evaluation } map
\be\label{evaluation}
\ev_{X,Y}:(X\multimap Y)\otimes X\to Y
\ee
obtained by correspondence (\ref{monoidal closure}) from $\id_{X\multimap Y}$.

In the category of cowordisms evaluation is represented by the following picture.
$$
\begin{tikzpicture}[scale=.7]
\node[left] at(-.7,.5){$\ev_{X,Y}:$};
  \draw[thick](0,1) to[out=0, in=90](.5,.75)
  to[out=-90, in=0](0,.5);
  \draw[thick](0,0) to[out=0, in=180](1,.5);
  \node[left] at(0,1){$X$};
  \node[left] at(0,.5){$X^\bot$};
\node[left] at(0,.0){$Y$};
\node[right] at(1,.5){$Y$};
\end{tikzpicture}
$$
Furthermore there is a natural map
$$(X\multimap Y)\otimes Z\to X\multimap (Y\otimes Z)$$
obtained from $\ev_{X,Y}\otimes\id_Z$  using symmetries and (\ref{monoidal closure}), which, in the case when the category is $*$-autonomous, gives rise to the {\it linear distributivity} map \cite{CocketSeely}
\be\label{linear distributivity}
\delta_{X,Y,Z}:(X\wp Y)\otimes Z\to X\wp(Y\otimes Z)
\ee
(substituting $X$ for $X^\bot$).

Linear distributivity by compositions with itself and symmetries gives rise to the {\it internal tensor}
\be\label{internal tensor}
\delta_{X,Y,Z,T}:(X\wp Y)\otimes (Z\wp T)\to X\wp(Y\otimes Z)\wp T.
\ee
and the {\it internal cotensor}
\be\label{internal cotensor}
\epsilon_{X,Y,Z,T}=\delta^\bot_{X^\bot,Y^\bot,Z^\bot,T^\bot}:X\otimes (Y\wp Z)\otimes T\to (X\otimes Y)\wp (Z\otimes T)
\ee
maps.

In the case of cowordisms, where cotensor and tensor differ only by factors ordering, all three maps are just permutation maps obtained as compositions of symmetries.

\subsubsection{Generalized compositions}
Any $*$-autonomous category, and the category of cowordisms in particular, has a number of operations that can be considered as generalized composition in disguise.

First, in any monoidal category, given morphisms
$$\sigma:X\to Y, \quad\tau: Y\otimes T\to Z,$$
we can construct the {\it partial composition} of $\sigma$ and $\tau$ {\it over} $T$
$$\tau\circ_Y\sigma:X\otimes T\to Z,$$
given as $$\tau\circ_Y\sigma=\tau\circ(\sigma\otimes\id_T).$$

In the case of cowordisms, partial composition is represented by the following picture.
$$\begin{tikzpicture}[scale=.7]

\node[left] at(-1,.25){$\sigma:$};
\draw[thick](0,.25)--(-0.5,.25);

\node[left] at(-0.5,.25){$X$};

\draw[thick](1,.25)--(1.5,.25);
\node[right] at(1.5,0.25){$Y$};

\draw[draw=black,fill=gray!10](0.,0)rectangle(1,.5);\node at(.5,.25){$\sigma$};

\begin{scope}[shift={(4.7,-.25)}]

\node[left] at(-1,.5){$\tau:$};
\draw[thick](0,.25)--(-0.5,.25);
\draw[thick](0,.75)--(-0.5,.75);
\node[left] at(-.5,.75){$T$};
\node[left] at(-0.5,.25){$Y$};

\draw[thick](1,.5)--(1.5,.5);
\node[right] at(1.5,0.5){$Z$};

\draw[draw=black,fill=gray!10](0.,0)rectangle(1,1);\node at(.5,.5){$\tau$};
\end{scope}
\begin{scope}[shift={(2,-.25)}]
\begin{scope}[shift={(8.5,0)}]
\node[left] at(-1,.5){$\tau\circ_Y\sigma:$};
\draw[thick](0,.25)--(-0.5,.25);

\node[left] at(-0.5,.25){$X$};

\draw[thick](1,.25)--(1.5,.25);

\draw[draw=black,fill=gray!10](0.,0)rectangle(1,.5);\node at(.5,.25){$\sigma$};
\draw[thick](-0.5,.75)--(2.5,0.75);
\node[left] at(-0.5,.75){$T$};
\begin{scope}[shift={(1.5,0)}]

\draw[thick](1,.5)--(1.5,.5);
\node[right] at(1.5,0.5){$Z$};

\draw[draw=black,fill=gray!10](0.,0)rectangle(1,1);\node at(.5,.5){$\tau$};
\end{scope}
\end{scope}

\end{scope}
\end{tikzpicture}
$$
%
%
%
%
%
%
Note that usual composition is a particular case: $\tau\circ\sigma=\tau\circ_{\bf 1}\sigma$.

If the category is furthermore $*$-autonomous, then,
given  two morphisms
 $$\sigma:{\bf 1}\to X\wp Y,\quad\tau:{\bf 1}\to Y^\bot\wp Z,$$
we can define the {\it  partial pairing} or the {\it cut}
$$\langle\sigma,\tau\rangle_Y:{\bf 1}\to X\wp Z$$
 of $\sigma$ and $\tau$ over $Y$
 by
$$
\langle\sigma,\tau\rangle_Y=(\id_X\wp (\llcorner\id_Y\lrcorner)\wp\id_Z)\circ\delta_{X,Y,Y^\bot,Z}\circ(\sigma\otimes\tau).
$$

In the case of cowordisms, cut is pictured as follows.
$$
\begin{tikzpicture}[scale=.7]
\node[left] at(-.2,0.5){$\sigma:$};

\draw[thick](1,.25)--(1.5,.25);
\draw[thick](1,.75)--(1.5,.75);
\node[right] at(1.5,0.75){$X$};
\node[right] at(1.5,0.25){$Y$};

\draw[draw=black,fill=gray!10](0.,0)rectangle(1,1);\node at(.5,.5){$\sigma$};

\begin{scope}[shift={(0,-1.2)}]

\node[left] at(-.2,0.5){$\tau:$};

\draw[thick](1,.25)--(1.5,.25);
\draw[thick](1,.75)--(1.5,.75);
\node[right] at(1.5,0.75){$Y^\bot$};
\node[right] at(1.5,0.25){$Z$};

\draw[draw=black,fill=gray!10](0.,0)rectangle(1,1);\node at(.5,.5){$\tau$};

\end{scope}

\begin{scope}[shift={(6,0)}]
\node[left] at(-.2,-.15){$\langle\sigma,\tau\rangle_Y:$};

\draw[thick](1,.25)
to[out=0,in=90](1.5,-.1)
to[out=-90,in=0](1.,-.45);
\draw[thick](1,.75)--(1.5,.75);
\node[right] at(1.5,0.75){$X$};

\draw[draw=black,fill=gray!10](0.,0)rectangle(1,1);\node at(.5,.5){$\sigma$};

\begin{scope}[shift={(0,-1.2)}]

\draw[thick](1,.25)--(1.5,.25);
\node[right] at(1.5,0.25){$Z$};

\draw[draw=black,fill=gray!10](0.,0)rectangle(1,1);\node at(.5,.5){$\tau$};

\end{scope}
\end{scope}
\end{tikzpicture}
\bigskip
$$

Partial pairing or cut can be understood as generalized composition in disguise because of the following observation.
\bp\label{cut is composition}

 For any two cowordisms
  $$\sigma:X\to Y,\quad\tau:Y\to Z,$$ it holds that
    $$\ulcorner\tau\circ\sigma\urcorner=\langle\ulcorner\sigma\urcorner,\ulcorner\tau\urcorner\rangle_{Y}. \quad\Box$$
\ep
\begin{proof}
  Evident from geometric representation.
\end{proof}

From (linear) logic point of view, partial composition corresponds to a cut between two-sided sequents, and partial pairing, to a cut between one-sided sequents. The above proposition expresses basic relation between one-sided and two-sided sequent calculus formulations.

%

\begin{remark}
  Propositions of this section are, of course, true for any $*$-autonomous, let alone compact, category.
  In fact, it can be shown that graphical reasoning is valid for any compact category, see \cite{Selinger}.
  However, in the case of cowordisms, which are {\it defined} as geometric objects, graphical representation is {\it literal} and does not require further justification.

  Also, in the end of this abstract section it might be reasonable to try to put our category of cowordisms into the general abstract context of compact categories.

  It is easy to see that  labeling and ordering boundary vertices is purely decorative; it allows us having better pictures. Eventually, if we define the size of a boundary $X$ as a pair consisting  of the size of the left part $X_l$ and of the right part $X_r$, then it is the only invariant: all boundaries with the same size are (canonically) isomorphic.

   Thus, up to equivalence of categories, we can consider boundaries to be just pairs of natural numbers. It can be shown then that the category of cowordisms over an alphabet $T$
 is equivalent to the {\it free compact category} (using definition in  \cite{AbramskyFreeTraced}) generated by the free monoid $T^*$, where  $T^*$ is considered as a category with one object.
 \end{remark}

\section{Linear logic grammars}
\subsection{Linear logic}
   Strictly speaking, the system discussed below is {\it multiplicative linear logic}, a fragment of full linear logic. However, since we do not consider other fragments, the prefix ``multiplicative'' will be omitted. A more detailed introduction to linear logic can be found in \cite{Girard}, \cite{Girard2}.

Given a set $N$ of  {\it positive literals}, we define the set $N^\bot$ of {\it negative literals}
as
$$N^\bot=\{X^\bot|\mbox{ }X\in N\}.$$
Elements of $N\cup N^\bot$ will be called {\it literals}.

The set  $Fm(N)$ of  ${\bf LL}$ formulas (over the alphabet $N$) is defined by the following induction.
\begin{itemize}
  \item Any $X\in N\cup N^\bot$ is a formula;
  \item if $X$, $Y$ are  formulas, then $X\wp Y$ and $X\otimes Y$ are   formulas;
  \end{itemize}

Connectives $\otimes$ and $\wp$ are called respectively {\it times} (also {\it tensor}) and {\it par} (also {\it cotensor}).

{\it Linear negation} $A^\bot$ of a
formula $A$ is defined
inductively as
$$(P^\bot)^\bot=P,\mbox{ for }P\in N,$$
$$(A\otimes B)^\bot=A^\bot\wp B^\bot,\quad
(A\wp B)^\bot=A^\bot\otimes B^\bot.$$
It should be understood that linear negation $(.)^\bot$ is not a connective of linear logic, but a shorthand notation for formulas, i.e. it belongs to the metalanguage.

{\it Linear implication} is defined as
\be\label{linear implication}
A\multimap B=A^\bot\wp B.
\ee
Linear implication is not a connective either.

An  ${\bf LL}$ {\it sequent} (over the alphabet $N$) is a finite sequence of ${\bf LL}$ formulas (over $N$).

The {\it sequent calculus} for ${\bf LL}$ is given by the following rules:
$$\vdash X^\bot,X~(\rm{Id}),\quad
\frac{\vdash \Gamma,X\quad\vdash
X^\bot,\Delta}{\vdash\Gamma,\Delta} ~(\rm{Cut}),$$
$$\frac{\vdash
\Gamma,X,Y,\Delta}{\vdash
\Gamma,Y,X,\Delta}
~(\rm{Ex}),$$
$$\frac{\vdash \Gamma,X,Y}
{\vdash\Gamma,X\wp Y}~
(\wp)\quad\frac{\vdash\Gamma, X \quad\vdash
Y,\Delta}{\vdash\Gamma,X\otimes Y,\Delta}~{ }
(\otimes).$$
\smallskip

The {\it tensor} and {\it cotensor introduction} rules, respectively $(\otimes)$ and $(\wp)$ are called {\it logical} rules. Other rules are {\it structural}.

 Notation (Ex) in the name of a rule stands for Exchange.

Linear logic  enjoys the fundamental property of {\it cut-elimination}: any proof with cuts can be transformed to its {\it cut-free} form;  cut-elimination is algorithmic and always terminates.
This allows computational and categorical interpretations  in the {\it proofs-as-programs} or {\it proofs-as-functions}  paradigm.

\subsection{Semantics}
Categorical interpretation of proof theory is based on the idea that formulas should be understood as objects, and proofs as morphisms in a category, while composition of morphisms corresponds to cut-elimination.

In a two-sided sequent calculus, formulas are interpreted as objects in a monoidal category, and a proof of the sequent
$$X_1,\ldots, X_n\vdash X$$
 is interpreted as a morphism of type
 $$X_1\otimes\ldots\otimes X_n\to X.$$
 This includes the case $n=0$, with the usual convention that the tensor of the empty collection of objects is the monoidal unit ${\bf 1}$.

Then the Cut rule corresponds to (partial)  composition. A crucial requirement is that the interpretation should be invariant under cut-elimination; a proof and its cut-free form are interpreted the same.

In the case of linear logic, whose sequents are one-sided, the appropriate setting for categorical interpretation   is   {\it $*$-autonomous categories} \cite{Seely}, \cite{Mellies_categorical_semantics}.

In this setting, a proof of the  sequent
$$
\vdash X_1,\ldots, X_n$$
 is interpreted as a morphism of type
 $${\bf 1}\to X_1\wp\ldots\wp X_n.$$
The Cut rule  corresponds to partial pairing, which can be understood as a symmetrized composition (Proposition \ref{cut is composition}).

We first describe interpretation of {\bf LL} in a general $*$-autonomous category and then give explicit rules for the case of cowordisms.

Given a $*$-autonomous category ${\bf C}$ and an alphabet $N$ of positive literals, an interpretation of ${\bf LL}$ in ${\bf C}$ consists in assigning to any positive literal $A$ an object $[A]$ of $\bf C$.
The assignment of objects  extends to all formulas in $Fm(N)$ by the obvious induction
$$[A\otimes B]=[A]\otimes [B],\quad [A^\bot]=[A]^\bot.$$
A sequent $\Gamma=A_1,\ldots A_n$ is interpreted as the object $[\Gamma]=[A_1]\wp\ldots\wp[A_n]$.

Given interpretation of formulas, proofs are interpreted by induction on the rules so that a proof $\pi$ of the sequent $\Gamma$ is interpreted as a morphism $[\pi]$ of the type ${\bf 1}\to[\Gamma]$.

The axiom $\vdash A^\bot,A$ is interpreted as   copairing  map (\ref{copairing})
$$\ulcorner \id_{[ A]}\urcorner:{\bf 1}\to [ A]^\bot\wp [ A].$$

The Cut rule  corresponds to partial pairing, as stated above.

The Exchange rule  corresponds to a symmetry transformation.

The $(\wp)$ rule does  strictly nothing. A proof of the sequent $\vdash\Gamma,A,B$ has already been interpreted
as a morphism  of the type ${\bf 1}\to[\Gamma]\wp[A]\wp[B]$, and this serves as the interpretation of the obtained proof of $\vdash\Gamma,A\wp B$ as well.

The $(\otimes)$ rule is interpreted using internal tensor (\ref{internal tensor}). Namely, given proofs $\pi$ and $\rho$ of sequents $\vdash\Gamma,A$ and $\vdash B,\Delta$
respectively, the proof of $\vdash \Gamma, A\otimes B, \Delta$ obtained from $\pi$ and $\rho$
is interpreted as $\delta_{\Gamma,A,B,\Delta}\circ\pi\otimes\rho$.

In general, it is quite customary in the literature to omit square brackets and denote a formula and its interpretation by the same expression, and  we will follow this practice when convenient.

\subsubsection{Interpretation in the category of cowordisms}
We now specialize to the case of the $*$-autonomous category ${\bf Cow}_T$ of cowordisms over an alphabet $T$.

A proof $\sigma$ of the sequent $\vdash A_1,\ldots,A_n$ is interpreted as a cowordism
\be\label{flip in the sequent}
\sigma:{\bf 1}\to A_1\wp\ldots\wp A_n=A_n\otimes\ldots\otimes A_1
\ee
 (we use the common convention that a formula or a proof is denoted same as its interpretation).

We will use {\it vertical} representation for depicting $\sigma$, with inputs below, and outputs above. Since there are no inputs, $\sigma$ will be a box with only outgoing wires. Since, in the vertical representation, vertices of the outgoing boundary are depicted in the order from right to left, the wires to $A_1,\ldots,A_n$ will be depicted {\it in the same order} as in the sequent, as is evident from (\ref{flip in the sequent}).

Rules for interpreting proofs in the category of cowordisms can be easily computed from the general prescription  described above, using graphical representations of various natural maps computed in the preceding section.
The result is represented in four pictures below. (No picture for the $(\wp)$ rule, because it does not do anything, as noted above.)
$$
 
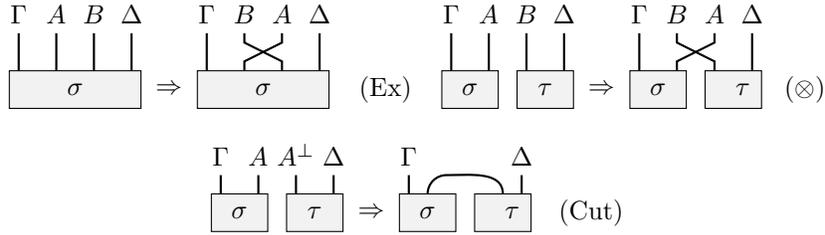
\tikz[scale=.5]{
 \draw[draw=black,fill=gray!10](-5,0)rectangle(-1.5,1);\node at(-3.25,.5){$\sigma$};
 \draw[thick,-](-4.75,1)--(-4.75,2);

 \node[above] at(-4.75,2.) {$\Gamma$};
 \draw[thick,-](-3.75,1)--(-3.75,2);
 \node[above] at(-3.75,2) {$A$};

 \draw[thick,-](-1.75,1)--(-1.75,2);

 \node[above] at(-1.75,2) {$\Delta$};
 \draw[thick,-](-2.75,1)--(-2.75,2);
 \node[above] at(-2.75,2) {$B$};
 \node at (-.75,.5) {$\Rightarrow$};

         \draw[draw=black,fill=gray!10](0,0)rectangle(3.5,1);\node at(1.75,.5){$\sigma$};
 \draw[thick,-](.25,1)--(.25,2);

 \node[above] at(.25,2.) {$\Gamma$};
 \draw[thick,-](1.25,1)--(1.25,1.25)--(2.25,1.75)--(2.25,2);
 \node[above] at(2.25,2) {$A$};

 \draw[thick,-](3.25,1)--(3.25,2);

 \node[above] at(3.25,2) {$\Delta$};
 \draw[thick,-](2.25,1)--(2.25,1.25)--(1.25,1.75)--(1.25,2);
 \node[above] at(1.25,2) {$B$};
 \node at (5,.5) {$(\rm{Ex})$};
 \begin{scope}[shift={(11.5,0)}]
 \draw[draw=black,fill=gray!10](-5,0)rectangle(-3.5,1);\node at(-4.25,.5){$\sigma$};
  \draw[draw=black,fill=gray!10](-3,0)rectangle(-1.5,1);\node at(-2.25,.5){$\tau$};
 \draw[thick,-](-4.75,1)--(-4.75,2);

 \node[above] at(-4.75,2.) {$\Gamma$};
 \draw[thick,-](-3.75,1)--(-3.75,2);
 \node[above] at(-3.75,2) {$A$};

 \draw[thick,-](-1.75,1)--(-1.75,2);

 \node[above] at(-1.75,2) {$\Delta$};
 \draw[thick,-](-2.75,1)--(-2.75,2);
 \node[above] at(-2.75,2) {$B$};
 \node at (-.75,.5) {$\Rightarrow$};

         \draw[draw=black,fill=gray!10](0,0)rectangle(1.5,1);\node at(.75,.5){$\sigma$};
 \draw[thick,-](.25,1)--(.25,2);

 \node[above] at(.25,2.) {$\Gamma$};
 \draw[thick,-](1.25,1)--(1.25,1.25)--(2.25,1.75)--(2.25,2);
 \node[above] at(2.25,2) {$A$};

         \draw[draw=black,fill=gray!10](2.,0)rectangle(3.5,1);\node at(3,.5){$\tau$};

 \draw[thick,-](3.25,1)--(3.25,2);

 \node[above] at(3.25,2) {$\Delta$};
 \draw[thick,-](2.25,1)--(2.25,1.25)--(1.25,1.75)--(1.25,2);
 \node[above] at(1.25,2) {$B$};
 \node at (4.65,.5) {$(\otimes)$};
  \end{scope}
  }
  $$
 $$
 \tikz[scale=.5]{
 \begin{scope}[shift={(21.5,0)}]
 \begin{scope}[shift={(1,0)}]

 \draw[draw=black,fill=gray!10](-5,0)rectangle(-3.5,1);\node at(-4.25,.5){$\sigma$};
  \draw[draw=black,fill=gray!10](-3,0)rectangle(-1.5,1);\node at(-2.25,.5){$\tau$};
 \draw[thick,-](-4.75,1)--(-4.75,1.5);

 \node[above] at(-4.75,1.5) {$\Gamma$};
 \draw[thick,-](-3.75,1)--(-3.75,1.5);
 \node[above] at(-3.75,1.5) {$A$};

 \draw[thick,-](-1.75,1)--(-1.75,1.5);

 \node[above] at(-1.75,1.5) {$\Delta$};
 \draw[thick,-](-2.75,1)--(-2.75,1.5);
 \node[above] at(-2.75,1.5) {$A^\bot$};
 \node at (-.75,.5) {$\Rightarrow$};
 \end{scope}

\begin{scope}[shift={(-4,0)}]
     \draw[draw=black,fill=gray!10](5,0)rectangle(6.5,1);\node at(5.75,.5){$\sigma$};
 \draw[thick,-](5.25,1)--(5.25,1.5);

 \node[above] at(5.25,1.5) {$\Gamma$};
 \draw[thick,-](5.75,1)to[out=90,in=180](6.25,1.5)
 to[out=0,in=90](7.75,1);

         \draw[draw=black,fill=gray!10](7.,0)rectangle(8.5,1);\node at(8,.5){$\tau$};

 \draw[thick,-](8.25,1)--(8.25,1.5);

 \node[above] at(8.25,1.5) {$\Delta$};
 \node at (10.1,.5) {$(\rm{Cut})$};
 \end{scope}
 \end{scope}
 }$$
(Note that the order of wires on the righthand side of the $(\otimes)$ rule picture is correct. The outgoing boundary is $\Gamma\wp(A\otimes B)\wp\Delta=\Delta\otimes(A\otimes B)\otimes\Gamma$, and according to our convention tensor factors are depicted in the opposite order.)

Observe that
  interpretations of proofs  do not depend on  the alphabet $T$ at all. So it would be more honest to say that this is  an interpretation  in the category ${\bf Cob}$ of {\it cobordisms}. The alphabet comes into play if we add new axioms to the logic, which gives us a {\it logic grammar}.

  \begin{remark}
     It is important to note that the model of linear logic in the category of cowordisms, as well as in any other compact category, is {\it degenerate}. It identifies (up to a natural isomorphism) tensor and cotensor.
  \end{remark}

\subsection{Adding lexicon}
An ${\bf LL}$ grammar is an interpretation of ${\bf LL}$ in a category of cowordisms supplied with a set of axioms together with cowordisms representing their ``proofs''.

So assume that an interpretation of ${\bf LL}$ in a category of cowordisms over an alphabet $T$ is given.
We will use the common convention that a formula, a sequent or a proof is denoted the same as its interpretation.

 Let us say that a {\it cowordism typing judgement} (over $N$ and $T$) is an expression of the form $$\frac{\sigma}{\vdash F},$$ where $F$ is an ${\bf LL}$ sequent (over the alphabet $N$), and $$\sigma:{\bf 1}\to F$$ is a cowordism over $T$.

\bd
A {\it cowordism signature} $\Sigma$ is a tuple $\Sigma=(N,T,\Xi)$, where
\begin{itemize}
  \item $N$ is a finite set of positive literals interpreted as boundaries;
  \item $T$ is a finite alphabet;
  \item $\Xi$ is a set of cowordism typing judgements over $N$ and $T$, called {\it signature axioms} or, simply, {\it axioms}.
  \end{itemize}
\ed

Given a cowordism signature  $\Sigma$, any {\bf LL} derivation from signature  axioms  gets an interpretation in the category of cowordisms by induction on the rules.

Let $\Gamma$ be an ${\bf LL}$ sequent (or formula), and let $\sigma:{\bf 1}\to \Gamma$ be a   cowordism.

We say  that
the cowordism typing judgement $$\frac{\sigma}{\vdash \Gamma}$$ is {\it derivable in $\Sigma$}, or that
$\Sigma$ {\it generates  cowordism $\sigma$  of type $\Gamma$}
if there exists a derivation of
$\vdash \Gamma$ from axioms of $\Sigma$ whose interpretation is $\sigma$.
The {\it cowordism type} $\Gamma$ {\it generated by $\Sigma$}, or, simply, the {\it cowordism type $\Gamma$ of $\Sigma$}, is the set of all cowordisms of type $\Gamma$ generated by $\Sigma$.

It is convenient to have an alternative characterization of derivability in an LLG.

For a sequent $\Gamma=A_1,\ldots,A_n$ let $\Gamma^\bot$ denote the formula $\Gamma^\bot=A_1^\bot\otimes\ldots\otimes A_n^\bot$.

Then we have the following.
\bp\label{alternative derivability}
Let $\Sigma$ be a cowordism signature.

A typing judgement $$\cfrac{\sigma}{\vdash\Gamma} $$ is derivable in $\Sigma$ iff there exist axioms
$$\cfrac{\sigma_1}{\vdash\Gamma_1},\ldots,\cfrac{\sigma_k}{\vdash\Gamma_k} $$ of $\Sigma$ and an {\bf LL} proof $\pi$ (not using  axioms of $\Sigma$) of the sequent $$\vdash \Gamma_1^\bot,\ldots,\Gamma^\bot_k,\Gamma$$ such that
$$\sigma=\pi\circ(\sigma_1\otimes\ldots\otimes\sigma_k)$$
(where $\pi$ is identified with the corresponding cowordism).
\ep
\begin{proof}
  By induction on derivation.
\end{proof}
\bc\label{ts-free=par-free}
If $\Sigma=(N,T,\Xi)$ is a cowordism signature and there is an axiom in $\Xi$  of the form
\be\label{sequent with par}
\frac{\sigma}{\vdash\Gamma,A\wp B,\Delta},
\ee
then we have a well-defined cowordism signature $\Sigma'$ obtained from $\Sigma$ by replacing (\ref{sequent with par}) with
 $$
\frac{\sigma}{\vdash\Gamma,A, B,\Delta},
$$
which generates same cowordism types as
 $\Sigma$. $\Box$
 \ec

\bd
A {\it linear logic grammar (LLG)} $G$ is a tuple $G=(N,T,Lex,S)$, where $(N,T,Lex)$ is a cowordism signature, and
 $S\in N$, the {\it standard type}, is a literal interpreted as a boundary with exactly one left and one right endpoint.
\ed

The set $Lex$ of the underlying signature axioms of $G$ is called a {\it lexicon} of $G$.

We say that a typing judgement, respectively, cowordism type is generated by $G$ if it is generated by the underlying cowordism signature.

 Now, any regular cowordism of the standard type $S$ generated by $G$ is an edge-labeled graph containing  a single edge labeled with a word over $T$. Thus the set of type $S$ regular cowordisms can be identified with a set of words in $T^*$.

The {\it language $L(G)$ generated by $G$} is the set of words labeling type $S$ regular cowordisms generated by $G$.

\begin{remark}
  We noted above that the model of linear logic in the category of cowordisms is degenerate.

  However, given an LLG (or just a cowordism signature) $G$, we can consider a {\it refined} category of cowordisms, whose objects are types of $G$ and morphisms are cowordisms generated by $G$. (More accurately, a morphism between types $A$ and $B$ is a cowordism
  $\sigma:A\to B$ such that the name $\ulcorner\sigma\urcorner$ is in the type $A\multimap B$ of $G$.)

  This gives us another model, which, in a generic case, will be nondegenerate.

  Thus, from the point of view of {\bf LL} model theory, an LLG (or a cowordism signature) is a way of refining compact closed structure in order to remove  degeneracy of the model. This  approach to constructing models is very common, see, for example, \cite{Tan}, \cite{HS}.
\end{remark}

\subsubsection{Example}
Consider the following alphabets  $T$, $N$ of terminal symbols  and positive literals respectively:
$$T=\{\rm{John},\rm{Mary},\rm{leaves},\rm{loves},\rm{madly},\rm{who}\},\quad
N=\{NP,S\}.$$
Let us construct a grammar to derive the sentence
\be\label{toy sentence}
\mbox{Mary who John loves madly leaves.}
\ee

We interpret both nonterminal symbols as the two-point boundary $X$ with $|X|=2$, $X_l=\{1\}$.

Then we introduce the following lexicon:
\[\left\{\frac{JOHN}{\vdash NP},\frac{MARY}{\vdash NP},\frac{LEAVES}{\vdash NP^\bot,S},\right.\]
\[\left.\frac{LOVES}{\vdash NP^\bot,NP^\bot,S},\frac{MADLY}{\vdash NP\otimes S^\bot, NP^\bot,S},\frac{WHO}{\vdash NP\otimes S^\bot,NP^\bot,NP}\right\}.
\]
(The types in the lexicon become less puzzling, if we identify every sequent with the cotensor of its formulas, and then rewrite it using linear implication (\ref{linear implication}). We get
$NP$, $NP\multimap S$, $NP\multimap NP\multimap S$, $(NP\multimap S)\multimap NP\multimap S$ etc.)

Cowordisms of the lexicon are depicted below.
$$
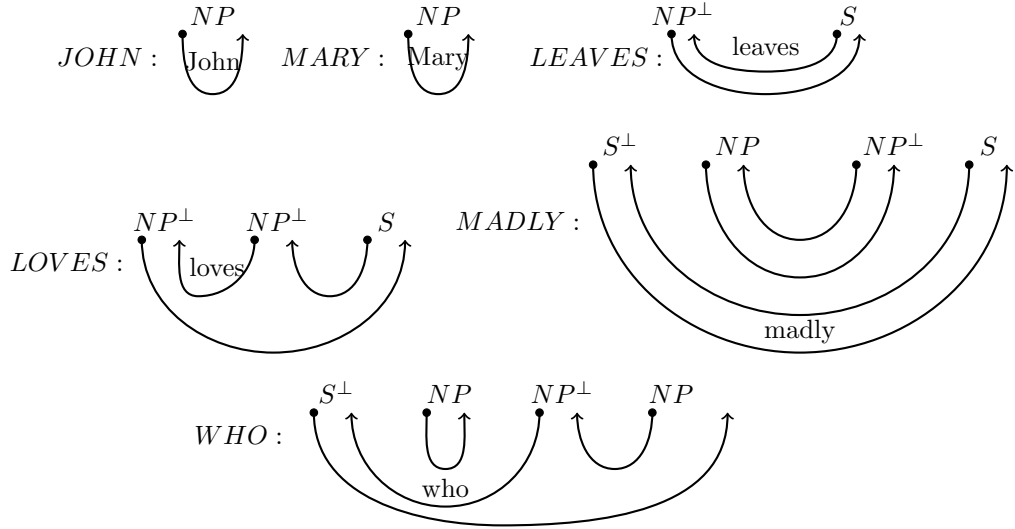
\begin{tikzpicture}
 \draw[thick,->](0,0) to  [out=-90,in=180] (0.4,-0.8) to  [out=0,in=-90] (.8,0);
 \draw [fill] (0,0) circle [radius=0.05];
  \node[above ] at (0.4,0) {$NP$};
   \node[above] at (0.4,-.6) {$\rm{John}$};
\node at (-1,-0.3){$JOHN:$};

\begin{scope}[shift ={(3,0)}]
 \draw[thick,->](0,0) to  [out=-90,in=180] (0.4,-0.8) to  [out=0,in=-90] (0.8,0);
 \draw [fill] (0,0) circle [radius=0.05];
  \node[above ] at (0.4,0) {$NP$};
   \node[above] at (0.4,-.6) {$\rm{Mary}$};
\node at (-1,-0.3){$MARY:$};

 \begin{scope}[shift={(3.5,0)}]
 \draw[thick,->](0,0) to  [out=-90,in=180] (1.25,-.8) to  [out=0,in=-90] (2.5,0);
 \draw [fill] (0,0.) circle [radius=0.05];
\draw[thick,<-](0.3,0) to  [out=-90,in=180] (1.25,-.5) to  [out=0,in=-90] (2.2,0);
\draw [fill] (2.2,0) circle [radius=0.05];
\node[above] at (1.25,-.4) {$\rm{leaves}$};
  \node[above] at (0.15,0) {$NP^\bot$};
\node[above] at (2.35,0) {$S$};

\node at (-1,-0.3){$LEAVES:$};
\end{scope}
 \end{scope}

\end{tikzpicture}
$$
$$
\begin{tikzpicture}
 \draw[thick,->](0,0) to  [out=-90,in=180] (1.75,-1.5) to  [out=0,in=-90] (3.5,0);
 \draw [fill] (0,0) circle [radius=0.05];
\draw[thick,<-](0.5,0) to  [out=-90,in=180] (0.75,-0.75) to  [out=0,in=-90] (1.5,0);
\draw [fill] (3,0) circle [radius=0.05];
\draw[thick,<-](2,0) to  [out=-90,in=180] (2.5,-0.75) to  [out=0,in=-90] (3,0);
\draw [fill] (1.5,0) circle [radius=0.05];
\node[above] at (1,-.6) {$\rm{loves}$};
  \node[above] at (0.3,0) {$NP^\bot$};
\node[above] at (1.8,0) {$NP^\bot$};
\node[above ] at (3.25,0) {$S$};
\node at (-1,-0.3){$LOVES:$};

\begin{scope}[shift={(6,)}]
 \draw[thick,->](0,0) to  [out=-90,in=180] (2.75,-2.5) to  [out=0,in=-90] (5.5,0);
\draw[thick,<-](0.5,0) to  [out=-90,in=180] (2.75,-2) to  [out=0,in=-90] (5,0);
\draw [fill] (0,0) circle [radius=0.05];
\draw [fill] (5,0) circle [radius=0.05];

\draw[thick,->](1.5,0) to  [out=-90,in=180] (2.75,-1.5) to  [out=0,in=-90] (4,0);
\draw[thick,<-](2,0) to  [out=-90,in=180] (2.75,-1) to  [out=0,in=-90] (3.5,0);
\node[above] at (2.75,-2.5) {$\rm{madly}$};
\draw [fill] (1.5,0) circle [radius=0.05];
\draw [fill] (3.5,0) circle [radius=0.05];

  \node[above right] at (0,0) {$S^\bot$};
\node[above left] at (5.5,0) {$S$};
\node[above right] at (1.5,0) {$NP$};
\node[above] at (4,0) {$NP^\bot$};
\node at (-1,-0.75){$MADLY:$};
\end{scope}
\end{tikzpicture}
$$
$$
\begin{tikzpicture}
\begin{scope}[shift={(6,0)}]
\draw[thick,->](0,0) to  [out=-90,in=180] (2.5,-1.5) to  [out=0,in=-90] (5.5,0);
\draw [fill] (0,0) circle [radius=0.05];
\node[above] at (0.3,0) {$S^\bot$};

 \draw[thick,->](1.5,0) to  [out=-90,in=180] (1.75,-.75) to  [out=0,in=-90] (2,0);
 \draw [fill] (1.5,0) circle [radius=0.05];
\draw[thick,<-](0.5,0) to  [out=-90,in=180] (1.75,-1.25) to  [out=0,in=-90] (3,0);
\draw [fill] (3,0) circle [radius=0.05];
\draw[thick,<-](3.5,0) to  [out=-90,in=180] (4,-0.75) to  [out=0,in=-90] (4.5,0);
\draw [fill] (4.5,0) circle [radius=0.05];
\node[above] at (1.75,-1.25) {$\rm{who}$};
  \node[above] at (1.8,0) {$NP$};
\node[above] at (3.3,0) {$NP^\bot$};
\node[above ] at (4.75,0) {$NP$};
\node at (-1,-0.3){$WHO:$};
\end{scope}
\end{tikzpicture}
$$

This defines an LLG.

Let us start with the following  derivation.
    \[
    \scalebox{.75}[1]
        {
        \text
            {
$
\cfrac
                            {
                            \cfrac
                                {
                                \cfrac
                                    {
                                    JOHN
                                    }
                                    {
                                    \vdash NP
                                    }
                                ~
                                \cfrac
                                    {
                                    \cfrac
                                        {
                                        \cfrac
                                            {
                                            \cfrac
                                                {
                                                LOVES
                                                }
                                                {
                                                \vdash NP^\bot,NP^\bot,S
                                                }
                                            }
                                            {
                                            \vdash NP^\bot,NP^\bot\wp S
                                            }(\wp)
                                        ~
                                        \cfrac
                                            {
                                            MADLY
                                            }
                                            {
                                            \vdash NP\otimes S^\bot,NP^\bot,S
                                            }
                                        }
                                        {
                                        \vdash NP^\bot,NP^\bot,S
                                        }(\rm{Cut})
                                    }
                                    {
                                    \vdash NP^\bot,NP^\bot,S
                                    }(\rm{Ex})
                                }
                                {
                                \vdash NP^\bot,S
                                }(\rm{Cut})
                            }
                            {
                            \vdash NP^\bot\wp S
                            }(\wp)
            $
            }
        }
    \]

We start its cowordism interpretation from axioms corresponding to $LOVES$ and $MADLY$.

Then the first step is a cut  between these two. This is interpreted by partial pairing, and the picture is the following.
$$
\begin{tikzpicture}[xscale=.7]
 \draw[thick,->](0,0) to  [out=-90,in=180] (1.75,-1.5) to  [out=0,in=-90] (3.5,0)
 to [out=90,in=180] (4,0.5) to [out=0,in=90]
 (4.5,0) to  [out=-90,in=180] (7.25,-2.5) to  [out=0,in=-90] (10,0);
 \draw [fill] (0,0) circle [radius=0.05];
\draw[thick,<-](0.5,0) to  [out=-90,in=180] (0.75,-0.75) to  [out=0,in=-90] (1.5,0)
 to [out=90,in=180] (4,2.5) to [out=0,in=90]
(6.5,0) to  [out=-90,in=180] (7.25,-1) to  [out=0,in=-90] (8,0);
;


\node[above] at (1,-.6) {$\rm{loves}$};
  \node[above] at (0.3,0) {$NP^\bot$};

\draw[thick,-](2,0) to  [out=-90,in=180] (2.5,-0.75) to  [out=0,in=-90] (3,0)
to [out=90,in=180] (4,1) to [out=0,in=90]
(5,0) to  [out=-90,in=180] (7.25,-2) to  [out=0,in=-90] (9.5,0);
\draw[thick,-](2,0) to  [out=90,in=180] (4, 2) to  [out=0,in=90] (6,0);

\draw [fill] (9.5,0) circle [radius=0.05];

\draw[thick,->](6,0) to  [out=-90,in=180] (7.25,-1.5) to  [out=0,in=-90] (8.5,0);
\node[above] at (7.25,-2.5) {$\rm{madly}$};
\draw [fill] (8,0) circle [radius=0.05];
\draw [fill] (8,0) circle [radius=0.05];

\node[above left] at (10,0) {$S$};
\node[above] at (8.5,0) {$NP^\bot$};
\node at(11.5,-1){$=$};

\begin{scope}[shift={(13,0)}]
\draw[thick,->](0,0) to  [out=-90,in=180] (1.75,-1.5) to  [out=0,in=-90] (3.5,0);
 \draw [fill] (0,0) circle [radius=0.05];
\draw[thick,<-](0.5,0) to  [out=-90,in=180] (0.75,-0.75) to  [out=0,in=-90] (1.5,0);
\draw [fill] (3,0) circle [radius=0.05];
\draw[thick,<-](2,0) to  [out=-90,in=180] (2.5,-0.75) to  [out=0,in=-90] (3,0);
\draw [fill] (1.5,0) circle [radius=0.05];
\node[above] at (1,-.6) {$\rm{loves}$};
  \node[above] at (0.3,0) {$NP^\bot$};
\node[above] at (1.8,0) {$NP^\bot$};
\node[above ] at (3.25,0) {$S$};
\node[above] at (1.75,-1.5) {$\rm{madly}$};
\end{scope}
\end{tikzpicture}
$$
Applying to this the Exchange rule and cutting with $JOHN$, we get the following cowordism interpreting the above derivation.
$$
\begin{tikzpicture}
 \draw[thick,-]
 (3,1) to  [out=90,in=0] (2,2) to  [out=180,in=90]
 (0,0.5) to  [out=-90,in=180] (0.75,-1.5) to  [out=0,in=-90] (1.5,0.5)
  to  [out=90,in=180] (2,1.5) to  [out=0,in=90] (2.5,1);
 \node[above] at (0.75,-1.5) {$\rm{John}$};

    \draw[thick,->](4,1)--(2.5,0) to  [out=-90,in=180] (4.25,-1.5) to  [out=0,in=-90] (6,1);

\draw[thick,<-](4.5,1)--(3,0) to  [out=-90,in=180] (3.25,-0.75) to  [out=0,in=-90] (4,0)--(2.5,1);

\draw[thick,-](3,1)--(4.5,0) to  [out=-90,in=180] (4.8,-0.75) to  [out=0,in=-90] (5.5,1);
\draw [fill] (5.5,1) circle [radius=0.05];
\node[above] at (3.5,-.6) {$\rm{loves}$};
\node[above ] at (5.75,1) {$S$};
\draw [fill] (4,1) circle [radius=0.05];
\node[above ] at (4.25,1) {$NP^\bot$};
\node[above] at (4.35,-1.5) {$\rm{madly}$};
\node at (6.25,-0.25) {$=$};

\draw[thick,->](6.5,.75) to  [out=-90,in=180] (9,-1.25) to  [out=0,in=-90] (11.5,.75);
\draw[thick,<-](7,.75) to  [out=-90,in=180] (9,-.75) to  [out=0,in=-90] (11,.75);
\draw [fill] (6.5,.75) circle [radius=0.05];
\draw [fill] (11,.75) circle [radius=0.05];
\node[above right ] at (11,.75) {$S$};
\node[above right ] at (6.5,.75) {$NP^\bot$};
\node[above] at (9,-.6) {$\rm{John}\mbox{ }\rm{loves}$};
\node[above] at (9,-1.22) {$\rm{madly}$};

\end{tikzpicture}
$$
We denote it as $JOHN\_LOVES\_MADLY$.

Now we continue with the derivation.
    \[
    \scalebox{.75}[1]
        {
        \text
            {
            $
            \cfrac
                {
                \cfrac
                    {
                    \cfrac
                        {
                        MARY
                        }
                        {
                        \vdash NP
                        }
                    ~
                    \cfrac
                        {
                        \cfrac
                            {
                            \cfrac
                                {
                                \cfrac
                                    {
                                   WHO
                                    }
                                    {
                                    \vdash NP^\bot\otimes S, NP^\bot,NP
                                    }
                                }
                                {
                                \vdash NP^\bot,NP\otimes S^\bot,NP
                                }(\rm{Ex})
                            }
                            {
                            \vdash NP^\bot,NP,NP\otimes S^\bot
                            }(\rm{Ex})
                        ~
                        \cfrac
                            {
                            JOHN\_LOVES\_MADLY
                            }
                            {
                            \vdash NP^\bot\wp S
                            }(\wp)
                        }
                        {
                        \vdash NP^\bot, NP
                        }(\rm{Cut})
                    }
                    {
                    \vdash NP
                    }(\rm{Cut})
                ~
                \cfrac
                    {
                    LEAVES
                    }
                    {
                    \vdash NP^\bot,S
                    }
                }
                {
                \vdash S
                }(\rm{Cut})
            $
            }
        }
    \]

We continue the cowordism interpretation  applying the series of Exchange rules to $WHO$, which results in the following.
$$
\begin{tikzpicture}
\draw[thick,->](5.3,0) to  [out=-90,in=180] (5.75,-0.75) to  [out=0,in=-90] (6.2,0);
\draw [fill] (5.3,0) circle [radius=0.05];
\node[above] at (5.8,0) {$NP$};
\node[above] at (2.75,-1.5) {$\rm{who}$};

 \draw[thick,->](1.,0) to  [out=-90,in=180] (2.75,-1.5) to  [out=0,in=-90] (4.5,0);
 \draw [fill] (1.,0) circle [radius=0.05];
\draw[thick,<-](1.5,0) to  [out=-90,in=180] (1.75,-0.75) to  [out=0,in=-90] (2.5,0);
\draw [fill] (2.5,0) circle [radius=0.05];
\draw[thick,<-](3,0) to  [out=-90,in=180] (3.5,-0.75) to  [out=0,in=-90] (4,0);
\draw [fill] (4.,0) circle [radius=0.05];

  \node[above] at (1.3,0) {$NP^\bot$};
\node[above] at (2.8,0) {$NP$};
\node[above ] at (4.25,0) {$S^\bot$};
\end{tikzpicture}
$$
We cut the result with $JOHN\_LOVES\_MADLY$, which yields the following cowordism.
$$
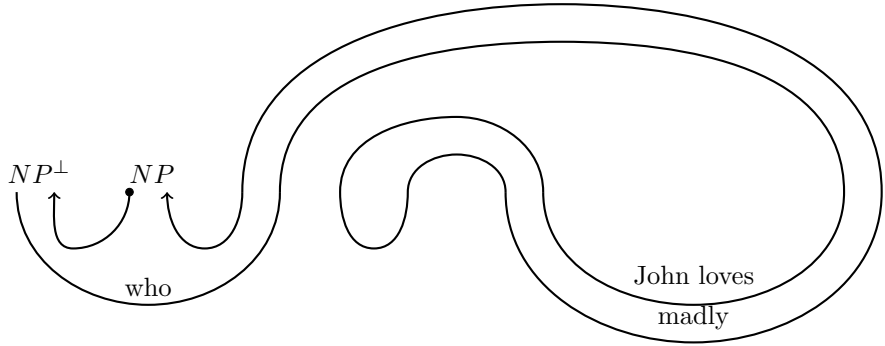
\begin{tikzpicture}
\draw[thick,-]
(8,0) to  [out=90,in=0] (6.85,1) to  [out=180,in=90]
(5.3,0) to  [out=-90,in=180] (5.75,-0.75) to  [out=0,in=-90] (6.2,0)
to  [out=90,in=180] (6.85,0.5) to  [out=0,in=90] (7.5,0);

\node[above] at (2.75,-1.5) {$\rm{who}$};

 \draw[thick,-](1.,0) to  [out=-90,in=180] (2.75,-1.5) to  [out=0,in=-90] (4.5,0);
\draw[thick,<-](1.5,0) to  [out=-90,in=180] (1.75,-0.75) to  [out=0,in=-90] (2.5,0);
\draw [fill] (2.5,0) circle [radius=0.05];
\draw[thick,<-](3,0) to  [out=-90,in=180] (3.5,-0.75) to  [out=0,in=-90] (4,0);

  \node[above] at (1.3,0) {$NP^\bot$};
\node[above] at (2.8,0) {$NP$};


\draw[thick,-](7.5,0) to  [out=-90,in=180] (10,-2) to  [out=0,in=-90] (12.5,0)
to  [out=90,in=0] (8.25,2.5) to  [out=180,in=90](4,0);

\draw[thick,-](8,0) to  [out=-90,in=180] (10,-1.5) to  [out=0,in=-90] (12,0)
to  [out=90,in=0] (8.25,2) to  [out=180,in=90](4.5,0);
;
\node[above] at (10,-1.35) {$\rm{John}\mbox{ }\rm{loves}$};
\node[above] at (10,-1.95) {$\rm{madly}$};

\end{tikzpicture}
$$
The above  has a rather simple shape.
$$
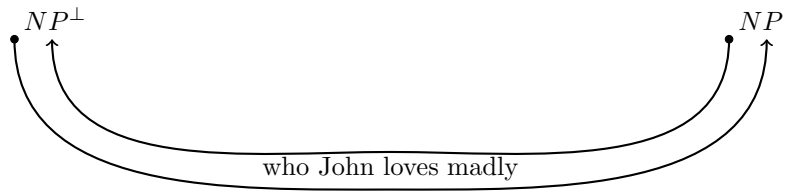
\begin{tikzpicture}
 \draw[thick,->](0,0) to  [out=-90,in=180] (5,-2) to  [out=0,in=-90] (10,0);
\draw[thick,<-](0.5,0) to  [out=-90,in=180] (5,-1.5) to  [out=0,in=-90] (9.5,0);
\draw [fill] (0,0) circle [radius=0.05];
\draw [fill] (9.5,0) circle [radius=0.05];
\node[above right ] at (9.5,0) {$NP$};
\node[above right ] at (0,0) {$NP^\bot$};
\node[above] at (5,-2) {$\mbox{who John loves madly}$};
\end{tikzpicture}
$$
Now it must be clear that after remaining steps in the derivation we obtain (\ref{toy sentence}), as desired.

\begin{remark}
    We deliberately chose a rather non-optimal derivation in order to show  symmetry transformations.  The same string could be generated in the same grammar with fewer Exchange rules. Also, we chose types for lexicon entries which appear as translations of linear implicational types. If, on the other hand, we had, as the type for the transitive verb (i.e. $LOVES$), the sequent $\vdash NP^\bot,S,NP^\bot$,  the analogous derivation would be shorter as well.
\end{remark}

\section{Encoding multiple context-free grammars}
In this section,  as an example, we establish a relationship between LLG and {\it multiple context-free grammars}.

\subsection{Multiple context-free grammars}
 Multiple context-free grammars were introduced in \cite{Seki}. We follow (with minor variations in notation) the presentation in \cite{Kanazawa}.

\bd\label{MCFG def} A {\it multiple context free grammar (MCFG)} $G$ is a tuple $G=(N,T,S,P)$ where
\begin{itemize}
  \item $N$ is a finite  alphabet of nonzero arity predicate symbols  called {\it nonterminal symbols} or {\it nonterminals};
  \item $T$ is a finite alphabet of {\it terminal symbols} or {\it terminals};
  \item $S\in N$, the {\it start symbol},  is unary;
  \item $P$ is a finite set of sequents, called {\it productions},  of the form
  \be\label{production}
   B_1(x^1_1,\ldots,x^1_{k_1}),\ldots,B_n(x^n_1,\ldots,x^n_{k_n})\vdash A(s_1,\ldots,s_k),
  \ee
  where
  \begin{enumerate}[(i)]
    \item $n\geq0$ and $A,B_1,\ldots, B_n$ are nonterminals  with arities $k,k_1,\ldots,k_n$ respectively;
  \item $\{x^j_i\}$ are pairwise distinct variables not from $T$;
  \item $s_1,\ldots, s_k$ are  words built of terminals and $\{x_i^j\}$;
  \item each of the variables $x^j_i$ occurs exactly once in exactly one of the words $s_1,\ldots s_k$.
  \end{enumerate}
 \end{itemize}
\ed
\smallskip

{\bf Remark} Productions are often written in the opposite order in the literature; with $A$ on the left and $B_1,\ldots,B_n$ on the right.

Also, our ``non-erasing'' condition (iv) in the definition of an MCFG, namely, that all  $x_i^j$ occurring on the left occur exactly once on the right, is too strong compared with  original definitions in \cite{Seki}, \cite{Kanazawa}. Usually it is required only that each $x_i^j$ should occur at most once on the right. However, it is known \cite{Seki} that adding the non-erasing condition does not change the expressive power of MCFG, in the sense that the class of generated languages (see below) remains the same.

\bd
The set of {\it  predicate formulas derivable in $G$} is defined by the following induction.
\begin{enumerate}[(i)]
 \item {
 If a production $\vdash A(s_1,\ldots,s_k)$ is in $P$, then  $A(s_1,\ldots,s_k)$ is derivable.
 }
 \item
 {
 For every production {(\ref{production})} in $P$,
 if
 \begin{itemize}
   \item  $B_1(s^1_1,\ldots,s^1_{k_1}),\ldots,B_n(s^n_1,\ldots,s^n_{k_n})$
are derivable,
   \item $t_m$ is  the result of substituting the word $s^j_i$ for every variable $x^j_i$ in $s_m$,  for  $m=1,\ldots,k$,
 \end{itemize}
 then the formula $A(t_1,\ldots,t_k)$ is derivable.
 }
\end{enumerate}
\ed
\begin{remark}
  Technically, condition (i) above is a particular case of (ii) when $n=0$.
\end{remark}

We use notation $\vdash_G F$ to express that $F$ is derivable in $G$.

\bd
The {\it language generated } by an MCFG $G$ is the set of words $w$ for which $\vdash_G S(w)$.

{\it Multiple context-free language} is a language generated by some MCFG.
\ed

When all predicate symbols in $N$ are unary, the above definition reduces to the more familiar case of a {\it context free grammar} (CFG).

\subsection{Cowordism representation of MCFG}
\subsubsection{Productions as cowordisms}\label{MCFG -> cowordisms}
\bd\label{predicate representation}
Let $C$ be a predicate symbol of arity $k>0$.

We define the {\it carrier} $[ C]$ of $C$ as the boundary with cardinality $2k$ and the set of left endpoints
\be\label{MFCG graph(A)}
[C]_l=\{l,3,\ldots,2k-1\}.
\ee
The {\it pattern} $Pat(C)$ of $C$ is the cobordism $$Pat(C):{\bf 1}\to [C]$$
given by edges $$(2i-1,\epsilon,2i),\quad i=1,\ldots,k.$$

If
\be\label{just a formula}
C(s_1,\ldots,s_k),
\ee
is a predicate formula,
 where $s_1,\ldots,s_k$ are words over some alphabet $T$, then
 the cowordism
 $$[C(s_1,\ldots,s_k)]:{\bf 1}\to [C]$$
 over ${T}$
 {\it representing formula }(\ref{just a formula}) is
 given by the set of edges
 $$(2i-1,s_i,2i),\quad i=1,\ldots,k.$$
\ed

The above definition allows representing predicate formulas with variables as well. If some of $s_1,\ldots,s_k$ in (\ref{just a formula}) are not words, but variables from
some set $V$, we simply consider the set $T\cup V$ as a new alphabet.

Now assume that we are given an MCFG $G=(N,T,S,P)$.

We are going to represent productions of $G$ as cowordisms as well.
\bd\label{production representation}
Given a production $p$ of form (\ref{production}), a  cowordism
  $$[p]:{[B_1]}\otimes\ldots\otimes [B_n]\to [ A]$$
  {\it represents} $p$  if the following equation holds
\be\label{equation for MCFG}
[p]\circ([B_1(x_1^1,\ldots,x^1_{k_1})]\otimes\ldots\otimes [{B_n}(x^n_1,\ldots,x^n_{k_n})])=[{A}(s_1,\ldots,s_k)].
\ee
\ed

Before giving a general construction of $[ p]$ we consider a simple example.
\begin{ex}
Let  $T=\{a,b\}$ be the alphabet of terminals.

Consider nonterminals $P$, $Q$ and $S$ of arities 2, 2 and 1 respectively, and the following six productions
$$\vdash P(\epsilon,\epsilon)~(1),\quad\vdash Q(\epsilon,\epsilon)~(2),\quad P(x,y)\vdash P(xa,yb)~(3),$$
$$Q(z,t)\vdash Q(za,ta)~(4),\quad Q(z,t)\vdash Q(zb,tb)~(5),\quad
Q(z,t),P(x,y)\vdash S(zxty)~(6).$$

It is easy to see that the constructed MCFG generates the language
\be\label{MCFL}
\{wa^nwb^n|w\in T^*, n\geq 0\}.
\ee

Now we construct cowordisms representing $G$.
We will omit square brackets and use the same notation for a predicate symbol or a production and  for its graphical representation.

Six cowordisms representing the productions are shown on the picture below (for better readability, we label vertices with corresponding variables,  the subscripts $l,r$ denoting  left and right endpoints respectively)
$$
 \tikz[scale=.7]{
 \node[right] at(-2.25,.75) {$(1)$};


\node[right] at(-4,0) {$x_l$};
 \node[right] at(-4,.5) {$x_r$};
 \node[right] at(-4,1) {$y_l$};
\node[right] at(-4,1.5) {$y_r$};

\draw[dashed,-](-3.25,1.5)--(-3,1.5)--(-3,0)--(-3.25,0);
\node[right]at (-3,.75){$P$};

\draw [fill] (-4,0) circle [radius=0.05];
\draw [fill] (-4,1) circle [radius=0.05];

\draw[thick,->](-4,0)
to[out=180, in=-90] (-4.5,.25)
to[out=90, in=180] (-4,.5);

\draw[thick,->](-4,1)
to[out=180, in=-90] (-4.5,1.25)
to[out=90, in=180] (-4,1.5);

\node[right] at(2,.75) {$(2)$};

\draw[dashed,-](1,1.5)--(1.25,1.5)--(1.25,0)--(1,0);
\node[right]at (1.25,.75){$Q$};


\node[right] at(.25,0) {$z_l$};
 \node[right] at(.25,.5) {$z_r$};
 \node[right] at(.25,1) {$t_l$};
\node[right] at(.25,1.5) {$t_r$};

\draw [fill] (.25,0) circle [radius=0.05];
\draw [fill] (.25,1) circle [radius=0.05];

\draw[thick,->](.25,0)
to[out=180, in=-90] (-.25,.25)
to[out=90, in=180] (.25,.5);

\draw[thick,->](.25,1)
to[out=180, in=-90] (-.25,1.25)
to[out=90, in=180] (.25,1.5);

\node[right] at(9.25,.75) {$(3)$};

\draw[dashed,-](8.25,1.5)--(8.5,1.5)--(8.5,0)--(8.25,0);
\node[right]at (8.5,.75){$P$};

\draw[dashed,-](4.75,1.5)--(4.5,1.5)--(4.5,0)--(4.75,0);
\node[left]at (4.5,.75){$P$};

\node[left] at(5.5,0) {$x_l$};
 \node[left] at(5.5,.5) {$x_r$};
 \node[left] at(5.5,1) {$y_l$};
\node[left] at(5.5,1.5) {$y_r$};

\node[right] at(7.5,0) {$x_l$};
 \node[right] at(7.5,.5) {$x_r$};
 \node[right] at(7.5,1) {$y_l$};
\node[right] at(7.5,1.5) {$y_r$};

\draw [fill] (5.5,.5) circle [radius=0.05];
\draw [fill] (5.5,1.5) circle [radius=0.05];
\draw [fill] (7.5,0) circle [radius=0.05];
\draw [fill] (7.5,1) circle [radius=0.05];

\node[above] at(6.5,1.5) {$b$};
\node[above] at(6.5,0.5) {$a$};

\draw[thick,<-](5.5,0)--(7.5,0);
\draw[thick,->](5.5,.5)--(7.5,.5);
\draw[thick,<-](5.5,1)--(7.5,1);
\draw[thick,->](5.5,1.5)--(7.5,1.5);

}.
$$
$$
 \tikz[scale=.7]{

\node[right] at(3.75,.75) {$(4)$};

\draw[dashed,-](-.75,1.5)--(-1,1.5)--(-1,0)--(-.75,0);
\node[left]at (-1,.75){$Q$};

\draw[dashed,-](2.75,1.5)--(3,1.5)--(3,0)--(2.75,0);
\node[right]at (3,.75){$Q$};

\node[left] at(0,0) {$z_l$};
 \node[left] at(0,.5) {$z_r$};
 \node[left] at(0,1) {$t_l$};
\node[left] at(0,1.5) {$t_r$};

\node[right] at(2,0) {$z_l$};
 \node[right] at(2,.5) {$z_r$};
 \node[right] at(2,1) {$t_l$};
\node[right] at(2,1.5) {$t_r$};

\draw [fill] (0,.5) circle [radius=0.05];
\draw [fill] (0,1.5) circle [radius=0.05];
\draw [fill] (2,0) circle [radius=0.05];
\draw [fill] (2,1) circle [radius=0.05];

\node[above] at(1,1.5) {$a$};
\node[above] at(1,0.5) {$a$};

\draw[thick,<-](0,0)--(2,0);
\draw[thick,->](0,.5)--(2,.5);
\draw[thick,<-](0,1)--(2,1);
\draw[thick,->](0,1.5)--(2,1.5);

\node[right] at(11,.75) {$(5)$};

\draw[dashed,-](6.5,1.5)--(6.25,1.5)--(6.25,0)--(6.5,0);
\node[left]at (6.25,.75){$Q$};

\draw[dashed,-](10,1.5)--(10.25,1.5)--(10.25,0)--(10,0);
\node[right]at (10.25,.75){$Q$};

\node[left] at(7.25,0) {$z_l$};
 \node[left] at(7.25,.5) {$z_r$};
 \node[left] at(7.25,1) {$t_l$};
\node[left] at(7.25,1.5) {$t_r$};

\node[right] at(9.25,0) {$z_l$};
 \node[right] at(9.25,.5) {$z_r$};
 \node[right] at(9.25,1) {$t_l$};
\node[right] at(9.25,1.5) {$t_r$};

\draw [fill] (7.25,.5) circle [radius=0.05];
\draw [fill] (7.25,1.5) circle [radius=0.05];
\draw [fill] (9.25,0) circle [radius=0.05];
\draw [fill] (9.25,1) circle [radius=0.05];

\node[above] at(8.25,1.5) {$b$};
\node[above] at(8.25,0.5) {$b$};

\draw[thick,<-](7.25,0)--(9.25,0);
\draw[thick,->](7.25,.5)--(9.25,.5);
\draw[thick,<-](7.25,1)--(9.25,1);
\draw[thick,->](7.25,1.5)--(9.25,1.5);

}.
$$
$$
 \tikz[scale=.7]{

\node[right] at(4.5,2) {$(6)$};

\node[left] at(0,0) {$z_l$};
 \node[left] at(0,.5) {$z_r$};
 \node[left] at(0,1) {$t_l$};
\node[left] at(0,1.5) {$t_r$};
 \node[left] at(0,2.5) {$x_l$};
\node[left] at(0,3) {$x_r$};
 \node[left] at(0,3.5) {$y_l$};
\node[left] at(0,4) {$y_r$};

\draw[thick,<-](0,0)--(.5,0);
\draw[thick,-](0,.5)--(.5,.5);
\draw[thick,<-](0,1)--(.5,1);
\draw[thick,-](0,1.5)--(.5,1.5);
\draw[thick,<-](0,2.5)--(.5,2.5);
\draw[thick,-](0,3)--(.5,3);
\draw[thick,<-](0,3.5)--(.5,3.5);
\draw[thick,-](0,4)--(.5,4);

\draw[dashed,-](-.75,4)--(-1,4)--(-1,2.5)--(-.75,2.5);

\draw[dashed,-](-.75,1.5)--(-1,1.5)--(-1,0)--(-.75,0);

\node[left] at(-1,.75) {$Q$};
\node[left] at(-1,3.25) {$P$};

\node[right] at(3.75,2) {$S$};

\draw[thick,->](.5,4)
to[out=0,in=90](3.,2.5)
to[out=-90,in=180](3.75,2.25)
;

\draw[thick,-](.5,0)
to[out=0,in=-90](3.,1.5)
to[out=90,in=180](3.75,1.75)
;

\draw [fill] (3.75,1.75) circle [radius=0.05];

\draw[thick,-](0.5,3)--(1.5,1.5);
\draw[thick,-](0.5,1.5)--(1.5,3);
\draw[thick,-](0.5,2.5)--(1.5,1.);
\draw[thick,-](0.5,1.)--(1.5,2.5);


\draw[thick,-](1.5,1.5)
to[out=0,in=-90](2.,2)
to[out=90,in=0](1.5,2.5)
;

\draw[thick,-](1.5,3)
to[out=180,in=-90](2,3.25)
to[out=90,in=0](1.5,3.5)
--(.5,3.5)
;

\draw[thick,-](0.5,.5)--(1.5,.5)
to[out=180,in=-90](2,.75)
to[out=90,in=0](1.5,1)
;

 \draw [fill] (0,.5) circle [radius=0.05];
\draw [fill] (0,1.5) circle [radius=0.05];
\draw [fill] (0,3) circle [radius=0.05];
\draw [fill] (0,4) circle [radius=0.05];

 }.$$

 It is easy to see that the above cowordisms indeed are solutions of
 the corresponding equations of form (\ref{equation for MCFG}).

 Hopefully, it is  also easy to convince oneself that these six cowordisms  together with symmetry maps generate by composition and tensor product all cowordisms from ${\bf 1}$ to $S$ representing words of (\ref{MCFL}) and no other words.
\end{ex}

The solution of (\ref{equation for MCFG}) is constructed as follows.

Let $V$ be the set of all variables $x_i^j$ occurring in $p$.

Let
$$[B]=[B_1]\otimes\ldots\otimes[B_n],$$
and let $\phi_1,\ldots,\phi_n$ be the embeddings
$$\phi_j:{\bf I}([B_j])\to{\bf I}([B]),\quad i\mapsto i+\sum\limits_{\alpha=1}^{j-1}2k_j,\quad j=1,\ldots,n.$$

%
%


Each variable $y=x_i^j\in V$ corresponds to a pair of opposite polarity endpoints in $[B]$.
We denote
$$l(y)=\phi_j(2i-1),\quad r(y)=\phi_j(2i).$$

Now, each word $s_m$, $m=1,\ldots,k$, on the righthand side of $p$ is a concatenation of the form
$$s_m=w^0_my^1_mw^1_m\ldots y^{\alpha_m}_mw^{\alpha_m}_m,$$
 where all $$w^0_m,\ldots,w^{\alpha_m}_m$$
 are words in the alphabet $T$ (possibly empty), and
 $$y^1_m,\ldots,y^{\alpha_m}_m$$
  are variables from $V$. (With the convention that $\alpha_m$ may equal zero, in which case $s_m=w^0_m$.)

We define the  cowordism  $$[ p]:[B]\to[A]$$ as the regular multiword with the boundary $[A]\otimes[B]^\bot$ with the  edges
$$(2m-1,w^0_m,2k+l(y^1_m)),\quad
(2k+r(y^{i}_m),w^i_m,2k+l(y^{i+1}_m)),$$
$$(2k+r(y^{\alpha_m}_m),w^{\alpha_m}_m,2m),\quad i=1,\ldots \alpha_m-1,\quad m=1,\ldots,k.$$

This solves (\ref{equation for MCFG}).

An immediate consequence of (\ref{equation for MCFG}) is the following.
 \bp\label{production is a cowordism}
 Let
%
$$
 B(s_1^1,\ldots,s_1^{k_1}),\ldots,B(s_n^1,\ldots,s_n^{k_n})
 $$
 be predicate formulas, with arities
 $k_1,\ldots,k_m$ respectively.


 Let $t_m$ be  the result of substituting the word $s^j_i$ for every variable $x^j_i$ in $s_m$,  for  $m=1,\ldots,k$,

 Then
 we have the identity
 $$[p]\circ([B_1(s_1^1,\ldots,s^1_{k_1})]\otimes\ldots\otimes [{B_n}(s^n_1,\ldots,s^n_{k_n})])=[{A}(t_1,\ldots,t_k)].\mbox{ }\Box$$
   \ep

\subsubsection{Context-free cowordism grammars}
Now that we have represented MCFG productions as cowordisms, it will be convenient to reformulate (and slightly generalize) MCFG  in the  language of cowordisms directly.

\bd
A {\it  context-free cowordism grammar (cowordism CFG)} $G$ is a tuple $G=(N,T,P,S)$, where
\begin{itemize}
  \item $N$ is a finite set of  {\it types}, each type $A$ being assigned  a boundary $[A]$, the {\it carrier} of $A$;
  \item $T$ is a finite alphabet of {\it terminal symbols};
  \item $P$ is a finite set of rules of the form
\be\label{cowordism produxtion def}
\sigma:A_1\otimes\ldots\otimes  A_n\to  A,
\ee
Where $A_1,\ldots,A_n,A$ are elements of $N$, and
$$\sigma: [{A_1}]\otimes\ldots\otimes [{A_n}]\to [A].$$
is a cowordism;
  \item $S\in N$ is the {\it sentence type}, with  ${\bf{I}}([S])=2$, $([S])_l=\{1\}$.
\end{itemize}
\ed

Elements of $P$ are called {\it cowordism productions}.

Now, for any type $A\in N$, we will define  the set of {\it type $A$ cowordisms}  {\it  generated by $G$}, or, simply, the {\it cowordism type $A$ of} $G$. We will write $\vdash_G\sigma:A$ to express that $\sigma $ is in the cowordism type $A$ of $G$.

The set is defined by induction.
\begin{itemize}
  \item If  a cowordism production $\sigma:{\bf 1}\to A$ is in $P$, then $\vdash_G\sigma:A$.
  \item If a cowordism production $$\sigma:A_1\otimes\ldots\otimes A_n\to A$$ is in $P$, and $$\vdash_G\tau_i:A_i,\quad i=1,\ldots,n,$$
      then $\vdash_G\sigma\circ(\tau_1\otimes\ldots\otimes\tau_n):A$.
\end{itemize}

Any regular cowordism in the type $S$ of $G$ consists of exactly one edge labeled with a word in $T^*$. It is natural to identify the cowordism type $S$ of $G$ with the language generated by $G$.

\bd
The {\it  language generated by the cowordism CFG} $G$ or, simply, the {\it language of $G$} is the set of words occurring as edge labels in   type $S$ regular cowordisms generated by $G$.
\ed
We denote the language generated by $G$ as $L(G)$.

Now, given an  MCFG $G=(N,T,P,S)$, the cowordism representation described in the preceding section
gives us a cowordism CFG.

 We put $G'=(N,T,P',S)$, where each nonterminal $C\in N$  is assigned its carrier $[C]$ as in Definition \ref{predicate representation}, and $P'$ is the set of cowordism productions representing elements of $P$ in the sense of Definition \ref{production representation}:
$$ P'=\{[p]~|~p\in P\}.$$

It is easy to see that $G$ and $G'$ are equivalent in all conceivable ways.
\bp\label{MFCG2extended}
In notation as above, for any predicate symbol $C\in N$ of arity $k$ and words $s_1,\ldots,s_k\in T^*$
it holds that $\vdash_G C(s_1,\ldots,s_k)$ iff $\vdash_{G'} [C(s_1,\ldots, s_k)]$.

In particular, $L(G)=L(G')$.
\ep
\begin{proof}
  By induction on derivation using Proposition \ref{production is a cowordism}.
\end{proof}

\subsubsection{From cowordism CFG to  MCFG}
Let $G=(N,T,P,S)$ be a cowordism CFG.

For each $A\in N$ and regular cowordism $\sigma:{\bf 1}\to [A]$ such that  $\vdash_G\sigma:A$
let $Pat(\sigma)$ be the pattern of $\sigma$, i.e. the cobordism obtained from $\sigma$ by erasing all labels from edges.

We say that $Pat(\sigma)$ is a {\it possible pattern} of $A$.

We denote the set of  possible patterns of $A$ as $Patt(A)$. Note that this set is finite.

\bd
The  cowordism CFG $G$ is {\it simple}, if for any  type $A\in N$ the set $Patt(A)$ contains at most one element.
\ed
\smallskip

It is easy to see that a cowordism CFG constructed from an  MCFG is simple. Let us show that there is an inverse translation.
\bp
If a language is generated by a simple cowordism CFG, then it is also generated by  MCFG.
\ep
\begin{proof}
Let $G=(N,T,P,S)$ be a cowordism CFG.

If there are types without any possible pattern, they do not participate in generation of the language of $G$ and can be safely erased from the grammar. So,
without loss of generality we may assume that for any $A\in N$ the set $Patt(A)$ is nonempty.

Let  $A\in N$.  Let $\pi_A:{\bf 1}\to [A]$ be the unique regular cobordism in $Patt(A)$.

It follows that the cardinality of $[A]$ is an even number $2k$ and there are $k$ left endpoints and $k$ right endpoints, since $\pi_A$ is a perfect matching.

In particular we can write $$[A]_l=\{i_1,\ldots,i_k\},\quad[A]_r=\{j_1,\ldots,j_k\},$$
$$\pi_A=\{(i_\alpha,\epsilon,j_\alpha)|\mbox{ }\alpha=1,\ldots,k\}.$$

Consider the boundary $[A]'$ with cardinality $2k$ and the set of left endpoints
 $$[A]_l'=\{1,3,\ldots, 2k-1\}.$$

We have cobordisms
$$\rho_A:[A]\to[A]',\quad\tau_{A}:[A]'\to[A]$$
defined, respectively by the sets of edges
$$(2\alpha-1,\epsilon,4k-i_\alpha+1), \quad(4k-j_\alpha+1,\epsilon,2\alpha),\mbox{ }\alpha=1\ldots,k$$
and
$$(i_\alpha,\epsilon,4k-2\alpha+2), \quad(4k-2\alpha,\epsilon,j_\alpha),\mbox{ }\alpha=1\ldots,k,$$
satisfying $$\tau_A\circ\rho_A=\id_{[A]},\quad\rho_A\circ\tau_A=\id_{[A']}.$$

It follows that $[A]$ and $[A]'$ are isomorphic, and the cobordism
$$\rho_A\circ\pi_A:{\bf 1}\to[A]'$$ has edges
\be\label{good pattern}
(2\alpha-1,\epsilon,2\alpha),\quad\alpha=1,\ldots,k.
\ee
It follows that we  get an equivalent cowordism CFG by assigning to each $A\in N$ the carrier $[A]'$ and replacing, in each cowordism production of form (\ref{cowordism produxtion def}), the cowordism $\sigma$ with
$$\rho_A\circ\sigma\circ(\tau_{A_1}\otimes\ldots\otimes\tau_{A_n}): [{A_1}]'\otimes\ldots\otimes [{A_n}]'\to [A]'.$$
Moreover, in the obtained grammar, the pattern of $A$ will be precisely (\ref{good pattern}).

So, without loss of generality we can assume  that $[A]=[A]'$, and the pattern of $A$ is (\ref{good pattern}).

Then for each $A\in N$ we introduce a predicate symbol $A$ of arity $k$, where $2k$ is the cardinality of $[A]$.

Let $P_0\subseteq P$ be the set of  cowordism productions that participate in generation of $L(G)$.
All cowordisms occurring in $P_0$ are regular, because elements of $L(G)$ are defined by regular cowordisms.

For each element $p\in P_0$ of form (\ref{cowordism produxtion def}) we write an
MCFG production $p'$ as follows.

For each $j=1,\ldots, n$, let  the cardinality of $[A_j]$ be $2k_j$.

Let $V$ be the set of variables $x_i^j$, $i=1,\ldots,k_j$, $j=1,\ldots,n$.

For each $j=1,\ldots,n$ we define the regular cowordism over $T\cup V$
$$\pi_j:{\bf 1}\to [A_j]$$ by the set of edges
$$(2i-1,x_i^j,2i),\quad i=1,\ldots,k_j.$$
Then we have $Pat(\pi_j)=Pat(A_j)$ and, in fact,
$$\pi_j=[A_j(x_1^j,\ldots,x_{k_j}^j)]$$
in the sense of Definition \ref{predicate representation}.

It follows that the pattern of $$\tau=\sigma\circ(\pi_1\otimes\ldots\otimes\pi_n)$$ coincides with $Pat(A)$,
hence $\tau$ represents, in the sense of Definition \ref{predicate representation} some formula $$A(s_1,\ldots,s_k),$$ where $s_1,\ldots,s_k\in T\cup V$ and the cardinality of $A$ is $2k$.

We write $p'$ as
$$A_1(x^1_1,\ldots,x^1_{k_1}),\ldots,A_n(x^n_1,\ldots,x^n_{k_n})\vdash A(s_1,\ldots,s_k).$$,

Thus we obtain an MCFG $G'=(N,T,P',S)$, where
$$P'=\{p'|~p\in P_0\}.$$

By induction on derivations similar to Proposition \ref{MFCG2extended}, we establish that $G'$ generates the same language as $G$.
\end{proof}

Now we generalize the above to arbitrary cowordism CFG $G$.

\bl
For any cowordism CFG $G$ there exists a simple cowordism CFG $G'$ generating the same language.
\el
\begin{proof}
  Let $G=(N,T,P,S)$.
Without loss of generality we may assume that $L(G)$ is nonempty, otherwise the statement is obvious. (We can take $G'$ with no productions at all.)

We construct a new cowordism CFG $G'$ as follows.

For any  type $A\in N$   and any   possible pattern $\pi$ of $A$
we introduce a new symbol $(A,\pi)$.

We define the set $N'$ of types of $G'$ as
$$N'=\{(A,\pi)|\mbox{ }A\in N,\pi\in Patt(A)\}.$$
Interpretation of types is given by
$$[A,\pi]=[ A].$$

For any cowordism production $$\sigma:A_1\otimes\ldots\otimes A_n\to A$$ of $G$ we consider all possible cowordism productions  of the form
\be\label{new lexicon}
\sigma':(A_1,\pi_1)\otimes\ldots\otimes(A_n,\pi_n)\to(A,\pi),
\ee
where
$$\pi_i\in Patt(A_i),\mbox{ } i=1,\ldots, n,$$
 and $\pi\in Patt(A)$ is constructed as the composition
$$\tau=Pat(\sigma)\circ(\pi_1\otimes\ldots\otimes\pi_n).$$

The  set $P'$ of productions for $G'$ consists of all cowordism productions of form (\ref{new lexicon}). Again, there are only finitely many of them.

Since the set $L(G)$ is assumed nonempty,  the set $Patt(S)$ is a singleton. We denote $S'=(S,e)$, where $e$ is the only element of $Patt(S)$.

We define $G'$ as $G'=(N',T,L',S')$.

It is immediate that $G'$ is simple and generates the same  language as $G$.
\end{proof}

Combining the above with the
preceding lemma, we obtain the following.
\bl\label{MCFG=cowordisms}
A language is multiple context-free iff it is generated by a cowordism CFG. $\Box$
\el

\subsection{Embedding in LLG}
\subsubsection{From cowordism CFG to LLG}\label{encoding MCFG}
Any cowordism CFG (hence any MCFG) gives rise to an LLG.

Let
$G=(N,T,P,S)$ be a cowordism CFG.

We consider
each type $A\in N$ as a positive literal, whose interpretation is  its carrier $[A]$.

Then we define an LLG  $G'=(N,T,P',S)$
where the lexicon $P'$ consists of expressions
$$\frac{\ulcorner\sigma\urcorner}{\vdash\Gamma}$$
such that
\be\label{lex from MCFG}
\Gamma=A_1^\bot,\ldots, A_n^\bot,A
\ee
and there is a  production $\sigma:A_1\otimes\ldots\otimes A_n\to A$  in $P$.

Note that the grammar is well defined. The sequent $\Gamma$ is interpreted as the boundary
$[\Gamma]=[{A_1}]^\bot\wp\ldots\wp[{A_n}]^\bot\wp[ A]$, and the name $\ulcorner\sigma\urcorner$, by definition, is a cowordism from the empty boundary ${\bf 1}$ to $[ \Gamma]$.

  It is immediate
  from Proposition \ref{production is a cowordism}
  that the language generated by $G$ is  a subset of the language generated by $ G'$.

Let us prove the opposite inclusion.

Note that the constructed LLG $G'$ has a particularly simple lexicon: formulas in the lexicon do not contain any logical  connective. Let us call such lexicons {\it logic-free}.
 \bp\label{no logical rules}
 Let $H$ be an LLG with a logic-free lexicon.

 Let $\Gamma$ be a sequent not containing logical connectives, and let
 $$\cfrac{\sigma}{\vdash\Gamma}$$ be a typing judgement derivable in $H$.

 Then $\sigma$ is an interpretation of a derivation $\pi$ of $\Gamma$ from axioms of $H$ such that $\pi$ does not contain any logical rules.
 \ep
\begin{proof}
  Let $\pi$ be a derivation of $\Gamma$ from axioms of $H$ whose interpretation is $\sigma$.

  Since interpretation is invariant under cut-elimination, and cut-elimination for  {\bf LL} always terminates, we may assume that no step of the cut-elimination algorithm can be applied to transform $\pi$.

 Assume that there is an occurrence of a logical rule in $\pi$.

  Any logical rule introduces a logical connective. It follows that there is an application of the Cut rule afterwards that eliminates this connective, because $\Gamma$ is connective-free. This cut is on a compound formula, so none of its premises is an axiom. But then the cut-elimination algorithm applies, and $\pi$ can be transformed.
\end{proof}

Coming back to our cowordism CFG $G$ and corresponding LLG $G'$, we have the following.
\bc
Let
$$\Gamma=A_1^\bot,\ldots,A_n^\bot,A
$$
where $A_1,\ldots,A_n,A\in N$, and let
$$
\cfrac{\sigma}{\vdash\Gamma}
$$
be a typing judgement derivable in $G'$,

If
$$\sigma_i:{\bf 1}\to A_i$$  are cowordisms such that
$$\vdash_G\sigma_i:A_i,\quad i=1,\ldots,n,$$
then  $$\vdash_G\langle\sigma_1\otimes\ldots\otimes\sigma_n,\sigma\rangle_{A_1\otimes\ldots\otimes A_n}:A.$$
\ec
\begin{proof}
  By induction on a $G'$ derivation without logical rules.
\end{proof}
\bc
The language of $G'$ coincides with the language of $G$. $\Box$
\ec

Thus we have established the following.
\bl\label{extended->LLG}
If a language is generated by a cowordism CFG, then it is generated by an LLG. $\Box$
\el
\smallskip

The converse of Lemma \ref{extended->LLG} is false, as we will show later.

It is not hard to see though that it holds when the lexicon of an LLG is logic free.

In fact, we can be less restrictive and require that the lexicon has no occurrence of the $\otimes$-connective.
We will call such lexicons {\it $\otimes$-free}.

\subsubsection{From  $\otimes$-free lexicon to cowordism CFG }
An $\otimes$-free lexicon is essentially a logic free lexicon.
 \bp\label{no logical connectives}
  Any $\otimes$-free LLG  $G$ can
  be replaced with an LLG $G'$ whose lexicon is logic free and which generates same types as $G$.
 \ep
 \begin{proof}
   Immediate from Corollary \ref{ts-free=par-free}.
 \end{proof}

So it is sufficient to consider an LLG $G=(T,N,Lex,S)$  whose lexicon $Lex$ is logic free.

For each axiom $$\frac{\sigma}{\vdash A_1,\ldots,A_n}\in Lex,$$ where $A_1,\ldots,A_n$ are literals,
consider $n$ cowordisms
\be\label{new productions}
\sigma_i:A_{i+1}^\bot\otimes\ldots\otimes A_{n}^\bot\otimes A_1^\bot\otimes\ldots\otimes A_{i-1}^\bot \to A_i,\quad i=1,\ldots,n
\ee
obtained from $\sigma$ using correspondence (\ref{*-autonomy}) and symmetry transformations.

We take the set of literals
$$N'=N\cup N^\bot$$
as the set of types for a cowordism CFG
and the set $P$ of all cowordisms of form (\ref{new productions}) as the set of cowordism productions.

This gives us a cowordism CFG $G'$.

By an easy induction on derivation $G'$ we establish that $L(G')\subseteq L(G)$.

For the opposite inclusion we have the following.
\bp
Let $$\Gamma=A_1,\ldots,A_n$$
where $A_1,\ldots,A_n$ are literals, and let
$$
\frac{\sigma}{\vdash\Gamma}
$$
be a typing judgement derivable in $G$,

If
$$\sigma_{i}:{\bf 1}\to A_{i}$$  are cowordisms such that
$$\vdash_{G'}\sigma_{i}:A_{i},\quad i=1,\ldots,n-1$$
then  $$\vdash_{G'}\langle\sigma_{1}\otimes\ldots\otimes\sigma_{n-1},\sigma\rangle_{A_{1}\otimes\ldots\otimes A_{n-1}}:A_{n}.$$
\ep
\begin{proof}
  By induction on a $G$ derivation without logical rules.
\end{proof}
\bc
The language of $G'$ coincides with the language of $G$. $\Box$
\ec

Summing up, we have established the following.
\bl
A language $L$ is generated by a cowordism CFG iff $L$ is generated by an ${\otimes}$-free LLG. $\Box$
\el
\smallskip

Putting the above and Lemma  \ref{MCFG=cowordisms}  together we obtain the following.
\bt\label{tens. free lex is MCFG}
A language is multiple context-free iff it is generated by an LLG with a $\otimes$-free lexicon.
$\Box$
\et

\section{Representing abstract categorial grammars}
Abstract categorial grammars (ACG) were introduced in \cite{deGroote}. They are based on the purely implicational fragment of linear logic, and {\bf LL} grammars of this paper can be seen as a representation and extension of ACG.

In this section we assume that the reader is familiar with basic notions of  $\lambda$-calculus, see \cite{Barendregt} for a reference. We use \cite{HylandDePaiva} as a reference for syntax and semantics of linear $\lambda$-calculus and intuitionistic linear logic. We note though that we consider only the simplest, implicational fragment, while definitions and results  in \cite{HylandDePaiva}  are formulated for the full system. In fact, theorems of \cite{HylandDePaiva} that we cite are rather straightforward in the purely implicational case.

\subsection{Implicational logic and Linear $\lambda$-calculus}
\subsubsection{Linear $\lambda$-calculus}
Given a  set $N$ of {\it literals}, the set $Fm_\multimap(N)$ of {\it implicational linear logic ({\bf ILL})} formulas (over $N$), also denoted as  $Tp(N)$  and called the set of {\it linear implicational types} (over $N$),   is defined by induction.
\begin{itemize}
\item Any $A\in N$ is in $Fm_\multimap(N)$;
\item if $A,B\in Fm_\multimap(N)$, then $(A\multimap B)\in Fm_\multimap(N)$.
\end{itemize}

{\bf ILL} formulas, seen as  types, are used for typing {\it linear $\lambda$-terms},
which are $\lambda$-terms where each  variable occurs exactly once.

More accurately,   given a  set $X$ of {\it variables} and a   set $C$ of  {\it constants}, with $C\cap X=\emptyset$, the set $\Lambda(X,C)$ of {\it linear $\lambda$-terms} is defined by the following.
\begin{itemize}
\item Any $a\in X\cup C$ is in $\Lambda(X,C)$;
\item if $t,s\in\Lambda(X,C)$ are linear $\lambda$-terms whose sets of free variables are disjoint then $(t\cdot s)\in\Lambda(X,C)$;
\item if $t\in\Lambda(X,C)$, and $x\in X$ occurs freely in $t$ exactly once then $(\lambda x.t)\in \Lambda(X,C)$.
\end{itemize}

We use  common notational conventions such as omitting dots and outermost brackets and writing  iterated applications as
\be\label{iterated application}
(tsk)=(ts)k.
\ee

A {\it typing judgement} is a sequent of the form
$$x_1:A_1,\ldots,x_n:A_n\vdash t:A,$$
where $x_1,\ldots x_n\in X$ are pairwise distinct ($n$ may be zero), $t\in\Lambda(X,C)$, and $A_1,\ldots,A_n,A\in Tp(N)$.

 Typing judgements are derived from the following type inference rules.
$${x:A\vdash x:A}~(\mbox{Id}),$$
$$\frac{\Gamma\vdash t:A\multimap B\quad\Delta\vdash s:A}{\Gamma,\Delta\vdash (ts):B}~(\multimap\rm{E}), \quad\frac{\Gamma,x:A,\Delta\vdash t:B}{\Gamma,\Delta\vdash (\lambda x.t):A\multimap B}~(\multimap\rm{I}).$$

An {\it {\bf ILL} sequent} is an expression of the form $A_1,\ldots,A_n\vdash A$, where $A_1,\ldots,A_n$ and $A$ are {\bf ILL} formulas.

{\it Natural deduction} formulation of {\bf ILL} is given by rules of linear $\lambda$-calculus with all terms erased.

Thus, linear $\lambda$-calculus, from a proof-theoretical pot of view, is a way to label ${\bf ILL}$ proofs. The following is a ``linear'' case of {\it Curry-Howard isomorphism}.
\bp[\cite{HylandDePaiva}]\label{Curry-Howard} There is a one-to one correspondence between {\bf ILL} natural deduction proofs and derivable typing judgements. $\Box$
\ep

\bd
A {\it linear signature}, or, simply, a {\it signature}, $\Sigma$ is a triple $\Sigma=(N,C,\mathfrak{T})$, where $N$ is a finite set of atomic types, $C$ is a finite set of constants and $\mathfrak{T}$ is a function assigning to each constant $c\in C$ a linear implicational type $\mathfrak{T}(c)\in Tp(N)$.
\ed

Typing judgements of the form
\be\label{linear signature axiom}
\vdash c:\mathfrak{T}(c),
\ee
 where $c\in C$, are called {\it signature axioms} of $\Sigma$.

Given a signature  $\Sigma$, we say that a typing judgement is {\it derivable in} $\Sigma$ if it is derivable from axioms of $\Sigma$ by rules of linear $\lambda$-calculus. We write in this case $\Gamma\vdash_\Sigma t:A$.


\subsubsection{Semantics}\label{lambda-calc. semantics}
Let ${\bf C}$ be a symmetric monoidal category, and
$N$ be a set of literals.

An interpretation of linear types (or {\bf ILL} formulas) over $N$  in ${\bf C}$ consists in assigning to each atomic type $A\in N$ an object $[A]\in{\bf C}$. This is extended to all types in $Tp(N)$ by the obvious induction:
$$[A\multimap B]=[A]\multimap [B].$$

The interpretation extends to all  derivable typing judgements (respectively, all {\bf ILL} natural deduction proofs).

To each derivable typing judgement $\sigma$ of the form
$$x_1:A_1,\ldots,x_n;A_m\vdash A$$
(respectively, to an {\bf ILL} natural deduction proof of the corresponding sequent)
we assign a ${\bf C}$-morphism
$$[\sigma]:
[A_1]\otimes\ldots\otimes [A_n]\to [A],$$
if $n>0$, or
$$[\sigma]:{\bf 1}\to [A],$$ if $n=0$, by induction on type inference rules.

If $\sigma$ is $x:A\vdash x:A$  is  the (Id) axiom, then $[\sigma]=\id_{[A]}$.

If $\sigma$ is obtained from a derivable judgement $\sigma'$ by the ($\multimap\rm{I}$) rule, then $[\sigma]$ is obtained from $[\sigma']$ using symmetry and  correspondence (\ref{monoidal closure}).

If $\sigma$ is obtained from  derivable judgements
$$\sigma_1=\Gamma_1\vdash t:A\multimap B,\quad \sigma_2=\Gamma_2\vdash s:A$$
 by the ($\multimap\rm{E}$) rule, then
$$[\sigma]=\ev_{[A],[B]}\circ([\sigma_1]\otimes[\sigma_2])$$
(see (\ref{evaluation}) for notation).

We now specialize the general construction concretely to the category of cowordisms, which is monoidal closed.

Interpretation of rules is shown on the following picture (in the horizontal representation; square brackets are omitted).
$$
 \tikz[scale=.3]{
         \draw[draw=black,fill=gray!10](0,0)rectangle(2,1);\node at(1,.5){$\sigma$};
 \draw[thick,-](2,.5)--(3,.5);

 \node[right] at(3,.5) {$A$};
 \draw[thick,-](-1,.5)--(0,.5);

 \node[left] at(-1,.5) {$\Delta$};

    \draw[draw=black,fill=gray!10](0,-2)rectangle(2,-.5);\node at(1,-1.25){$\tau$};
\draw[thick,-](2,-.8)--(3,-.8);

 \node[right] at(3,-.6) {$A^\bot$};
\draw[thick,-](2,-1.7)--(3,-1.7);
  \node[right] at(3,-1.9) {$B$};

 \draw[thick,-](-1,-1.25)--(0,-1.25);
 \node[left] at(-1,-1.25) {$\Gamma$};

\node at(7,-.5) {$\Longrightarrow$};
\begin{scope}[shift={(5,0)}]
         \draw[draw=black,fill=gray!10](7,0)rectangle(9,1);\node at(8,.5){$\sigma$};
 \draw[thick,-](9,.5)--(10,.5)
 to[out=0,in=90](10.5,-.25)
 to[out=-90,in=0](10,-1)
 ;

 \draw[thick,-](6,.5)--(7,.5);

 \node[left] at(6,.5) {$\Delta$};

         \draw[draw=black,fill=gray!10](7,-2)rectangle(9,-.5);\node at(8,-1.25){$\tau$};
\draw[thick,-](9,-1)--(10 ,-1);

\draw[thick,-](9,-1.5)--(10,-1.5);
  \node[right] at(10,-1.5) {$B$};

 \draw[thick,-](6,-1.25)--(7,-1.25);
 \node[left] at(6,-1.25) {$\Gamma$};

\node[right] at(11,-.5) {$(\multimap E)$};

\end{scope}

 }.$$
 $$
 \tikz[scale=.5]{
 \begin{scope}[shift={(1,0)}]
 \begin{scope}[shift={(.5,0)}]
         \draw[draw=black,fill=gray!10](0,.55)rectangle(2,2.05);\node at(1,1.3){$\sigma$};
 \draw[thick,-](-1,1.8)--(0,1.8);
 \draw[thick,-](-1,1.3)--(0,1.3);
 \draw[thick,-](-1,.8)--(0,.8);

 \node[left] at(-1,2) {$\Delta$};
 \node[left] at(-1,1.3) {$A$};
 \node[left] at(-1,.6) {$\Gamma$};

\draw[thick,-](2,1.3)--(3,1.3);

 \node[right] at(3,1.3) {$B$};
 \end{scope}
 \node at(6,1.3) {$\Longrightarrow$};
 \end{scope}

\begin{scope}[shift={(3,0)}]
         \draw[draw=black,fill=gray!10](8,.25)rectangle(10,1.75);\node at(9,1){$\sigma$};
 \draw[thick,-](6.25,1)--(6.75,1.)--(7.5,1.5)--(8,1.5);
 \draw[thick,-](6.75,1.5)--(7.5,1.)--(8,1.);
 \draw[thick,-](6.25,.5)--(8,.5);
\draw[thick,-](6.75,1.5)
to[out=180,in=-90](6.25,2.)
to[out=90,in=180](6.75,2.25)
--(11,2.25);

 \node[left] at(6.25,1) {$\Delta$};
 \node[right] at(11,2.25) {$A^\bot$};
 \node[left] at(6.25,.3) {$\Gamma$};

\draw[thick,-](10,1)--(11,1);

 \node[right] at(11,1) {$B$};
 \node[left] at(14,1.3) {$(\multimap I)$};
\end{scope}

 }.$$

Now, assume that we are given a signature $\Sigma=(N,C,\mathfrak{T})$.

An interpretation of the signature $\Sigma$  in a monoidal closed category ${\bf C}$
    consists of an interpretation of $Tp(N)$ in $\bf C$ together with a function assigning to each signature axiom (\ref{linear signature axiom}) a morphism
\be\label{semantic signature axiom}
[c]:{\bf 1}\to [\mathfrak{T}(c)].
\ee

The interpretation then extends to all typing judgements derivable in $\Sigma$ by the same rules.

In accordance with notation in (\ref{semantic signature axiom}), we will denote the cowordism interpretation of a derivable typing judgement of the form $\vdash t:F$  as $[t]=[\vdash t:F]$.

(The notation is slightly ambiguous, because the type $F$, in general, is not uniquely determined by the term $t$, but this will not lead to a confusion.)

It might be evident from the geometric representation that the interpretation is invariant, say, under $\beta$-reduction: an introduction rule ($\multimap\rm{I}$) followed by  an elimination ($\multimap\rm{E}$) amounts to partial composition of cowordisms, i.e. to substitution. In any case we have the following.

\bl[\cite{HylandDePaiva}]\label{beta-eta}
Let $\Sigma$ be a signature interpreted in a monoidal closed category $\bf C$.

 If $$\Gamma\vdash t:A,\quad \Gamma\vdash s:A$$ are  typing judgements derivable in $\Sigma$, and the terms $t,s$ are $\beta\eta$-equivalent, then the interpretations coincide,
    $$[\Gamma\vdash t:A]=[\Gamma\vdash s:A];$$
\el

\subsubsection{Sequent calculus formulation}
The {\it sequent} formulation of linear $\lambda$-calculus  is given by the following rules:
$$x:A\vdash x:A~(\rm{Id}),\quad \frac{\Gamma\vdash t:A\quad{x:A,\Delta\vdash s:B}}{\Gamma,\Delta\vdash s[x:=t]: B}~(\rm{Cut}),$$
$$ \frac{\Gamma,x:A,y:B,\Delta\vdash t:C}{\Gamma,y:B,x:A,\Delta\vdash C}~(\rm{Ex}),$$
$$ \frac{\Gamma\vdash x:A\quad{y:B,\Delta\vdash s:C}}{\Gamma,f:A\multimap B,\Delta\vdash s[y:=(fx)]C}~(\multimap\rm{L}),\quad
 \frac{\Gamma,x:A\vdash t:B}{\Gamma\vdash \lambda x.t:A\multimap B}~(\multimap\rm{R}).$$

The {\it sequent calculus} formulation of  ${\bf ILL}$ is obtained from the above by erasing all terms in typing judgements.

 The sequent calculus formulation of {\bf ILL} is cut-free with algorithmic cut-elimination, similarly to ordinary ({\it classical}) ${\bf LL}$.

It is well known that natural deduction and sequent calculus formulations  are equivalent.

\bl[\cite{HylandDePaiva}]
A typing judgement (respectively, an ${\bf ILL}$ sequent) is derivable in natural deduction iff it is derivable in sequent calculus.

Moreover, there exists a translation from natural deduction proofs to sequent calculus proofs and vice versa. $\Box$
\el

Any interpretation of linear types ({\bf ILL} formulas)   in a monoidal closed category gives rise to an interpretation of sequent calculus proofs as well.

Specifically, the Cut is interpreted as partial composition, the Exchange, as a symmetry transformation, while the ($\rm{Id}$) axiom and the $(\multimap\rm{R})$ rule are already present in the natural deduction.

For the sake of illustration let us compute the interpretation of the $(\multimap\rm{L})$ rule. (We specialize to the case of cowordisms below, although the statement is true for any $*$-autonomous category). We will use the usual convention and not distinguish notationally a formula or a proof from its interpretation.

 \bl\label{left introduction ILL}
 If $\sigma$ and $\tau$ are cowordisms representing  proofs of sequents
 $$\Gamma\vdash A,\quad B,\Delta\vdash C$$ respectively,
 then the proof of the sequent $$\Gamma,A\multimap B,\Delta\vdash C$$ obtained from $\sigma$ and $\tau$ by the $(\multimap\rm{L})$ rule is interpreted as the cowordism
 \be\label{cob. left into rule}
 (\llcorner\sigma\lrcorner\wp\tau)\circ\epsilon_{\Gamma,A^\bot,B,\Delta}
 \ee
  (see (\ref{internal cotensor}) for notation), shown in the following picture.
$$
 \tikz[scale=.3]{
         \draw[draw=black,fill=gray!10](0,0)rectangle++(2,1);\node at(1,.5){$\sigma$};
 \draw[thick,-](2,.5)--(3,.5);

 \node[right] at(3,.5) {$A$};
 \draw[thick,-](-1,.5)--(0,.5);
 \node[left] at(-1,.5) {$\Gamma$};
\begin{scope}[shift={(-6,-1.5)}]
         \draw[draw=black,fill=gray!10](6,-.25)rectangle(8,1.25);\node at(7,.5){$\tau$};
\draw[thick,-](8,.5)--(9,.5);

 \node[right] at(9,.5) {$C$};

 \draw[thick,-](5,0)--(6,0);
 \node[left] at(5,0) {$B$};
 \draw[thick,-](5,1)--(6,1);
 \node[left] at(5,1) {$\Delta$};
 \end{scope}

 \begin{scope}[shift={(14,-1.25)}]
 \node[right]at(-8,1.){$\Longrightarrow$};

          \draw[draw=black,fill=gray!10](0,-1)rectangle(2,0);\node at(1,-.5){$\sigma$};
 \draw[thick,-](2,.-.5)--(3,-.5)
 to[out=0,in=-90](3.5,0)
 to[out=90,in=0](2,.5)--(-.5,.5)--(-1.5,1.5)--(-2,1.5);

 \node[left] at(-2,1.5) {$A^\bot$};
 \draw[thick,-](-2,-.5)--(0,.-.5);
 \node[left] at(-2,-.5) {$\Gamma$};

         \draw[draw=black,fill=gray!10](0,1.25)rectangle(2,2.75);\node at(1,2){$\tau$};
\draw[thick,-](2,2)--(3,2);

 \node[right] at(3,2) {$C$};

 \draw[thick,-]
 (-2,.5)--(-1.5,.5)--(-.5,1.5)--(0,1.5);
 \node[left] at(-2,.5) {$B$};
 \draw[thick,-](-2,2.5)--(0,2.5);
 \node[left] at(-2,2.5) {$\Delta$};
 \node[right]at(4,1.){($\multimap\rm{L}$)};
 \end{scope}
 }.$$
 \el
 \begin{proof}
 The $(\multimap\rm{L})$ rule is emulated using Exchange, Cut and the natural deduction rule $(\multimap\rm{E})$ by means of the following derivation:
$$
\cfrac
    {
        \cfrac
            {
                     \cfrac
                        {
                        \begin{array}{c}
                        \\
                        A\multimap B\vdash A\multimap B
                        \end{array}
                        \quad
                        \cfrac
                            {
                            {\sigma}
                            }
                            {
                            \Gamma\vdash A
                            }
                        }
                        {
                            A\multimap B,\Gamma\vdash B
                        }(\multimap E)
            }
            {
              \Gamma,A\multimap B\vdash B
            }(\rm{Ex})
        \quad
            \cfrac{\tau}
            {B,\Delta\vdash B}
    }
    {
    \Gamma,A\multimap B,\Delta\vdash C.
    }(\rm{Cut})
$$

The interpretation of the above proof is computed in the following picture.

$$
 
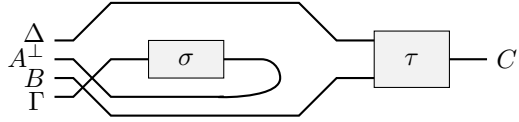
\tikz[scale=.5]{
         \draw[draw=black,fill=gray!10](0,0)rectangle++(2,1);\node at(1,.5){$\sigma$};
 \draw[thick,-](2,.5)--(3,.5)
 to[out=0,in=90](3.5,0)
 to[out=-90,in=180](2,-.5)
 --(-1,-.5)--(-2,.5)--(-2.5,.5);

 \node[left] at(-2.5,.6) {$A^\bot$};
 \draw[thick,-](-2.5,-.5)--(-2,-.5)--(-1,.5)--(0,.5);
 \node[left] at(-2.5,-.6) {$\Gamma$};

         \draw[draw=black,fill=gray!10](6,-.25)rectangle(8,1.25);\node at(7,.5){$\tau$};
\draw[thick,-](8,.5)--(9,.5);

 \node[right] at(9,.5) {$C$};

 \draw[thick,-](-2.5,0)--(-2,0)--(-1,-1)--(4,-1)--(5,0)--(6,0);
 \node[left] at(-2.5,0) {$B$};
 \draw[thick,-](-2.5,1)--(-2,1)--(-1,2)--(4.,2)--(5,1)--(6,1);
 \node[left] at(-2.5,1.2) {$\Delta$};
 }.$$

%
%
%
%
%
%
%

On the other hand, computing the cowordism in (\ref{cob. left into rule})  directly from its expression we get the  alternative picture.
$$
 
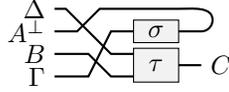
\tikz[scale=.3]{
         \draw[draw=black,fill=gray!10](0,0)rectangle(2,1);\node at(1,.5){$\sigma$};
 \draw[thick,-](2,.5)--(3,.5)
 to[out=0,in=-90](3.5,1)
 to[out=90,in=0](2,1.5)--(-1.5,1.5)--(-2.5,.5)--(-3.5,.5);

 \node[left] at(-3.5,.5) {$A^\bot$};
 \draw[thick,-](-3.5,-1.5)--(-2.,-1.5)--(-1.,.5)--(0,.5);
 \node[left] at(-3.5,-1.5) {$\Gamma$};

         \draw[draw=black,fill=gray!10](0,-1.75)rectangle(2,-.25);\node at(1,-1){$\tau$};
\draw[thick,-](2,-1)--(3.,-1);

 \node[right] at(3,-1) {$C$};

 \draw[thick,-]
 (-3.5,-.5)--(-2,-.5)--(-1,-1.5)--(0,-1.5);
 \node[left] at(-3.5,-.5) {$B$};
 \draw[thick,-](0,-.5)--(-1,-.5)--(-2.5,1)--(-2.5,1)--(-3,1.5)--(-3.5,1.5);
 \node[left] at(-3.5,1.5) {$\Delta$};
 }.$$
Both alternatives represent the same graph as depicted in the formulation of the lemma.
\end{proof}

  Similarly to the case of classical {\bf LL}, interpretations of {\bf ILL} in monoidal closed categories, in particular, in categories of cowordisms are invariant under cut-elimination \cite{HylandDePaiva}.

  The following proposition, similar to Proposition  \ref{alternative derivability}, will be useful.

  \bp\label{alternative derivability ILL}
Let $\Sigma$ be a linear signature and let an interpretation of $\Sigma$ in the category of cowordisms be given.

A typing judgement $\sigma$ of the form ${\Gamma}\vdash t:A $ is derivable in $\Sigma$ iff there exist signature axioms
$$\vdash c_1:\mathfrak{T}(c_1),\ldots,\vdash c_n:\mathfrak{T}(c_n), $$ a term $t'$ and a
typing judgement $\sigma'$ of the form $$x_1:\mathfrak{T}(c_1),\ldots,\vdash x_n:\mathfrak{T}(c_n),\Gamma\vdash t':A$$
(where $x_1,\ldots,x_n$ do not occur in $\sigma$) derivable
in linear $\lambda$-calculus  (not using axioms of $\Sigma$) such that
$$\sigma=\sigma'\circ(\sigma_1\otimes\ldots\otimes\sigma_n),$$
where $\sigma_i=[ c_i]$, $i=1,\ldots,n$, (see (\ref{semantic signature axiom}) for notation).
\ep
\begin{proof}
  By induction on derivation.
\end{proof}

\subsection{Embedding into {\bf LL}}
{\bf ILL} can also be faithfully represented as a fragment of classical {\bf LL}.

Given an alphabet $N$ of literals, the set $Fm_\multimap(N)$ of implicational formulas over $N$ is mapped to the set $Fm(N)$ of classical formulas over $N$ by means of translation (\ref{linear implication}).

Next, an {\bf ILL} sequent $$A_1,\ldots,A_n\vdash A$$ is mapped to the {\bf LL} sequent $$\vdash A_1^\bot,\ldots,A_n^\bot,A.$$

{\bf ILL} sequent proofs translate to {\bf LL} proofs in a straightforward way, with the $(\multimap \rm{R})$ rule corresponding to the $(\wp)$ rule,  the $(\multimap \rm{L})$ rule corresponding to the $(\otimes)$ rule, while Identity, Cut and Exchange of {\bf ILL} translate to Identity, Cut and Exchange of {\bf LL} respectively.

It is very easy to check  that this embedding is {\it conservative} \cite{Schelinx}: a sequent is cut-free derivable in {\bf ILL} iff its translation is cut-free  derivable in {\bf LL}. Moreover,  {\bf ILL} proofs of an {\bf ILL} sequent and {\bf LL} proofs of its translation are in a one-to-one correspondence modulo some inessential permutations of rules. On the semantic side this is reflected in the following.

\bl\label{semantic conservativity}
Let an alphabet $N$ of literals and an interpretation of elements of $N$ as boundaries be given.

Let $A_1,\ldots,A_n,A\in Fm_\multimap(N)$.

There exists an {\bf LL} proof $\pi$ of the sequent
$$\vdash A^\bot_1,\ldots,A_n^\bot,A$$
iff there is an ${\bf ILL}$ proof $\pi'$ of the sequent $$A_1,\ldots,A_n\vdash A$$
such that the interpretation $\sigma$ of $\pi$ equals the name of the interpretation $\sigma'$ of $\pi'$,
$\sigma=\ulcorner\sigma'\urcorner$.
\el
\begin{proof}
  By induction on cut-free derivations using straightforward manipulations with pictures as in Lemma \ref{left introduction ILL}.
\end{proof}

\subsubsection{Representing linear signatures}\label{representing linear signatures}
Thanks to the embedding of {\bf ILL} to {\bf LL}, we can represent linear signatures in {\bf LL} as well.

Let $\Sigma=(N,C,\mathfrak{T})$ be a linear signature and assume that an interpretation of $\Sigma$ in the category of cowordisms over some alphabet $T$ is given.

As previously, we consider  atomic types of $\Sigma$ as literals and this gives us an interpretation of {\bf LL} as well.

We then translate $\Sigma$ to an equivalent cowordism signature $\Sigma'$ as follows.

We define the set  $\Xi$ of  axioms as
$$\Xi=\{\cfrac{[c]}{\vdash\mathfrak{T}(c)}|~c\in C\}.$$
(See (\ref{semantic signature axiom}) for notation.)
In other words, elements  of  $\Xi$ are signature   axioms of $\Sigma$ taken together with their interpretations.


The cowordism signature $\Sigma'$ is defined as $\Sigma'=(N,T,\Xi)$.

\bl\label{conservativity over ACG}
Let $A_1,\ldots,A_n,A\in Tp(N)$, and $$\sigma:{\bf 1}\to [A_1]^\bot\wp\ldots\wp [A_n]^\bot\wp [A]$$ be a cowordism over $T$.

The typing judgement $$\cfrac{\sigma}{{\vdash A_1^\bot,\ldots,A_n^\bot,A}}$$ is derivable in $\Sigma'$ iff
there is some term $t$ such that
\be\label{ACG derivable judgement}
x_1:A_1,\ldots,x_n:A_n\vdash_{\Sigma} t:A
\ee
 and $\sigma$ is the name of the cowordism
$$[x_1:A_1,\ldots,x_n:A_n\vdash t:A]:[A_1]\otimes\ldots\otimes [A_n]\to [A]$$
corresponding to (\ref{ACG derivable judgement}).
\el
\begin{proof}
Proposition \ref{alternative derivability ILL}, Lemma \ref{semantic conservativity} and Proposition \ref{alternative derivability}.
\end{proof}

We say that the cowordism signature $\Sigma'$ is the {\it cowordism representation} of the linear signature $\Sigma$ induced by the given interpretation of types.

\subsection{Abstract categorial grammars}
Given two linear signatures $\Sigma_i=(N_i,C_i,\mathfrak{T}_i)$, $i=1,2$, a {\it map of signatures} $$\phi:\Sigma_1\to\Sigma_2$$ is a
pair $\phi=(F,G)$, where
\begin{itemize}
  \item $F:Tp(\Sigma_1)\to  Tp(\Sigma_2)$ is a function satisfying the homomorphism property $$F(A\multimap
      B)=F(A)\multimap F(B),$$
  \item $G:C_1\to\Lambda(X,C_2)$ is a function such that
       for any $c\in C_1$ it holds that $\vdash_{\Sigma_2}G(c):F(\mathfrak{T}(c))$.
\end{itemize}

The map $G$ above extends inductively to a map
$$G:\Lambda(X,C_1)\to\Lambda(X,C_1)$$
  by
  $$G(x)=x,\mbox{ }x\in X,$$
   $$G(ts)=(G(t)G(s)),\quad G(\lambda x.t)=(\lambda x.G(t)).$$

For economy of notation, we write $\phi(A)$ for $F(A)$ when $A\in Tp(C_1)$, and we write $\phi(t)$ for $G(t)$ when $t\in\Lambda(X,C_1)$.

\bd
An {\it  abstract categorial grammar (ACG)} $G$ is a tuple $G=(\Sigma_{abstr},\Sigma_{obj}, \phi,S)$, where
\begin{itemize}
  \item $\Sigma_{abstr}$, $\Sigma_{obj}$, are linear signatures, respectively, the {\it abstract} and the {\it object} signature;
    \item $\phi:\Sigma_{abstr}\to \Sigma_{obj}$, the {\it lexicon}, is a map of signatures;
  \item $S$, {\it the standard type}, is an atomic type of $\Sigma_{abstr}$.
\end{itemize}
\ed

The {\it abstract} and the {\it object} language, respectively, $L_{abstr}(G) $ and $L_{obj}(G)$ generated by $G$ are defined as the sets of terms for which
$$
\vdash_{\Sigma_{abstr}}t:S,\mbox{ respectively, }\vdash_{\Sigma_{obj}}t:\phi(S).$$

In the context of formal languages, the object language generated by an ACG should be understood as a representation of  syntax.

\subsection{Representing as LLG}\label{translating ACG general}
Let $G=(\Sigma_{abstr},\Sigma_{obj}, \phi,S)$ be an ACG, and assume that an interpretation of the object signature $\Sigma_{obj}$ in the category of cowordisms is given.

This gives us immediately an interpretation of the abstract signature as well.

Indeed,  types and signature axioms of $\Sigma_{abstr}$, are mapped by the lexicon $\phi$ to  types and derivable typing judgements of $\Sigma_{obj}$, and these are mapped to  boundaries and cowordisms  by the given interpretation. Composing the two we get the desired interpretation of $\Sigma_{abstr}$.

Thus, to any type $A\in Tp(\Sigma_{abstr})$ we assign the boundary
$$[ A]=[\phi(A)],$$
and to any signature axiom $\vdash c:\mathfrak{T}(c)$ of $\Sigma_{abstr}$
we assign the cowordism
$$[c]=[\phi(c)]:{\bf 1}\to [\phi(\mathfrak{T}(c))]=[\mathfrak{T}(c)].$$

\bp\label{commutation with the lexicon}
Let $\sigma$ be a  typing judgement
$$
x_1:A_1,\ldots,x_n:A_n\vdash t:A
$$
derivable in $\Sigma_{abstr}$.
Then the interpretation $[\sigma]$ of $\sigma$  coincides with the interpretation of
the typing judgement
$$
x_1:\phi(A_1),\ldots,x_n:\phi(A_n)\vdash \phi(t):\phi(A)
$$
(which is derivable in $\Sigma_{obj}$).
\ep
\begin{proof}
  By induction on derivation.
\end{proof}

It follows that we have a cowordism signature $\Sigma'$ representing the abstract signature $\Sigma_{abstr}$ as described in Section \ref{representing linear signatures}.

We define the LLG $G'$ {\it representing $G$} by $G'=(\Sigma',S)$.

Lemma \ref{conservativity over ACG} and Proposition \ref{commutation with the lexicon} immediately imply that the LLG $G'$ is a translation of the ACG $G$.
\bl\label{ACG representation in general}
In notation as above, let $A_1,\ldots,A_n, A\in Tp(\Sigma_{abstr})$.

  A typing judgement  $$\frac{\sigma}{\vdash A_1^\bot,\ldots, A_n^\bot, A}$$ for some cowordism $\sigma$  is derivable in $G'$ iff there exists a term $t$ such that the typing judgement $$x:A_1,\ldots, A_n\vdash t:A$$ is derivable in the abstract signature $\Sigma_{abstr}$.

Furthermore, if all cowordisms of type $S$ generated by $G'$ are regular, then there are translations $$L_{abstr}(G)\to L(G),\quad L_{obj}(G)\to L(G)$$
of both the abstract and the object language of $G$ to the language of $G'$,
sending  a term $t\in L_{abstr}(G)$ (respectively, $t\in L_{obj}(G)$) to the interpretation $[t]$ of the corresponding typing judgement $\vdash t:S$ (respectively, $\vdash t:\phi(S)$).
$\Box$
\el
\smallskip

Whether the translation of grammars is {\it faithful} and the generated languages are {\it isomorphic} depends, of course, on the chosen interpretation of the object signature in the cowordism category. (In general, the interpretation can be {\it degenerate}, i.e. not injective, or it can map the object language to singular cowordisms.)

We now discuss the cases of {\it string} and {\it tree} ACG, relevant for analysis of syntax generation, where the translation is indeed a transparent isomorphism.

\section{Encoding string and tree ACG}
\subsection{Encoding string ACG}
\subsubsection{String signature}
Let $T$ be a finite alphabet.

The {\it string signature} $Str_T$ {\it over} $T$ is the linear signature with a single atomic type $O$, the alphabet $T$ as the set of constants and the  typing assignment
$$\mathfrak{T}(c)=O\multimap O\mbox{ }\forall c\in T.$$

We denote the type $O\multimap O$ as $str$.

If $t$ is a closed (i.e. not having free variables) term such that $\vdash_{Str_T}t:str$, we say that $t$ is a
{\it string term}.

Any word $a_1\ldots a_n$ in the alphabet $T$ can be represented as the string term
\be\label{string term}
\rho(a_1\ldots a_n)=(\lambda x.a_1(\ldots(a_n(x))\ldots)).
\ee

It is not hard to see that, if we identify $\beta\eta$-equivalent terms, the map $\rho$ has an inverse.
\bl\label{typeable string terms}
Any  string term $t$ is $\beta\eta$-equivalent to the term $\rho(w)$ for some $w\in T^*$.
\el
\begin{proof}
\begin{enumerate}[(i)]
\item There is no $\lambda$-term $t$ such that $\vdash_{Str_T}t:O$ (for example, because any derivable typing judgement has an even number of $O$ occurrences).
\item Using (i), we prove by induction on type inference that if $t$ is a $\beta$-normal term such that $\vdash_{Str_T}t:F$ for some type $F$, then $t$ is either a constant $t\in T$, or an abstraction, $t=(\lambda x.t')$ for some variable $x$ and term $t'$.
\item Using (ii), we prove by induction on type inference that for any derivable typing judgement $x:O\vdash_{Str_T}t:O$, where $t$ is a $\beta$-normal term, it holds that $t=c_1(\ldots(c_n(x))\ldots)$ for some constants $c_1,\ldots,c_n\in T$.
\end{enumerate}

Now without loss of generality we can assume that $t$ in the hypothesis of the lemma is $\beta$-normal.

If  $t$ is a constant then $t$ is   $\beta\eta$-equivalent to $\rho(t)$. Otherwise claim (ii) implies that the last rule in the derivation of the typing judgement $t:O\multimap O$ is ($\multimap\rm{E}$). Then the statement of the lemma follows from (iii).
\end{proof}

In the setting of Lemma \ref{typeable string terms}, we say that the word $w$ is {\it represented} by the string term $t$.
%
%

\subsubsection{Cowordism representation of the string signature}
Let us interpret  the atomic type $O$ as the one-point boundary $[O]$ with cardinality 1 and $[O]_l=\{1\}$.

By induction this gives us an interpretation  of all types in $Tp(Str_T)$ as boundaries.

We note that any regular cowordism $\sigma:O\to O$  is a graph consisting of a single edge labeled with some word $w\in T^*$. We denote such a cowordism as $[w]$ and call it the {\it cowordism representation} of $w$.

We interpret each signature axiom $\vdash c:O\multimap O$, where $c\in T$, as the corresponding regular cowordism $[c]:{\bf 1}\to [O]\multimap [O]$.

This gives us an interpretation of the signature $Str_T$ in the category of cowordisms over $T$.

We call this interpretation  the {\it cowordism representation of the string signature}.

Lemma \ref{typeable string terms} easily leads us to the following.
\bp\label{string typing judgements}
For any  string term $t$, the cowordism representation $[ t]$ of the corresponding typing judgement $\vdash t:str$ coincides with the cowordism representation $[w]$ of  the word $w\in T^*$ represented by $t$, $[t]=[w]$.
\ep
\begin{proof}
  By Lemma \ref{beta-eta}, the cowordism representation does not distinguish $\beta\eta$-equivalent terms, hence, by Lemma \ref{typeable string terms}, we can identify $t$ with $\rho(w)$.

  Induction on the length of $w$.

  For the empty word the statement is clear.

  Otherwise $w=w'c$, where $w'\in T^*$, $c\in T$.

  Then the term $t=\rho(w)$ is $\beta\eta$-equivalent to $t'=\lambda x.(\rho(w)c)$, and the typing judgement $\vdash_{Str_T}t':str$ is derived as follows
  $$
  \cfrac
    {
      \cfrac
        {
                \begin{array}{c}
                \\
                \vdash \rho(w'):O\multimap O
                \end{array}
                \quad
                \cfrac{
                \vdash c:O\multimap O\quad x:O\vdash x:O
                }
                {
                x:O\vdash  cx:O
                }
                ~(\multimap\rm{E})
        }
        {
        x:O\vdash (\rho(w'))(cx): c:O
        }~(\multimap\rm{E})
      }
      {
      \vdash t':O\multimap O
      }~(\multimap\rm{I}).
  $$

  Applying the induction hypothesis to $\rho(w')$  we immediately compute the interpretation of $\vdash t:str$ as $[w]$.
  Computation is shown in the picture below (with two elimination rules done at one step).
 $$
 \tikz[scale=.6]{

\node[right] at(1,.5) {$O^\bot$};
\node[right] at(1,0) {$O$};
\draw [fill] (1,0) circle [radius=0.05];
\draw[thick,<-](1,.5) to [out=180,in=90](0,.25) to  [out=-90,in=180] (1,0);
\node[above] at(.5,-.2){$w'$};

\node[right] at(1,1.5) {$O^\bot$};
\node[right] at(1,1) {$O$};
\draw [fill] (1,1) circle [radius=0.05];
\draw[thick,<-](1,1.5) to [out=180,in=90](0,1.25) to  [out=-90,in=180] (1,1);
\node[above] at(.5,.9){$c$};

\draw [fill] (1,2) circle [radius=0.05];
\draw[thick,<-](0,2)--(1,2);
      \node[right] at(1,2) {$O$};
 \node[left] at(0,2) {$O$};

\node[right] at(2.5,1) {$\Rightarrow$};

 \begin{scope}[shift={(4,0)}]

\node[right] at(1,0) {$O$};
\draw [fill] (1,.0) circle [radius=0.05];
\draw[thick,-](1,.5) to [out=180,in=90](0,.25) to  [out=-90,in=180] (1,0);
\node[above] at(.5,-0.2){$w'$};

\draw[thick,-](1,.5) to [out=0,in=-90](1.25,.75) to  [out=90,in=0] (1,1);

\draw[thick,-](1,1.5) to [out=180,in=90](0,1.25) to  [out=-90,in=180] (1,1);
\node[above] at(.5,.9){$c$};

\draw[thick,-](1,1.5) to [out=0,in=-90](1.25,1.75) to  [out=90,in=0] (1,2);

\draw[thick,<-](0,2)--(1,2);
      \node[left] at(0,2) {$O$};

\node[right] at(1.5,1) {$=$};

\begin{scope}[shift={(1,0)}]
\draw[thick,<-] (2,1)--(3.2,1);
\draw [fill] (3.2,1) circle [radius=0.05];
\node[left] at(2,1) {$O$};
\node[right] at(3.2,1) {$O$};
\node[above]at(2.6,.9){$w'c$};
\node[right] at(4,1) {$\Rightarrow$};
\end{scope}
 \begin{scope}[shift={(6,0)}]
\node[right] at(1,1.25) {$O^\bot$};
\node[right] at(1,.75) {$O$};
\draw [fill] (1,.75) circle [radius=0.05];
\draw[thick,<-](1,1.25) to [out=180,in=90](0,1) to  [out=-90,in=180] (1,.75);
\node[above] at(.5,.65){$w'c$};

 \end{scope}

 \end{scope}
 }.$$
  Since cowordism representation does not distinguish $\beta\eta$-equivalent terms, the statement is proven.
  \end{proof}

\subsection{Encoding string ACG}
A {\it string ACG} $G$ is an ACG whose object signature is the string signature over some alphabet $T$, $G=(\Sigma,Str_T,\phi,S)$.

The {\it string language} $L(G)$ generated by a string ACG $G$ is the set of words
$$L(G)=\{w\in T^*|~\rho(w)\in L_{obj}(G)\}.$$

Now let a string ACG $G$ be given.

We have an LLG $G'$ representing the ACG $G$ as described in Section \ref{translating ACG general}.

It follows from Lemma \ref{string typing judgements} that all cowordisms of type $str$ generated by $G'$ are regular. Hence by Lemma \ref{ACG representation in general} we have a translation of $L_{abstr}$ and $L_{obj}(G)$ to $L(G)$. Lemma \ref{string typing judgements} guarantees that this translation is one-to-one on the level of strings, i.e.  the string language $L(G)$ generated by $G$ coincides with the language $L(G')$ generated by $G'$.

We summarise.
\bt\label{encoding ACG}
If a   language is generated by a string ACG then it is also generated by an LLG. $\Box$
\et
\smallskip

It seems an interesting question whether the converse is true or not. We would expect that the two types of grammars generate the same class of languages.

\begin{remark}
 Since MCFG embed into string ACG \cite{deGrootePogodalla}, Theorem \ref{encoding ACG} on encoding  ACG in LLG grammars implies that MCFG embed into LLG. However it does not imply  the converse statement (Theorem  \ref{tens. free lex is MCFG}, that any $\otimes$-free lexicon gives rise to an MCFG).

 On the other hand it is not hard to see that Theorem \ref{tens. free lex is MCFG} together with Theorem \ref{encoding ACG} do imply the known result \cite{Salvati} that any second order string ACG generates a multiple context-free language. Thus we gave another, possibly  more ``geometric''  proof of this result.
 \end{remark}

\subsection{Tree ACG}
Linguists often represent natural language in terms of {\it labeled trees} rather than just strings and consider formal {\it tree languages} and corresponding formal grammars, such as {\it tree adjoining grammars} \cite{JoshiTAG}. We now discuss how to represent tree languages and tree ACG in the setting of cowordisms.

\subsubsection{Tree languages}
In what follows, a tree means a {\it rooted planar tree}.

We assume that we are given an alphabet $N$ of {\it node labels}, and we will consider {\it labeled trees} whose nodes are labeled with elements of $N$. We denote the set of all such trees  as  $Tr(N)$.
\begin{remark}
The alphabet of node labels is often subdivided further into alphabets of terminal and nonterminal symbols, but this is a purely cosmetic detail, and constructions that we discuss below can be easily adapted to this setting.
\end{remark}

It is convenient to introduce ``elementary building blocks'' for constructing elements of $Tr(N)$.

Let us say that an {\it elementary  tree} $A_k$, where $A\in N$, is a planar rooted  tree consisting of the root  labeled with $A$ and its $k$ children, which are unlabeled.

Each elementary tree $A_k$ can be identified with a $k$-ary functional symbol.
Then  trees over $N$ can be identified with terms constructed from these functional symbols.

Namely, the tree with a single vertex labeled with $A\in N$ corresponds to the constant $A_0$, and the tree $\alpha\in N$
 obtained by attaching trees $\alpha_1,\ldots,\alpha_k$   to a  root labeled with $A\in N$ corresponds to the term $A_k(\alpha_1,\ldots,\alpha_k)$.

Now, given a tree $\alpha$, we say that
the {\it branching} of a node $v$ of $\alpha$ is the number of children of $v$.
The {\it branching} of $\alpha$ is the maximal branching of its nodes.

We denote the set of all labeled trees over  $N$ with branching less or equal to $n$ as  $Tr(N,n)$.

A {\it tree language} over alphabet $N$  is a set of planar rooted trees whose nodes are labeled with elements of $N$.

We say that a tree language $L$ has {\it bounded branching} if there is $n$ such that the branching of all trees in $L$ have branching less or equal to $n$, $L\subseteq Tr(N,n)$. In this case we say that $n$ is the {\it maximal branching} of $L$, if all trees in $L$ have branching less or equal to $n$ and there is $\alpha\in L$ whose branching is $n$.

\subsubsection{Encoding tree languages}\label{encoding tree languages}
The set $Tr(N)$ of labeled trees
can be conveniently encoded into  words over a certain alphabet.
If, furthermore, we restrict to the set $Tr(N,n)$ of trees with bounded branching, then the alphabet for the encoding is finite.

Given an alphabet $N$ of node labels and a maximal branching $n$, let us put $$T(N,n)=\{A_i^k|~A\in N,0\leq i\leq k,0\leq k\leq n\}.$$

We encode each element $\alpha\in Tr(N,n)$ as a word $[\alpha]\in(T(N,n))^*$
 inductively as follows (identifying labeled trees with functional terms):
$$[A_0]=A_0^0,\quad[A_k(\alpha_1,\ldots,\alpha_k)]=A^0_k[\alpha_1]A^1_k[\alpha_2]\ldots A^{k-1}_k[\alpha_k]A^k_k.$$

Given a tree language $L$ over $N$ with branching $n$, we denote as $[L]$ the image of $L$ in $(T(N,n))^*$ under the above encoding.
Note that the encoding map from  $L$ to $[L]$ is a bijection.

Apart from representing ``finished'' trees we also need to represent trees  ``under construction'' and operations on them. For that purpose we will encode elementary trees as cowordisms over $T(N,n)$.

We introduce the boundary $\mathcal{T}$ defined by
\be\label{cowordism tree type}
|\mathcal{T}|=2,\quad \mathcal{T}_l=\{1\}
\ee
and map each elementary tree $A_k$ to the {\it elementary tree cowordism}
$$[A_k]:\underbrace{\mathcal{T}\otimes\ldots\otimes\mathcal{T}}_{k\mbox{ times}}\to \mathcal{T}$$
over $T(N,n)$, whose edges are labeled with the words  $A_k^0,\ldots,A_k^k$, defined by the following picture (in the vertical representation).
$$
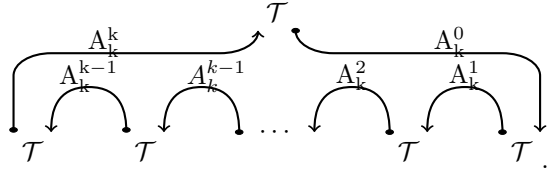
\tikz[yscale=.6]{

\draw [fill] (-3.5,0.05) circle [radius=0.05];
\draw [fill] (-2,0.05) circle [radius=0.05];
\draw [fill] (-.5,0) circle [radius=0.05];
\draw [fill] (1.5,0) circle [radius=0.05];
\draw [fill] (3,0) circle [radius=0.05];
\draw [fill] (.25,2.25) circle [radius=0.05];

\draw[thick,->](-3.5,0)--(-3.5,1.25)
to[out=90, in =180]
(-3,1.75)--(-.75,1.75)
to[out=0, in =-90]
(-.25,2.25);

\draw[thick,<-](-3,0)
to[out=90, in =180](-2.5,1)
to[out=0, in =90](-2,0);
\node at (-2.5,1.25){$\rm{A^{k-1}_k}$};

\node at (-.8,1.25){$A_k^{k-1}$};

\draw[thick,<-](-1.5,0)
to[out=90, in =180](-1,1)
to[out=0, in =90](-.5,0);

\node at (0,0){$\cdots$};

\draw[thick,<-](.5,0)
to[out=90, in =180](1,1)
to[out=0, in =90](1.5,0);
\node at (1.,1.25){$\rm{A_k^{2}}$};

\draw[thick,<-](2,0)
to[out=90, in =180](2.5,1)
to[out=0, in =90](3,-0);
\node at (2.5,1.25){$\rm{A_k^1}$};

\draw[thick,<-](3.5,0)--(3.5,1.25)
to[out=90, in =0]
(3,1.75)--(.75,1.75)
to[out=180, in =-90]
(.25,2.25);

\node at (2.3,2){$\rm{A_k^0}$};
\node at (-2.3,2){$\rm{A_k^k}$};

\node[above] at (0,2.25){$\mathcal{T}$};
\node[below] at (3.25,0){$\mathcal{T}$};
\node[below] at (1.75,0){$\mathcal{T}$};
\node[below] at (-1.75,0){$\mathcal{T}$};
\node[below] at (-3.25,0){$\mathcal{T}$};
}.
$$

The case $k=0$ corresponds to the tree $A_0$ with a single vertex labeled with $A$, and the cowordism  $$[A_0]:{\bf 1}\to \mathcal{T}$$ consists of a single edge labeled with $A_0^0$.

For a ``finished tree'' $\alpha\in Tr(N,n)$ we identify its encoding word $[\alpha]\in (T(N,n))^*$ with the corresponding regular cowordism from ${\bf 1}$ to $\mathcal{T}$ whose only edge is labeled with $[\alpha]$.

It is easy to see then that this encoding is a sound representation of elementary trees considered as functional symbols, i.e.
$$[A_k(\alpha_1,\ldots,\alpha_k)]=[A_k]\circ([\alpha_1]\otimes\ldots\otimes[\alpha_k]).$$

\subsection{Tree signature}
(Tree signatures and tree ACG were defined in \cite{deGrooteTAG} without special names.)

Labeled trees with bounded branching, being, essentially, functional terms, can be represented as linear $\lambda$-terms in the obvious way.

Let an alphabet $N$ and a maximal branching $n$ be given.

The {\it tree signature} $Tree_{N,n}$ has unique atomic type $\mathcal{T}$ (notational clash with (\ref{cowordism tree type}) is intentional), the set of constants $C=\{A_k|~A\in N, k\leq n\}$ (another intentional abuse of notation!)
and signature axioms
\be\label{tree signature axioms}
\vdash A_k:\underbrace{\mathcal{T}\multimap\ldots\multimap\mathcal{T}}_{k\mbox{ times}}\multimap \mathcal{T},
\ee
where $A\in N$, $k\leq n$.

If $t$ is a closed (i.e. not having free variables) term such that $\vdash_{Tree_{N,n}}t:\tau$, we say that $t$ is a
{\it tree term}.

A labeled tree $\alpha\in Tr(N,n)$ is represented as a tree term $\rho(\alpha)$ defined by induction:
$$\rho(A_0)=A_0\quad
  \rho(A_k(\alpha_1,\ldots,\alpha_k))=A_k\rho(\alpha_1)\ldots \rho(\alpha_k),$$
  where notational convention (\ref{iterated application}) is used.

It  is not hard to show that any tree term  represents an element of $Tr(N,n)$ up to a $\beta\eta$-equivalence.

Let us introduce some terminology.

We define {\it open tree types} by the following induction:
\begin{itemize}
\item $\mathcal{T}$ is the open tree type with $0$ {\it hanging vertices};
\item if $F$ is the open tree type with $m$ { hanging vertices}, then $\mathcal{T}\multimap F$ is the open tree type with $m+1$ {\it hanging vertices}.
\end{itemize}

And {\it open tree terms} are defined by the following:
\begin{itemize}
\item any $\beta$-normal term built from constants and variables using only application is an open tree term with $0$ {\it bound leaves};
\item a term $\lambda x.t$, where $t$ is an open tree term with $k$ {bound leaves},  is an open tree term with $k+1$ {\it bound leaves}.
\end{itemize}

\bl\label{typeable tree terms}
Any tree term $t$,  is $\beta\eta$-equivalent to the term $\rho(\alpha)$ for some $\alpha\in Tr(N,n)$.
\el
\begin{proof}
Without loss of generality, the term $t$ is $\beta$-normal.

 Then the statement follows from the following more general statement, which is proven by induction on derivation.

Assume that
$$x_1:\mathcal{T},\ldots, x_n:\mathcal{T}\vdash_{Tree_{N,n}} t:F,$$ where $t$ is $\beta$-normal.
\begin{enumerate}[(i)]
\item If $F$ is an open tree type with $m$ hanging vertices, then there exists $k\leq m$ such that $t$ is an open tree term with $k$ bound leaves.
\item If $t$ is not of the form $t=\lambda y.t'$, then $F$ is an open tree type and $t$ is an open tree term with no bound vertices.
\end{enumerate}
\end{proof}

In the setting of Lemma {\ref{tree typing judgements}} we say that $\alpha$ is the tree {\it represented} by the term $t$.

\subsubsection{Cowordism representation of the tree signature}
An interpretation of the tree signature, essentially, has already been constructed in Section \ref{encoding tree languages}.

Let us interpret  the atomic type  $\mathcal{T}$ as the corresponding boundary $\mathcal{T}$ defined in (\ref{cowordism tree type}).

By induction this gives us an interpretation  of all types in $Tp(\mathcal{T})$ as boundaries.

We interpret each signature axiom of  form (\ref{tree signature axioms}) as the elementary tree cowordism $[A_i]$.

This gives us an interpretation of the signature $Tree_{N,n}$ in the category of cowordisms over $T(N,n)$.

We call this interpretation  the {\it cowordism representation} of tree signature.

We have the following proposition, identical in form to Proposition \ref{string typing judgements}.
\bp\label{tree typing judgements}
For any tree term $t$, the cowordism representation $[t]$  of the corresponding typing judgement
$\vdash_{Tree_{N,n}} t:\mathcal{T}$ coincides with the cowordism representation $[\alpha]$ of the tree $\alpha\in Tr(N,n)$ represented by $t$, $[t]=[\alpha]$.
\ep
\begin{proof}
  Similar to Proposition \ref{string typing judgements}.

  Without loss of generality, the term $t$ can be assumed $\beta$-normal, then, by Lemma \ref{typeable tree terms}, $t$ is $\beta\eta$-equivalent to $\rho(\alpha)$ for some $\alpha\in Tr(N,n)$, we use  induction on $\alpha$.
\end{proof}

\subsection{Encoding tree ACG}
 Tree abstract categorial grammars (tree ACG) are defined similarly to string ACG.

A {\it tree ACG} $G$ is an ACG whose object signature is the tree signature over some alphabet  $N$ of node labels with maximal branching $n$, $$G=(\Sigma,{Tree_{N,n}},\phi,S).$$

The {\it tree language} $L(G)$ generated by a tree ACG $G$ is the set of trees
$$L(G)=\{\alpha\in Tr(N,n)|~\rho(w)\in L_{obj}(G)\}.$$

Tree ACG are encoded in LLG in the same way as string ACG.

Let  a tree ACG $G$ be given.

The cowordism representation of the tree signature is an interpretation of the object signature in the category of cowordisms. Hence
we have an LLG $G'$ representing the ACG $G$ as described in Section \ref{translating ACG general}.

Lemma \ref{ACG representation in general} and Lemma \ref{tree typing judgements} immediately imply that the encoding $[L(G)]$ of the tree language $L(G)$ generated by $G$ coincides with the language $L(G')$ generated by $G'$.

We summarise.
\bt
If a tree  language $L$ is generated by a tree ACG then its encoding $[L]$ is  generated by an LLG. $\Box$
\et

\section{Encoding backpack problem}
It is known that ACG, in general, can generate NP-complete languages \cite{YoshinakaKanazawa} (and for unlexicalized ACG, even decidability is an open question).

It seems natural to expect  that LLG can generate NP-complete languages as well. In this last section we show how an LLG can generate solutions of the backpack problem.
Our purpose here is mainly  illustrative. We try to convince the reader that the geometric language of cowordisms is indeed intuitive and convenient for analysing language generation.

We consider backpack problem in the form of the {\it subset sum problem}.

\bd {\bf Subset sum problem (SSP)}: Given a finite sequence $s$ of integers, determine if there is a subsequence $s'\subseteq s$ such that $\sum\limits_{z\in s'}z=0$.
\ed

SSP is known to be NP-complete, see \cite{Martello}.

We now define a language representing solutions of SSP.

We  represent integers as words in the alphabet $\{+,-\}$, we call them {\it numerals}. An integer $z$ is represented (non-uniquely) as a word for which the difference of $+$ and $-$ occurrences equals $z$.

We say that a numeral is {\it irreducible}, if it consists only of pluses or only of minuses.

We represent finite sequences of integers as words in the alphabet $T=\{+,-,\bullet\}$, with $\bullet$ interpreted as a separation sign. Thus a word in this alphabet should be read as a list of numerals separated by bullets.

When all numerals in the list are irreducible, we say that the list is irreducible. Note that any sequence of integers has unique representation as an irreducible list.


Let us construct an {\bf LLG} that generates solutions of SSP.

We will use two positive literals $H,S$ and interpret each of them as the boundary $X$ of cardinality 2 with $X_l=\{1\}$.


First we construct a system of cowordisms which generates lists of numerals representing sequences that sum to zero.

We define  four cowordisms
$$cons:S\otimes S\to S,\quad  push:H\otimes H\to H\otimes H,\quad close:{\bf 1}\to H$$
in the graphical language as follows.
 $$
 
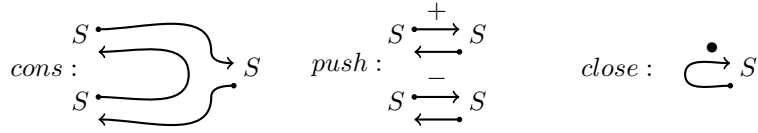
\tikz[scale=.6]
        {

         \draw[thick,->](0,0) to  [out=0,in=90] (2.5,-0.5) to  [out=-90,in=180] (3,-0.75);
        \draw [fill] (0,0) circle [radius=0.05];

        \draw[thick,<-](0,-0.5) to  [out=0,in=90] (2,-1) to  [out=-90,in=0] (0,-1.5);
        \draw [fill] (-0,-1.5) circle [radius=0.05];
          \draw[thick,<-](0,-2.)to  [out=0,in=-90](2.5,-1.5) to  [out=90,in=-180] (3,-1.25);
          \draw [fill] (3,-1.25) circle [radius=0.05];
          \node[above left] at (0,-0.5) {$S$};
        \node[above left] at (0,-2) {$S$};
        \node[above right] at (3,-1.25) {$S$};
        \node[above left] at(-.25,-1.25) {$cons:$};
        \begin{scope}[shift={(7,0)}]

        \draw [fill] (0,0) circle [radius=0.05];
        \draw[thick,->](0,0) -- (1,0);

        \draw [fill] (1,-0.5) circle [radius=0.05];
        \draw[thick,<-](0,-0.5) -- (1,-0.5);
        \node[above right] at (1,-.5) {$S$};
        \node[above left] at (0,-.5) {$S$};
        \node[above] at (0.5,0) {$+$};

        \draw [fill] (0,-1.5) circle [radius=0.05];
        \draw[thick,->](0,-1.5) -- (1,-1.5);
        \draw [fill] (1,-2) circle [radius=0.05];
        \draw[thick,<-](0,-2.) -- (1,-2.);
        \node[above right] at (1,-2) {$S$};
        \node[above left] at (0,-2) {$S$};
        \node[above] at (0.5,-1.5) {$-$};

        \node[above left] at(-.5,-1.25) {$push:$};
%
        \begin{scope}[shift={(6,-.75)}]
          \node[above] at (0.6,0) {$\bullet$};
        \draw[thick,<-](1,0) to  [out=-180,in=90] (0,-0.25) to  [out=-90,in=-180] (1,-.5);
        \draw [fill] (1,-0.5) circle [radius=0.05];

        \node[above right] at (1,-.5) {$S$};
        \node[above left] at(-.5,-.5) {$close:$};
        \end{scope}
%
        \end{scope}
         }\\
         $$

The cowordism $cons$, by iterated compositions with itself, generates lists with arbitrary many empty slots. Then $push$ fill the slots (always in pairs), and  $close$ closes a slot and puts a separation sign.

It is easy to see that all cowordisms from ${\bf 1}$ to $S$ generated by the above system together with symmetry transformations represent sequences of integers summing to zero, and vice versa, for any sequence summing to zero, its irreducible list representation is generated by the above.

Now, in order to generate solutions of SSP we need some extra, ``deceptive'' slots, which  contain elements not summing to zero. These slots  will be represented by the boundary $H$.

%

We define  cowordisms
$$open_H:H\otimes S\to S,\quad close_H:{\bf 1}\to H,$$
$$\quad push_+:H\to H,\quad push_-:H\to H$$
as follows.
%
 $$
 
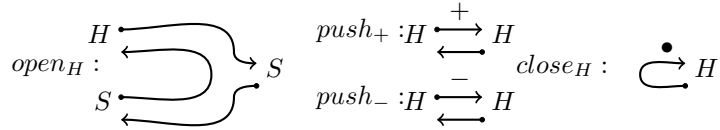
\tikz[scale=.6]{
\draw [fill] (0,0) circle [radius=0.05];
 \draw[thick,->](0,0) to  [out=0,in=90] (2.5,-0.5) to  [out=-90,in=180] (3,-0.75);

\draw [fill] (0,-1.5) circle [radius=0.05];
\draw[thick,<-](0,-0.5) to  [out=0,in=90] (2,-1) to  [out=-90,in=0] (0,-1.5);
\draw [fill] (3,-1.25) circle [radius=0.05];
  \draw[thick,<-](0,-2.)to  [out=0,in=-90](2.5,-1.5) to  [out=90,in=-180] (3,-1.25);
  \node[above left] at (0,-0.5) {$H$};
\node[above left] at (0,-2) {$S$};
\node[above right] at (3,-1.25) {$S$};
\node[above left] at(-.25,-1.25) {$open_H:$};
\begin{scope}[shift={(7,0)}]
\begin{scope}[shift={(4.5,-.75)}]
      \node[above] at (0.6,0) {$\bullet$};
    \draw[thick,<-](1,0) to  [out=-180,in=90] (0,-0.25) to  [out=-90,in=-180] (1,-.5);
    \draw [fill] (1,-.5) circle [radius=0.05];
    \node[above right] at (1,-.5) {$H$};
    \node[above left] at(-.5,-.5) {$close_H:$};
\end{scope}
  \draw [fill] (0,0) circle [radius=0.05];
\draw[thick,->](0,0) -- (1,0);
\draw [fill] (1,-0.5) circle [radius=0.05];
\draw[thick,<-](0,-0.5) -- (1,-0.5);
\node[above right] at (1,-.5) {$H$};
\node[above left] at (0,-.5) {$H$};
\node[above] at (0.5,0) {$+$};
\node[above left] at(-.5,-.5) {$push_+:$};
\begin{scope}[shift={(0,-1.5)}]
  \draw [fill] (0,0) circle [radius=0.05];
\draw[thick,->](0,0) -- (1,0);
\draw [fill] (1,-0.5) circle [radius=0.05];
\draw[thick,<-](0,-0.5) -- (1,-0.5);
\node[above right] at (1,-.5) {$H$};
\node[above left] at (0,-.5) {$H$};

\node[above] at (0.5,0) {$-$};
\node[above left] at(-.5,-.5) {$push_-:$};
\end{scope}
\end{scope}
 }$$
%
%
%
%
%
%

%
%
%
%
The cowordism $open_H$ adds deceptive slots to the list, $push_-$ and $push_+$ fill them with arbitrary numerals, and $close_H$ closes them.

%
%
%
%
%
%

Let us denote the set of cowordisms from ${\bf 1}$ to $S$ generated by the above system and symmetry as $L_0$.

It is easy to see that  $L_0$ membership problem is, essentially,  SSP.
In particular,  a sequence $s$ of integers is a solution of SSP iff the corresponding irreducible list is in $L_0$.
It follows that $L_0$ is NP-hard.

It is also easy to show that $L_0$ membership problem is itself in NP, hence $L_0$ is, in fact, NP-complete.

We define an LLG $G$ by a lexicon consisting of names of the above cowordisms.
Then it is easy to see that $G$  generates $L_0$. This can be established  by reasoning similar to that in Lemma \ref{conservativity over ACG}. We omit a proof because of space limitation.

(This example shows also that the  generative power of LLG is strictly greater than MCFG, because multiple context-free languages are effectively decidable.)


\begin{thebibliography}{45}
\bibitem{AbramskyFreeTraced}
S. Abramsky,
``Abstract scalars, loops, free traced and strongly compact closed categories'',
Proceedings of CALCO 2005, Springer Lecture Notes in Computer Science, 3629, pp. 1-31, 2005.
\bibitem{AbramskyCoecke}  S. Abramsky, B. Coecke. ``Categorical quantum mechanics'', in Handbook
of quantum logic and quantum structures: quantum logic, pp. 261-324,
2008.
\bibitem{Baez}  J. C. Baez and J. Dolan, ``Higher-dimensional Algebra and Topological Quantum Field Theory'', J.Math.Phys. 36,  pp. 6073-6105, 1995.
\bibitem{Barendregt} H.P. Barendregt,  ``The Lambda Calculus — Its Syntax and Semantics''. Studies in Logic and the Foundations of Mathematics. 103. Amsterdam: North-Holland, 1985.
\bibitem{Barr}  M. Barr. ``$*$-Autonomous Categories'', Lecture Notes in Mathematics 752, Springer,
1979.
\bibitem{HylandDePaiva} N. Benton, G. Bierman, J. Hyland, V. de Paiva,
"Term assignment for Intuitionistic Linear Logic",
Report 262, Computer Laboratory, University of Cambridge, 1992.
\bibitem{CocketSeely} R. Cockett, R. Seely, ``Weakly Distributive Categories'', Journal of Pure and Applied Algebra, 114(2), pp 133-173, 1997.
\bibitem{CoeckeSadrzadehClark} B. Coecke, M. Sadrzadeh, S. Clark, ``Mathematical Foundations for a Compositional Distributional Model of Meaning''. Lambek Festschirft, special issue of Linguistic Analysis, 2010.
\bibitem{Coecke_Lmabek_vs_Lambek}    B. Coecke, M. Sadrzadeh, M. Sadrzadeh, ``Lambek vs. Lambek: Functorial Vector Space Semantics and String Diagrams for Lambek Calculus'',
     Annals of Pure and Applied Logic, 164(11), 1079-1100, 2013.
\bibitem{Dalrymple}  M. Dalrymple., J. Lamping, F. Pereira, F., V. Saraswat,  “Linear logic for meaning assembly,”
in Proceedings of CLNLP, Edinburgh, South Queensferry: ELSNET, 1995.
\bibitem{deGroote}    P. de Groote,
  ``Towards Abstract Categorial Grammars'',
in  Proceedings of the 39th Annual Meeting on Association for Computational Linguistics, ACL '01, pp.148-155, 2001.
\bibitem{deGrooteTAG}
 P. de Groote, ``Tree-Adjoining Grammars as Abstract Categorial Grammars'', In:
TAG+6, Proceedings of the sixth International Workshop on Tree Adjoining Grammars and Related Frameworks, pp. 145-150, 2001.
\bibitem{deGrootePogodalla}   P. de Groote, S. Pogodalla,  ``On the expressive power of
abstract categorial grammars: Representing context-free formalisms''. Journal of
Logic, Language and Information 13(4):421-438, 2004.
\bibitem{Girard} Jean-Yves Girard, ``Linear logic'', Theoretical Computer Science, 50:1-102, 1987.
 \bibitem{Girard2} Jean-Yves Girard, ``Linear logic: its syntax and semantics'', in J.-Y.Girard, Y.Lafont and L.Regnier, eds. Advances in Linear Logic, 1-42, Cambridge University Press, 1995, Proc. of the Workshop on Linear Logic, Ithaca, New York, June, 1993.
\bibitem{GirardProofsAndTypes}
J.-Y. Girard, P. Taylor, and Y. Lafont,  ``Proofs and Types'', Cambridge University Press, New York, NY, USA, 1989.
\bibitem{HS} M. Hyland and A. Schalk, ``Glueing and
Orthogonality for Models of Linear Logic'', Theoretical
Computer Science 294, pp. 183--231, 2003.
 \bibitem{JoshiTAG} A.K. Joshi, L.S. Levy, M. Takahashi,
``Tree adjunct grammars'',
Journal of Computer Systems Science, 10 (1), pp. 136-163,  1975.
\bibitem{Kanazawa} M. Kanazawa, ``The Pumping Lemma for Well-Nested Multiple Context-Free Languages'', in Developments in Language Theory, 13th International Conference, {DLT}
               2009, Stuttgart, Germany, June 30 - July 3, 2009. Proceedings, Lecture Notes in Computer Science 5583, pp. 312-325, 2009.
\bibitem{KellyLaplaza}  G.M. Kelly, M.L. Laplaza,   ``Coherence for compact closed categories''. Journal of Pure and Applied Algebra. 19: 193-213, 1980.
\bibitem{KubotaLevine} Y. Kubota, R. Levine, ``Gapping as like-category coordination'', in
D. B\'echet \& A. Dikovsky, eds, `Logical Aspects of Computational Linguistics',
Vol. 7351 of Lecture Notes in Computer Science, Springer, Nantes,
pp. 135-150, 2012.
\bibitem{Lambek} J. Lambek, Joachim, ``The mathematics of sentence structure'', Amer. Math. Monthly, 65: 154-170, 1958.
\bibitem{Lambek_pregroups} J. Lambek. ``Type grammar revisited'', Logical Aspects of Computational Linguistics, 1582, 1999.
\bibitem{MacLane} S. Mac Lane, ``Categories for the working mathematician'', Springer-Verlag, 1971.
\bibitem{Martello} S. Martello, P. Toth,  ``4 Subset-sum problem''. Knapsack problems: Algorithms and computer interpretations'', Wiley-Interscience. pp. 105-136, 1990.
\bibitem{Mellies_categorical_semantics} P.-A. Melli\'es, ``Categorical semantics of linear logic'', in: Interactive Models of Computation and Program Behaviour, Panoramas et Synth\`eses 27, Soci\'et\'e Math\'ematique de France 1-196, 2009.
\bibitem{PollardMichalicek}  V. Mihali\v{c}ek, C. Pollard, ``Distinguishing phenogrammar from
tectogrammar simplifies the analysis of interrogatives.'' In  {Proceedings of the 15th and 16th International Conference on Formal Grammar} (FG'10/FG'11),
 130-145. Springer  2012.
\bibitem{Moortgat}  M.  Moortgat,  ``Categorial type logics'', in Johan van Benthem
and Alice ter Meulen, eds., Handbook of Logic and Language,
chapter 2, pp. 93-178, Elsevier, MIT Press, 1997.
\bibitem{Moot_comparing} R. Moot, ``Comparing and evaluating extended
Lambek calculi'',.
In Kubota, Y. and Levine, R., editors,
Proceedings for ESSLLI 2015
Workshop ‘Empirical Advances in
Categorial Grammar’, University of Tsukuba and Ohio State, pp. 108-131, 2015.
\bibitem{Moot_inadequacy}
R. Moot,  ``Hybrid type-logical grammars, 1rst-order linear logic and the
descriptive inadequacy of lambda grammars'', https://hal.archives-ouvertes.fr/
hal-00996724,  2014.
\bibitem{MootRetore} R. Moot, C. Retor\'e,  "The Logic of Categorial Grammars: A
Deductive Account of Natural Language Syntax and Semantics", Lecture Notes in
Artificial Intelligence, Springer,  2012.
\bibitem{Morrill_Displacement} G. Morrill, O. Valent\'{\i}n, M. Fadda, ``The displacement calculus'', Journal of Logic, Language and Information 20(1), 1-48, 2011.
\bibitem{Muskens} R. Muskens, ``Separating syntax and combinatorics in categorial grammar'',
Research on Language and Computation 5(3), 267-285,  2007.
\bibitem{Pentus} M. Pentus, ``Lambek Grammars Are Context Free'', in Proceedings of the Eighth Annual Symposium on Logic in Computer Science
               {(LICS} '93), Montreal, Canada, June 19-23, 1993, pp. 429--433, 1993.
\bibitem{Salvati} S. Salvati, ``Encoding second order string ACG with deterministic tree
walking transducers''. In S. Wintner, ed., Proceedings of FG 2006: The 11th
conference on Formal Grammar, FG Online Proceedings, pages 143-156. CSLI
Publications, 2017.
\bibitem{Seely}  R.A.G. Seely, ``Linear logic, $*$-autonomous categories and cofree coalgebras'', in: J.Gray and A.Scedrov (editors), Categories in Computer Science and Logic, Contemporary Mathematics 92, 371-382, Amer. Math. Soc., 1989.
\bibitem{Seki} H. Seki,  T. Matsumura, M. Fujii, and T. Kasami, ``On
multiple context-free grammars'', Theoretical Computer Science 88(2):191-229, 1991.
\bibitem{Selinger} P. Selinger,  ``A survey of graphical languages for monoidal categories'', in: B. Coecke (ed.), New Structures
for Physics,  275-337,  Springer-Verlag, 2011.
\bibitem{Schelinx} H. Schellinx, ``Some Syntactical Observations on Linear Logic'', Journal
of Logic and Computation 1, 4, 537-559, 1991.
\bibitem{Tan}  A.M. Tan, ``Full completeness for Models of Linear Logic', PhD. Thesis, Cambridge,
1997.
\bibitem{TroelstraSchwichtenberg}
A. S. Troelstra and H. Schwichtenberg,  ``Basic Proof Theory'', Cambridge University Press, New York, NY, USA, 1996.
\bibitem{YoshinakaKanazawa} R. Yoshinaka, M. Kanazawa,  ``The Complexity and Generative Capacity of Lexicalized Abstract Categorial Grammars''. In: Blache P., Stabler E., Busquets J., Moot R. (eds) Logical Aspects of Computational Linguistics. LACL. Lecture Notes in Computer Science, vol 3492. Springer, Berlin, Heidelberg, 2005.
\end{thebibliography}
 \end{document}